\documentclass[preprint,12pt,3p]{elsarticle}

\usepackage{subfigure,graphicx}
\usepackage{algorithm}
\usepackage{latexsym,amsmath,amssymb,graphicx}
\usepackage{algpseudocode}
\usepackage{multirow}
\usepackage{booktabs}
\usepackage{amsmath,amssymb}
\usepackage{epstopdf}
\usepackage{url}
\usepackage{multirow}
\usepackage{helvet}
\usepackage{courier}
\usepackage{booktabs}
\usepackage{latexsym,amsmath,amssymb,graphicx}
\usepackage{url}
\usepackage{color}
\usepackage{epstopdf}
\usepackage{bm}
\usepackage{geometry} 

\usepackage{color}

\newcommand{\pp}{\mathbf{p}}
\newcommand{\qq}{\mathbf{q}}

\newcommand{\llambda}{\bm{\lambda}}

\date{}
\begin{document}
	\begin{frontmatter}

\title{\textbf{A New Operator Splitting Method for  Euler's Elastica Model}}
\tnotetext[label0]{This is only an example}

% L.-J Deng(UESTC)
%\author{Liang-Jian Deng\thanks{School of Mathematical Sciences, University of Electronic Science and Technology of China, Chengdu, Sichuan, 611731, China.
%		(\email{liangjian.deng@uestc.edu.cn})}
%	\and Roland Glowinski\thanks{Department of Mathematics, University of Houston, Houston, Texas, 77204, USA; Department of Mathematics, Hong Kong Baptist University, Kowloon Tong, Hong Kong (\email{roland@math.uh.edu}).}
%	\and Xue-Cheng Tai\thanks{Department of Mathematics, Hong Kong Baptist University, Kowloon Tong, Hong Kong (\email{xuechengtai@hkbu.edu.hk}).}}

\author[label1]{Liang-Jian Deng}
\address[label1]{School of Mathematical Sciences, University of Electronic Science and Technology of China, Chengdu, Sichuan, 611731, China.}
%\address[label2]{Address Two\fnref{label4}}
%
%\cortext[cor1]{I am corresponding author}
%\fntext[label3]{I also want to inform about\ldots}
%\fntext[label4]{Small city}

\ead{liangjian.deng@uestc.edu.cn}

\author[label5]{Roland Glowinski}
\address[label5]{Department of Mathematics, University of Houston, Houston, Texas, 77204, USA; Department of Mathematics, Hong Kong Baptist University, Kowloon Tong, Hong Kong}
\ead{roland@math.uh.edu}

\author[label6]{Xue-Cheng Tai}
\address[label6]{Department of Mathematics, Hong Kong Baptist University, Kowloon Tong, Hong Kong}
\ead{xuechengtai@hkbu.edu.hk}

\begin{abstract}
	Euler's elastica  model has a wide range of applications in Image Processing and Computer Vision. However, the non-convexity, the non-smoothness and the nonlinearity of the associated energy functional make its minimization a challenging task, further complicated by the presence of high order derivatives in the model. In this article we propose a new operator-splitting algorithm to minimize the Euler elastica functional. This algorithm is obtained by applying an operator-splitting based time discretization scheme to an initial value problem (dynamical flow) associated with the optimality system (a system of multivalued equations). The sub-problems associated with the three fractional steps of the splitting scheme have either closed form solutions or can be handled by fast dedicated solvers. Compared with earlier approaches relying on ADMM (Alternating Direction Method of Multipliers), the new method has, essentially, only the time discretization step as free parameter to choose, resulting in a very robust and stable algorithm. The simplicity of the sub-problems and its modularity make this algorithm quite efficient. Applications to the numerical solution of smoothing test problems demonstrate the efficiency and robustness of the proposed methodology.
\end{abstract}


\end{frontmatter}

% REQUIRED
%\begin{AMS}
%	%  68Q25, 68R10, 68U05
%\end{AMS}

\section{Introduction}\label{sec:intro}
The Euler elastica energy of a planar curve is defined as the following functional,
\begin{equation}\label{eq:elastica_term}
E(\mathcal{C}) = \int_{\mathcal{C}}\left(a + b \kappa^2 \right)ds,
\end{equation}
where $\mathcal{C}$ represents the planar curve whose curvature is $\kappa$, $s$ stands for arc length, $a$ and $b$ are two positive parameters. Especially, if $a = 0$, $E(\mathcal{C})$ measures the twisting energy of the curve that is related to the curvature; if $b = 0$, then $E(\mathcal{C})$ measures the total length of the curve.

%\cite{horn,kallay1986,kallay1987}. In the field of image processing, the Euler's elastica energy was first utilized by Nitzberg, Mumford, and Shiota \cite{nitzberg1993} to segment an image into objects with different depths. Then, Masnou and Morel in \cite{masnou1998} extended the variational framework of \cite{nitzberg1993} to the level line structures, aiming to recover the missing parts of a gray level image which is named as \textit{disocclusion}. Following the work \cite{masnou1998}, to deal with the problems of image restoration, Ambrosio and Masnou \cite{ambrosio2003,ambrosio2004} advocated to regularizer all level curves of an image $v$ defined on a domain $\Omega$, then the Euler's elastica energy is defined as follows,
% To minimize (\ref{eq:elastica_term}), many methods are proposed in different decades, see
For \textit{imaging applications}, the generalized Euler elastica energy is defined by
\begin{equation}\label{eq:elastica1}
E(v) = \int_{\Omega}\left(a + b \left|\nabla \cdot \frac{\nabla v}{|\nabla v|} \right|^2 \right)|\nabla v|dx,
\end{equation}
where in (\ref{eq:elastica1}), $\Omega$ is a bounded domain of $\mathbf{R}^2$ (a rectangle, typically), $a$ and $b$ are like in (\ref{eq:elastica_term}), $v$ is a function of two variables belonging to an appropriate functional space containing (in principle) the underlying image, and $dx = dx_{1}dx_{2}$.

Recently, the Euler elastica energy defined by (\ref{eq:elastica1}) found applications in image processing, such as: denoising \cite{tai2011,zhang2017}, segmentation \cite{zhu2013,duan2014,zhangchen2016,bae2017}, inpainting \cite{chankang2002,britochen2010,tai2011,yashtini2016}, zooming \cite{tai2011}, illusory contour \cite{masnou2006,kang2014,taiduan2017}, segmentation with depth \cite{nitzberg1993,esedoglu2003,zhuchan2006}.  In \cite{chankang2002}, Chan et al. discuss the mathematical foundation of the Euler elastica model and its mathematical properties, motivated by applications to image inpainting. In addition, the authors of \cite{chankang2002} discussed also a numerical PDE method, in order to solve the associated nonlinear problem. In \cite{tai2011}, Tai et al. proposed an augmented Lagrangian method (ALM) to handle the Euler elastica energy and applied the resulting algorithm to the solution of imaging problems in denoising, inpainting, and zooming. More recently, Zhang et al. proposed in \cite{zhang2017} a linearized strategy to speed up ALM, and applied it to the solution of image denoising problems. In \cite{yashtini2016}, two numerical algorithms were proposed to solve inpainting related problems involving the Euler elastica energy (\ref{eq:elastica1}): The first algorithm is an improved variant of the ALM based algorithm discussed in \cite{tai2011}. The second algorithm is obtained by applying a split-Bregman method to a linearized elastica model proposed in \cite{bae2011}. Following an idea from \cite{masnou1998}, Masnou and Morel proposed in \cite{masnou2006} a novel method to handle the elastica energy functional and applied it to the solution of illusory contour problems. In \cite{kang2014}, Kang et al. used the Euler elastica energy as an effective tool to fuse the scattered corner bases. In \cite{taiduan2017}, Tai and Duan combined level set and binary representation of interfaces to solve, via the Euler elastica model, inpainting, segmentation and illusory control problems. In \cite{BrediesPock2015}, Bredies, Pock and Wirth suggested using as smoothing functional a convex, lower semi-continuous approximation of the Euler elastica energy and applied this approximation to the solution of some imaging problems: combined with tailored discretization of measures; the functional introduced in \cite{BrediesPock2015} has produced promising results.

%
%??? Add paper of Shapiro at al and Chambolle-Pock.
%
%In order to solve the image restoration problem by the Euler's elastica energy, we need to minimize the following functional:
Taking \textit{image restoration} as an illustration, in order to solve the image restoration problem, via Euler's elastica energy, we need to solve the following minimization problem
\begin{equation}\label{eq:elastica}
\min_{v} \left[ \int_{\Omega} \left(a + b \left|\nabla \cdot \frac{\nabla v}{|\nabla v|} \right|^2 \right)|\nabla v|dx + \frac{1}{2}\int_{\Omega}|f - v|^2 dx\right],
\end{equation}
where  $v$ is as in (\ref{eq:elastica1}), and $f$ is the image we are trying to denoise. The first term in the functional in (\ref{eq:elastica}) is a \textit{regularizing} one; it captures the image geometrical features. The second term is the \textit{fidelity} one; it enforces the underlying image to be close to $f$.

\textit{The main goal of this article is to develop a robust, stable and `almost' parameter free method to solve problem (\ref{eq:elastica}), and close variants of it}.

The non-convexity, the non-smoothness, and the high-order of the derivatives it contains, make the fast and robust solution of problem (\ref{eq:elastica}) a very challenging task. So far, there are only few methods to solve problems such (\ref{eq:elastica}); let us mention among them: two graph-cut based methods (\cite{el2010,bae2011}), an integer linear programming (ILP) method ( \cite{Schoenemann2009}), a method based on the approximation of the Euler elastica energy (\cite{bredies2015}), and the augmented Lagrangian method (ALM) (see, e.g., \cite{tai2011}). The method discussed in \cite{tai2011} is a particular realization of the Alternating Method of Multipliers (ADMM), a well-known method from Mathematical Programming (see, e.g., \cite{DGG2007} and references therein for more details). ADMM is a primal-dual method, closely related to the Douglas-Rachford alternating direction method (a well-known operator-splitting method). Following \cite{tai2011}, several extensions were proposed for solving, via the Euler elastica energy functional, a large variety of imaging problems (see \cite{zhu2013,duan2013,yashtini2015,bae2017}). Actually, readers can find further curvature based methods in \cite{bruckstein2001,Greer2004,esedoglu2005,Ballester2001,chambolle2017}. Primal-dual methods have been applied also to derive fast algorithms to handle the total variation (TV) imaging model, introduced in \cite{rof1992} by Rudin, Osher, and Fattorini (ROF). For instance, Droske and Bertozzi in \cite{Droske2010} combined the regularization techniques with active contour models to segment polygonal objects in aerial images. This method could avoid lossing features by using TV-based inverse scale-space techniques on the input data.
See more related literatures \cite{changolub1999,weickert2001,chambolle2004,zhuchan2008,yin2008,esser2009,esser20092,taiwu2009,Weiss2009,steidl2010,wutai2010,Bredieskunisch2010,zhang2011,ChambolleLevine2011,CasellesChambolle2015,ChambolleDuval2016,Kolmogorov2016}  and the references therein for more details.

In this article, we propose a novel and (relatively) simple operator-splitting method for the solution of problem (\ref{eq:elastica}). The principle of the method is very simple: (i) We introduce the vector-valued functions $\mathbf{q} (= \nabla v)$ and  $\bm{\mu} (= \mathbf{q}/|\mathbf{q}|)$. (ii) Using appropriate indicator functionals, we reformulate problem (\ref{eq:elastica}) as an unconstrained minimization problem with respect to the triple $(v,  \mathbf{q}, \bm{\mu})$. (iii) We derive an optimality system and associate with it an initial-value problem (gradient flow). (iv) We use the Lie scheme to time-discretize the above initial value problem and capture its steady state solutions. The sub-problems associated with the Lie scheme fractional steps have either closed form solutions or can be solved by fast dedicated algorithms (such as FFT).  Numerous applications to image smoothing show the efficiency of the proposed method.

When compared to ALM method of \cite{tai2011}, the method introduced in this article has the following advantages:

\begin{itemize}
	\item The time-discretization step is, essentially, the only parameter one has to choose, while ALM
	requires the balancing of three augmentation parameters.
	
	\item The results produced by the new method are less sensitive to parameter choice than those obtained
	by ALM.
	
	\item For the same stopping criterion tolerance, the new method needs less iterations than its ALM
	counterpart. Moreover, the new method has a lower cost per iteration than ALM.
	
\end{itemize}

This article is structured as follows: The novel method is described in Sections \ref{sec:1} and \ref{sect_pro}, while its finite difference implementation is discussed in Section \ref{sec:2}. Section \ref{sec:results} is dedicated to smoothing application, with some experiments designed to show the superiority of the new method. Some conclusions are drawn in Section \ref{sec:conclusion}.

\section{A reformulation of problem (\ref{eq:elastica})}\label{sec:1}

From Section \ref{sec:intro}, Euler's elastica problem reads as

\begin{equation}\label{eq:elastica44}
\min_{v\in \mathcal{V}} \left[\int_{\Omega} \left(a + b \left|\nabla \cdot \frac{\nabla v}{|\nabla v|} \right|^2 \right)|\nabla v|dx + \frac{1}{2} \int_{\Omega}|f - v|^2 dx\right],
\end{equation}
with $\mathcal{V}$ a functional space of the Sobolev's type, typically. If one is willing to replace min by inf in (\ref{eq:elastica44}), one can take $\mathcal{V} = \mathcal{H}^{2}(\Omega)$ if we assume that no specific condition is imposed a priori to $v$ on the boundary $\Gamma$ of $\Omega$. At any rate, the discrete problems largely ignore these functional analysis considerations, and we will say no more about the proper choice of $\mathcal{V}$. An important issue with formulation (\ref{eq:elastica44}) is that it makes no sense on those parts of $\Omega$ where $v$ vanishes. An obvious (and once popular) way to overcome this difficulty is to replace $|\nabla v|$ by $\sqrt{\epsilon^2 + |\nabla v|^2}$,  $\epsilon$ being a small parameter. A more sophisticated way, we borrow from \textit{viscoplasticity} (see, e.g., \cite{DL1976,DGG2007,glowinskiwachs2011}) is to replace $\frac{\nabla v}{|\nabla v|}$ by a vector-valued function $\bm{\mu}$ verifying

\begin{equation}\label{eq:relation_mu_v}
\bm{\mu}\cdot \nabla v = |\nabla v|, ~~ |\bm{\mu}| \leq 1,
\end{equation}
with $|\bm{\mu}| = \sqrt{\mu_{1}^{2} + \mu_{2}^{2}}$, $\forall \bm{\mu} = (\mu_{1}, \mu_{2})$, and then problem (\ref{eq:elastica44}) by

\begin{equation}\label{eq:relation_newform}
\min_{(v, \bm{\mu})\in \mathcal{W}} \left[\int_{\Omega} \left(a + b \left|\nabla \cdot \bm{\mu} \right|^2 \right)|\nabla v|dx + \frac{1}{2}\int_{\Omega}|f - v|^2 dx\right],
\end{equation}
where (formally)

\[
\mathcal{W} = \{ (v, \bm{\mu})\in \mathcal{H}^{1}(\Omega) \times  \mathcal{H}(\Omega, \text{div}), ~~ \bm{\mu}\cdot \nabla v = |\nabla v|, ~~ |\bm{\mu}| \leq 1 \},
\]
with
\[
\mathcal{H}(\Omega, \text{div}) = \{\bm{\mu}\in (\mathcal{L}^{2}(\Omega))^2, \nabla\cdot\bm{\mu}\in  \mathcal{L}^{2}(\Omega)\}.
\]

%A simple, but computationally important, result is provided by the following：
A simple, but computationally important, result is provided by the following:

\noindent\textbf{Proposition 1} ~Suppose that $(u, \bm{\lambda})$ is solution of problem (\ref{eq:relation_newform}). We have then

\begin{equation}\label{prop1}
\int_{\Omega} u dx =  \int_{\Omega} f dx.
\end{equation}

\noindent \textbf{Proof}. Consider the pair $(u + c, \bm{\lambda})$, where $c\in  \mathcal{R}$. Since $\nabla (u + c) = \nabla u$, the pair $(u + c, \bm{\lambda})$ belongs also to $\mathcal{W}$.  Let us denote by $J_{1}$ (resp., $J_{2}$) the left (resp., right) integral in (\ref{eq:relation_newform}). Since  $\nabla (u + c) = \nabla u$ we have $J_{1}(u + c, \bm{\lambda}) = J_{1}(u, \bm{\lambda})$. On the other hand
\begin{equation}\label{proof1}
J_{2}(u + c, \bm{\lambda}) = \frac{1}{2} \int_{\Omega} |u + c - f|^2 dx = J_{2}(u, \bm{\lambda}) + c \int_{\Omega} (u - f) dx + |\Omega| \frac{c^2}{2},
\end{equation}
with $|\Omega|  =$ measure of $\Omega$. The function $u$ being fixed, the quadratic function of $c$ in the right-hand side of (\ref{proof1}) takes its minimal value for $c = c_{m} = \frac{1}{|\Omega|} \int_{\Omega} (f - u) dx$. Suppose that $\int_{\Omega} (f - u) dx \neq 0$; then
\[
J_{2}(u + c_{m}, \bm{\lambda}) < J_{2}(u, \bm{\lambda}),
\]
implying that $(u, \bm{\lambda})$ is not a minimizer of $J_{1} + J_{2}$. We have thus necessarily $\int_{\Omega} u dx =  \int_{\Omega} f dx$. $\Box$

%\begin{remark}
	\textbf{Remark 2.1}: It is a common practice to assume periodicity when working with image processing problems. Proposition 1 still holds if one consider the minimization of the elastica functional in a space of sufficiently smooth periodic functions (functions defined over a two-dimensional torus). In this work, unless otherwise specified, we assume all the functions we are using are periodic in both the $x_1-$ and $x_2-$ directions.
%\end{remark}

Let us define the sets $\Sigma_{f}$ and $S$ by
\[
\Sigma_{f} = \{ \mathbf{q}\in (\mathcal{L}^{2}(\Omega))^2, ~ \exists~ v \in \mathcal{H}^{1}(\Omega)~ \text{such that} ~\mathbf{q} = \nabla v ~\text{and} \int_{\Omega} (v - f)dx = 0\},
\]
and
\[
S = \{ (\mathbf{q}, \bm{\mu})\in  (\mathcal{L}^{2}(\Omega))^2 \times (\mathcal{L}^{2}(\Omega))^2, ~\mathbf{q} \cdot \bm{\mu}= |\mathbf{q}|, ~|\bm{\mu}| \leq 1\}.
\]

There is then equivalence between problem (\ref{eq:relation_newform}) and
\begin{equation}\label{eq:transelastica22}
\min_{(\mathbf{q}, \bm{\mu})\in  (\mathcal{L}^{2}(\Omega))^2 \times (\mathcal{L}^{2}(\Omega))^2} \left[\int_{\Omega} \left(a + b |\nabla \cdot \bm{\mu} |^2 \right)|\mathbf{q}|dx + \frac{1}{2}\int_{\Omega}|v_{\mathbf{q}} - f|^2 dx + I_{\Sigma_{f}}(\mathbf{q}) + I_{S}(\mathbf{q},\bm{\mu})\right],
\end{equation}
where $I_{\Sigma_{f}}$ and $I_{S}$ are indicator functionals defined by

\begin{equation} \label{eq:Inabla_new}
I_{\Sigma_{f}} (\mathbf {q})
=
\left\{
\begin{aligned}
&0,       ~~~~~~~~~ \text{if} ~~  \mathbf{q} \in \Sigma_{f},\\
&+\infty,  ~~~~ \text{if} ~~  \mathbf{q} \in (\mathcal{L}^{2}(\Omega))^2 \backslash \Sigma_{f},\\
\end{aligned}
\right.
\end{equation}

and

\begin{equation} \label{eq:IS_new}
I_{S} (\mathbf{q}, \bm{\mu} )
=
\left\{
\begin{aligned}
&0,       ~~~~~~~~~~ \text{if}  ~~  (\mathbf{q}, \bm{\mu}) \in S,\\
&+\infty, ~~~~~ \text{if}  ~~  (\mathbf{q}, \bm{\mu}) \in (\mathcal{L}^{2}(\Omega))^2 \times (\mathcal{L}^{2}(\Omega))^2 \backslash S,\\
\end{aligned}
\right.
\end{equation}
$v_{\mathbf{q}}$ being the unique solution of the following %Neumann
problem
\begin{equation}\label{eq:eqpro}
\left\{
\begin{aligned}
&\nabla^{2} v_{\mathbf{q}} = \nabla\cdot \mathbf{q}   ~~ \text{in} ~ \Omega,\\
&
%(\nabla v_{\mathbf{q}} - \mathbf{q})\cdot \mathbf{n} = 0 ~~ \text{on} ~\partial\Omega,
\int_{\Omega} v_{\mathbf{q}} dx = \int_{\Omega} f dx.\\
\end{aligned}
\right.
\end{equation}
%$\mathbf{n}$ being the outward unit normal vector at the boundary $\Gamma$ of $\Omega$.
If we assume periodicity of all the functions, the function $v_{\mathbf{q}}$ is also periodic. Without periodicity, we will need to add a boundary condition to guarantee the uniqueness of the above problem (a typical one being $(\nabla v_{\mathbf{q}} - \mathbf{q})\cdot \mathbf{n} = 0$ on $\partial \Omega$).

\section{An Operator-Splitting Method for the Solution of Problem (\ref{eq:transelastica22})} \label{sect_pro}
\subsection{Optimality conditions and associated dynamical flow problem}

Let us denote by $J_1$ and $J_2$ the functionals defined by

\begin{equation} \label{eq:split}
\left\{
\begin{aligned}
&J_{1}(\mathbf{q}, \bm{\mu}) = \int_{\Omega} \left(a + b |\nabla \cdot \bm{\mu}|^2 \right)|\mathbf{q}|dx, \\
&J_{2}(\mathbf{q}) = \frac{1}{2}\int_{\Omega} |v_{\mathbf{q}} - f|^2dx,
\end{aligned}
\right.
\end{equation}
and suppose that $(\mathbf{p}, \bm{\lambda})$ is a minimizer of the functional in (\ref{eq:transelastica22}). We have then $u = v_{\mathbf{p}}$ and the following system of (necessary) optimality conditions holds:

\begin{equation} \label{eq:split2}
\left\{
\begin{aligned}
&\partial_{\mathbf{q}}J_{1}(\mathbf{p}, \bm{\lambda})  + D J_{2}(\mathbf{p}) + \partial I_{\Sigma_{f}}(\mathbf{p}) + \partial_{\mathbf{q}} I_{S}(\mathbf{p}, \bm{\lambda}) \ni \mathbf{0}, \\
&D_{\bm{\mu}} J_{1}(\mathbf{p}, \bm{\lambda}) + \partial_{\bm{\mu}} I_{S}(\mathbf{p}, \bm{\lambda}) \ni \mathbf{0},
\end{aligned}
\right.
\end{equation}
where the $\mathit{Ds}$ (resp., the $\partial s$) denotes classical differentials (resp., generalized differentials (subdifferentials in the case of non-smooth convex functionals, $I_{\Sigma_{f}}$ being a typical one)). We associate with (\ref{eq:split2}) the following initial value problem (dynamical flow):

\begin{equation} \label{eq:split3}
\left\{
\begin{aligned}
&\frac{\partial \mathbf{p}}{\partial t} + \partial_{\mathbf{q}}J_{1}(\mathbf{p}, \bm{\lambda}) + D J_{2}(\mathbf{p}) + \partial I_{\Sigma_{f}}(\mathbf{p}) + \partial_{\mathbf{q}}I_{S}(\mathbf{p}, \bm{\lambda}) \ni \mathbf{0}~ \text{in} ~ \Omega\times (0, +\infty),\\
&\gamma \frac{\partial \bm{\lambda}}{\partial t} + D_{\bm{\mu}}J_{1}(\mathbf{p}, \bm{\lambda}) +  \partial_{\bm{\mu}}I_{S}(\mathbf{p}, \bm{\lambda}) \ni \mathbf{0} ~ \text{in} ~ \Omega\times (0, +\infty),\\
& \left(\mathbf{p}(0), \bm{\lambda}(0) \right) = (\mathbf{p}_{0}, \bm{\lambda}_{0}),
\end{aligned}
\right.
\end{equation}
with  $\gamma > 0$ (the choice of $\gamma$ will be discussed in Section \ref{sec:3.5} ).

The rich structure of problems (\ref{eq:split2}), (\ref{eq:split3}) suggests to solve (\ref{eq:split2}) via the computation of the steady state solutions of (\ref{eq:split3}) using a time discretization method based on operator-splitting. This approach will be discussed in Section \ref{sec:3_2}.

%\begin{remark}
	\textbf{Remark 3.1}: We advocate taking $(\mathbf{p}_{0}, \bm{\lambda}_{0})\in S$ in (\ref{eq:split3}).
%\end{remark}

\subsection{An operator-splitting method for the solution of the dynamical flow problem (\ref{eq:split3}) } \label{sec:3_2}

Following \cite{glowinskiosher2017} (see also \cite{glowinski2017} for applications of operator-splitting methods to Imaging) we will use a \textit{Lie scheme} to time-discretize problem (\ref{eq:split3}). Let  $\tau (> 0)$ be a time discretization step; we denote $(n + \alpha)\tau$ by $t^{n+\alpha}$. Among the many possible splitting schemes of the Lie type one can employ to solve problem (\ref{eq:split3}), we advocate the one below:

\begin{equation} \label{eq:initial}
(\mathbf{p}^{0}, \bm{\lambda}^{0}) = (\mathbf{p}_{0}, \bm{\lambda}_{0}).
\end{equation}

%%%%%%%%%%1st Eq%%%%%%%%%%%%
\underline{\textit{$1^{st}$ Fractional step}}: Solve

\begin{equation} \label{eq:1st}
\left\{
\begin{array}{cc}

\left\{
\begin{aligned}
&\frac{\partial \mathbf{p}}{\partial t} + \partial_{\mathbf{q}}J_{1}(\mathbf{p}, \bm{\lambda}) \ni \mathbf{0} \\
&\gamma \frac{\partial \bm{\lambda}}{\partial t} + D_{\bm{\mu}}J_{1}(\mathbf{p}, \bm{\lambda}) = \mathbf{0} \\
\end{aligned}
\right.
&\text{in} ~~\Omega\times (t^{n}, t^{n+1}), \\&   \\

(\mathbf{p}(t^{n}), \bm{\lambda}(t^{n})) = (\mathbf{p}^{n}, \bm{\lambda}^{n}), &

\end{array}
\right.
\end{equation}

and set

\begin{equation} \label{eq:1staa}
(\mathbf{p}^{n+1/3}, \bm{\lambda}^{n+1/3}) = (\mathbf{p}(t^{n+1}), \bm{\lambda}(t^{n+1})).
\end{equation}

%%%%%%%%%%2nd Eq%%%%%%%%%%%%

\underline{\textit{$2^{nd}$ Fractional step}}: Solve

\begin{equation} \label{eq:2ndt}
\left\{
\begin{array}{cc}

\left\{
\begin{aligned}
&\frac{\partial \mathbf{p}}{\partial t}  + \partial_{\mathbf{q}}I_{S}(\mathbf{p}, \bm{\lambda}) \ni \mathbf{0}\\
&\gamma \frac{\partial \bm{\lambda}}{\partial t} +  \partial_{\bm{\mu}}I_{S}(\mathbf{p}, \bm{\lambda}) \ni \mathbf{0} \\
\end{aligned}
\right.
&\text{in} ~~\Omega\times (t^{n}, t^{n+1}), \\&   \\

(\mathbf{p}(t^{n}), \bm{\lambda}(t^{n})) = (\mathbf{p}^{n+1/3}, \bm{\lambda}^{n+1/3}), &

\end{array}
\right.
\end{equation}

and set

\begin{equation} \label{eq:2ndaa}
(\mathbf{p}^{n+2/3}, \bm{\lambda}^{n+2/3}) = (\mathbf{p}(t^{n+1}), \bm{\lambda}(t^{n+1})).
\end{equation}

%%%%%%%%%%%3rd Eq%%%%%%%%%%%%%%%%%

\underline{\textit{$3^{rd}$ Fractional step}}: Solve

\begin{equation} \label{eq:3rdt}
\left\{
\begin{array}{cc}

\left\{
\begin{aligned}
&\frac{\partial \mathbf{p}}{\partial t} + D J_{2}(\mathbf{p}) + \partial I_{\Sigma_{f}}(\mathbf{p})  \ni \mathbf{0}\\
&\gamma \frac{\partial \bm{\lambda}}{\partial t} = \mathbf{0} \\
\end{aligned}
\right.
&\text{in} ~~\Omega\times (t^{n}, t^{n+1}), \\&   \\

(\mathbf{p}(t^{n}), \bm{\lambda}(t^{n})) = (\mathbf{p}^{n+2/3}, \bm{\lambda}^{n+2/3}), &

\end{array}
\right.
\end{equation}

and set

\begin{equation} \label{eq:3rdaa}
(\mathbf{p}^{n+1}, \bm{\lambda}^{n+1}) = (\mathbf{p}(t^{n+1}), \bm{\lambda}^{n+2/3}).
\end{equation}

The Lie scheme (\ref{eq:1st})-(\ref{eq:3rdaa}) is only semi-discrete since we have not specified yet how to time-discretize the initial value problems (\ref{eq:1st}), (\ref{eq:2ndt}) and (\ref{eq:3rdt}). In order to so, we suggest using the following time discretization scheme (of the Marchuk-Yanenko type):

\begin{equation} \label{eq:initial2}
(\mathbf{p}^{0}, \bm{\lambda}^{0}) = (\mathbf{p}_{0}, \bm{\lambda}_{0}).
\end{equation}

Then, for $n \geq 0$, $(\mathbf{p}^{n}, \bm{\lambda}^{n}) \rightarrow (\mathbf{p}^{n+1/3}, \bm{\lambda}^{n+1/3}) \rightarrow (\mathbf{p}^{n+2/3}, \bm{\lambda}^{n+2/3}) \rightarrow (\mathbf{p}^{n+1}, \bm{\lambda}^{n+1})$ as follows:

\begin{equation} \label{eq:split4}
\begin{array}{cc}
\left\{
\begin{aligned}
&\frac{\mathbf{p}^{n+1/3} - \mathbf{p}^{n}}{\tau} + \partial_{\mathbf{q}}J_{1}(\mathbf{p}^{n+1/3}, \bm{\lambda}^{n}) \ni \mathbf{0}\\
&\gamma \frac{\bm{\lambda}^{n+1/3} - \bm{\lambda}^{n}}{\tau} + D_{\bm{\mu}}J_{1}(\mathbf{p}^{n+1/3}, \bm{\lambda}^{n+1/3} ) = \mathbf{0}
\end{aligned}
\right.
&\text{in} ~~\Omega\Rightarrow (\mathbf{p}^{n+1/3}, \bm{\lambda}^{n+1/3}), \\&   \\
\end{array}
\end{equation}

\begin{equation} \label{eq:split5}
\begin{array}{cc}
\left\{
\begin{aligned}
&\frac{\mathbf{p}^{n+2/3} - \mathbf{p}^{n+1/3}}{\tau} + \partial_{\mathbf{q}}I_{S}(\mathbf{p}^{n+2/3}, \bm{\lambda}^{n+2/3}) \ni \mathbf{0}\\
&\gamma \frac{\bm{\lambda}^{n+2/3} - \bm{\lambda}^{n+1/3}}{\tau} + \partial_{\bm{\mu}}I_{S}(\mathbf{p}^{n+2/3}, \bm{\lambda}^{n+2/3} ) \ni \mathbf{0}
\end{aligned}
\right.
&\text{in} ~~\Omega\Rightarrow (\mathbf{p}^{n+2/3}, \bm{\lambda}^{n+2/3}), \\&   \\
\end{array}
\end{equation}

\begin{equation} \label{eq:split6}
\begin{array}{cc}
\left\{
\begin{aligned}
&\frac{\mathbf{p}^{n+1} - \mathbf{p}^{n+2/3}}{\tau} + D J_{2}(\mathbf{p}^{n+1})  + \partial I_{\Sigma_{f}}(\mathbf{p}^{n+1}) \ni \mathbf{0}\\
&\gamma \frac{\bm{\lambda}^{n+1} - \bm{\lambda}^{n+2/3}}{\tau} = \mathbf{0}
\end{aligned}
\right.
&\text{in} ~~\Omega\Rightarrow (\mathbf{p}^{n+1}, \bm{\lambda}^{n+1}).\\&   \\
\end{array}
\end{equation}

In the following subsections we are going to discuss the solution of the various sub-problems encountered when applying scheme (\ref{eq:initial2})-(\ref{eq:split6}) to the solution of problem (\ref{eq:transelastica22}).

\subsection{Computing $\mathbf{p}^{n+1/3}$ from (\ref{eq:split4})}\label{sec:2.2}
The multi-valued equation verified by $\mathbf{p}^{n+1/3}$ in (\ref{eq:split4}) is nothing but the (formal) Euler-Lagrange equation of the following minimization problem

\begin{equation} \label{eq:q_n13pro}
\mathbf{p}^{n+1/3} = \arg\min_{\mathbf{q}\in (\mathcal{L}^{2}(\Omega))^2} ~ \left[\frac{1}{2} \int_{\Omega} |\mathbf{q} - \mathbf{p}^n|^2dx + \tau \int_{\Omega} \left(a + b |\nabla \cdot \bm{\lambda}^n|^2\right)|\mathbf{q}|dx\right].
\end{equation}

Problems such as (\ref{eq:q_n13pro}) are very common in Image Processing and Viscoplasticity. The closed form solution of problem (\ref{eq:q_n13pro}) is given by (see \cite{fortinglowinski1983,GlowLe1989,Donoho1995,wangyang2008,tai2011}):

\begin{equation} \label{eq:q_n13pro_solut}
\mathbf{p}^{n+1/3} = \textbf{\text{max}}\bigg \{0, 1 - \frac{c}{|\mathbf{p}^{n}|}\bigg\} \mathbf{p}^{n},
\end{equation}
where $c = \tau a + \tau b |\nabla \cdot \bm{\lambda}^n|^2$.\\

\subsection{Computing $\lambda^{n+1/3}$ from (\ref{eq:split4})}\label{sec:3_4}
The equation verified by $\bm{\lambda}^{n+1/3}$ in  (\ref{eq:split4}) is the (formal) Euler-Lagrange equation of the following minimization problem

%The $\bm{\lambda}^{n+1/3}$ problem is related to the second equation of (\ref{eq:11}). The solution  $\bm{\lambda}^{n+1/3}$ of (\ref{eq:11})  is a solution  of the following minimization problem
\begin{equation} \label{eq:lambda_n13pro}
\bm{\lambda}^{n+1/3} = \arg\min_{\bm{\mu}\in (\mathcal{L}^{2}(\Omega))^2} ~ \left[\gamma \int_{\Omega} \frac{|\bm{\mu} - \bm{\lambda}^n|^2}{2\tau}dx + J_{1}(\bm{\mu}, \mathbf{p}^{n+1/3})\right],
\end{equation}
where $\bm{\lambda}^n$ and $\mathbf{p}^{n+1/3}$ are known.

From the Euler-Lagrangian equation of (\ref{eq:lambda_n13pro}), we get that the solution $\bm{\lambda}^{n+1/3}$ is the solution of following linear equation:
% with an imposed periodic boundary condition:
\begin{equation} \label{eq:u_n13pro_solut1}
\gamma \frac{\bm{\lambda}^{n+1/3} - \bm{\lambda}^{n}}{\tau} - \nabla(2b|\mathbf{p}^{n+1/3}|\nabla\cdot{\bm{\lambda}}^{n+1/3}) = 0 \mbox{  in  }\Omega.
\end{equation}
This is a vector equation.
We shall use periodic boundary condition for the above equation and also other subproblems coming later. It is common to use this kind of boundary condition for image processing problems. It is easy to justify this approach by assuming that the image data is defined on a 2D torus, for example.
Under the periodic boundary condition,  (\ref{eq:u_n13pro_solut1}) can be efficiently and easily solved by the FFT,  see \cite[\S  3.2.5]{tai2011} and \cite{wangyang2008}.

\subsection{Computing $\left(\mathbf{p}^{n+2/3}, \lambda^{n+2/3}\right)$ from (\ref{eq:split5})}\label{sec:3.5}

\subsubsection{Decomposition of problem (\ref{eq:split5})}\label{sec:3.5.1} %(\ref{eq:split5})
One can view system (\ref{eq:split5}) as the Euler-Lagrange equation of the following minimization problem:

\begin{equation} \label{eq:q_n23pro}
\min_{(\mathbf{q}, \bm{\mu})\in S} ~ \left[\int_{\Omega} |\mathbf{q} - \mathbf{p}^{n+1/3}|^2 dx +\gamma \int_{\Omega} |\bm{\mu} - \bm{\lambda}^{n+1/3}|^2 dx \right].
\end{equation}

Problem (\ref{eq:q_n23pro}) can be solved point-wise, reducing, a.e. on $\Omega$, to the following finite dimensional constrained minimization problem:
\begin{equation} \label{eq:q_n23pro1}
(\mathbf{p}^{n+2/3}(x), \bm{\lambda}^{n+2/3}(x)) = \text{argmin}_{(\mathbf{q}, \bm{\mu})\in \sigma} j_{n+1/3}(\mathbf{q}, \bm{\mu}; x),
\end{equation}
where $\sigma =  \{ (\mathbf{q}, \bm{\mu})\in  \mathbf{R}^{2} \times \mathbf{R}^{2}, ~\mathbf{q} \cdot \bm{\mu}= |\mathbf{q}|, ~|\bm{\mu}| \leq 1\}$ and
\[
j_{n+1/3}(\mathbf{q}, \bm{\mu}; x) = \left|\mathbf{q} - \mathbf{p}^{n+1/3}(x)\right|^2 + \gamma\left|\bm{\mu} - \bm{\lambda}^{n+1/3}(x)\right|^2, \forall (\mathbf{q}, \bm{\mu})\in \mathbf{R}^2 \times \mathbf{R}^2.
\]
Let us define $\sigma_0$ and $\sigma_1$ by
\[
\sigma_{0} = \{(\mathbf{q}, \bm{\mu})\in \mathbf{R}^2 \times \mathbf{R}^2, \mathbf{q} = \mathbf{0},  |\bm{\mu}|\leq 1\},
~~ \sigma_{1} = \{(\mathbf{q}, \bm{\mu})\in \mathbf{R}^2 \times \mathbf{R}^2,  \mathbf{q} \neq \mathbf{0}, \mathbf{q} \cdot \bm{\mu}= |\mathbf{q}|, ~ |\bm{\mu}|= 1\}.
\]
We clearly have  $\sigma = \sigma_{0} \cup \sigma_{1}$, implying that to compute $\left(\mathbf{p}^{n+2/3}(x), \bm{\lambda}^{n+2/3}(x) \right)$, we may proceed as follows:

\textbf{(i)} Solve the following two uncoupled minimization problems
\begin{align}  \label{eq:tai1}
& \left(\mathbf{p}_{0}^{n+2/3}(x), \bm{\lambda}_{0}^{n+2/3}(x)\right) = \text{argmin}_{(\mathbf{q}, \bm{\mu})\in \sigma_{0}} j_{n+1/3}(\mathbf{q}, \bm{\mu}; x),
\\
\label{eq:tai2}
& \left(\mathbf{p}_{1}^{n+2/3}(x), \bm{\lambda}_{1}^{n+2/3}(x) \right) = \text{argmin}_{(\mathbf{q}, \bm{\mu})\in \sigma_{1}} j_{n+1/3}(\mathbf{q}, \bm{\mu}; x),
\end{align}

\textbf{(ii)} Choose the one that gives the smallest value to $j_{n+1/3}$ as the minimizer of (\ref{eq:q_n23pro1}), i.e.
\begin{align} \label{eq:q_n23pro1-a}
\left(\mathbf{p}^{n+2/3}(x), \bm{\lambda}^{n+2/3}(x)\right)
& = \text{argmin} \bigg [j_{n+1/3}( \mathbf{p}_{0}^{n+2/3}, \bm{\lambda}_{0}^{n+2/3}; x), \\
&  j_{n+1/3}(\mathbf{p}_{1}^{n+2/3}, \bm{\lambda}_{1}^{n+2/3}; x)\bigg ], \text{a.e. on}~\Omega. \nonumber
\end{align}

In what follows, we first introduce a strategy of adaptively choosing $\gamma$ in Section \ref{sect:S0a}. After that,  in Section \ref{sect:S0} and Section \ref{sect:S1}, we will discuss the minimization of the functional in (\ref{eq:q_n23pro1}) over $\sigma_0$ and $\sigma_1$, respectively.

\subsubsection{Selection of the parameter $\gamma$} \label{sect:S0a}

We intend to select the parameter $\gamma$ so that the two terms in $j_{n+1/3}$ are balanced. We note that
\[ \llambda (t) = \frac {\pp(t)}{|\pp(t) |}. \]
Thus
\begin{equation} \label{eq:dlambda}
\frac {\partial \llambda} {\partial t}  = \lim_{\tau \to 0} \frac 1 \tau \bigg ( \frac {\pp(t+\tau )}{|\pp(t+\tau ) |} - \frac {\pp(t )}{|\pp(t ) |}\bigg ).
\end{equation}
Due to the following relation:
\begin{align*}
\bigg | \frac \pp {|\pp| }  - \frac \qq {|\qq| } \bigg |^2 = \frac {|\pp|^2}{|\pp|^2 } + \frac {|\qq|^2}{|\qq|^2 } - \frac {2 \pp \cdot \qq }{|\pp| |\qq|}
=   2 \bigg ( 1 - \frac {\pp\cdot \qq} {|\pp| |\qq|} \bigg )
\le 2 \bigg ( 1 - \frac {2 \pp\cdot \qq} {|\pp|^2 +  |\qq|^2 }\bigg)
= 2 \frac {|\pp - \qq|^2}{|\pp|^2 +  |\qq|^2 },
\end{align*}
one has
\[ \bigg | \frac {\pp(t+\tau )}{|\pp(t+\tau )  |} - \frac {\pp(t )}{|\pp(t ) |} \bigg |
\le \frac {\sqrt{2} |\pp(t+\tau ) - \pp(t) |}{\sqrt{|\pp(t+\tau ) |^2 + |\pp(t )|^2 }}.
\]
Let $\tau \to 0$, we get from (\ref{eq:dlambda}) that
\[
\bigg |\frac {\partial \llambda} {\partial t}\bigg | \le \frac 1  {|\pp|}  {\bigg | \frac {\partial \pp } {\partial t}\bigg |}.
\]
For small $\tau$, the minimizer of (\ref{eq:q_n23pro}) verifies
\begin{equation} \label{eq:q_n23pro2-a}
\frac{|\mathbf{p}^{n+2/3} - \mathbf{p}^{n+1/3}|^2}{2\tau} + \gamma \frac{|\bm{\lambda}^{n+2/3} - \bm{\lambda}^{n+1/3}|^2}{2\tau}
\approx \frac \tau 2 \bigg (
\bigg | \frac{\partial \mathbf{p}}{\partial t} (t^{n+1/3})\bigg |^2 + \gamma \bigg |\frac{\partial \bm{\lambda}}{\partial t} (t^{n+1/3}) \bigg |^2  \bigg ).
\end{equation}
According to the above estimate, to balance these two terms, we just need to choose
\[\gamma = | \pp^{n+ 1/3}|^2. \]
In order to avoid the case $| \pp^{n+ 1/3}| \approx 0$, we choose in practice
\begin{equation}\label{eq:gamma}
\gamma =\max ( | \pp^{n+ 1/3}|^2, \hat{\alpha}), \end{equation}
where $\hat{\alpha}$ is a given small number. In this work, we empirically choose $\hat{\alpha} = \sqrt{\tau}$.

\subsubsection{Minimizing the functional in (\ref{eq:q_n23pro1}) over $\sigma_0$} \label{sect:S0}
Over $\sigma_0$ the minimization problem (\ref{eq:tai1}) reduces to
\begin{equation}\label{eq:tai3}
\min_{\bm{\mu}\in \mathbf{R}^2, |\bm{\mu}| \leq 1} \left|\bm{\mu} - \bm{\lambda}^{n+1/3}(x)\right|.
\end{equation}
Clearly, the solution of problem (\ref{eq:tai3}) is given by
\begin{equation}\label{eq:selectgamma_sol}
\bm{\lambda}_{0}^{n+1/3}(x) = \frac{\bm{\lambda}^{n+1/3}(x)}{\max(1, |\bm{\lambda}^{n+1/3}(x)|)}.
\end{equation}
Concerning $\mathbf{p}_{0}^{n+1/3}(x)$, we have, obviously,
\begin{equation}\label{eq:p013}
\mathbf{p}_{0}^{n+1/3}(x) = \mathbf{0}.
\end{equation}

\subsubsection{Minimizing the  functional in  (\ref{eq:q_n23pro1}) over $\sigma_1$} \label{sect:S1}
Over $\sigma_1$, the minimization problem (\ref{eq:tai2}) reduces to:
\begin{equation}\label{eq:sigma1}
\inf_{(\mathbf{q},\bm{\mu})\in \mathbf{R}^2\times \mathbf{R}^2, \mathbf{q} \neq \mathbf{0}, \mathbf{q}\cdot \bm{\mu}=|\mathbf{q}|, |\bm{\mu}| = 1} \left[|\mathbf{q} - \mathbf{p}^{n+1/3}(x)|^2 + \gamma|\bm{\mu} - \bm{\lambda}^{n+1/3}(x)|^2 \right].
\end{equation}

For notational simplicity, we introduce $\mathbf{x}$ and $\mathbf{y}$  defined by $\mathbf{x}=\mathbf{p}^{n+1/3}(x)$ and $\mathbf{y}=\bm{\lambda}^{n+1/3}(x)$, respectively. Using this notation and taking relation $|\bm{\mu}| = 1$ into account, problem (\ref{eq:sigma1}) takes the following simplified formulation
\begin{equation}\label{eq:sigma1_2}
\inf_{(\mathbf{q},\bm{\mu})\in \mathbf{R}^2\times \mathbf{R}^2, \mathbf{q} \neq \mathbf{0}, \mathbf{q}\cdot \bm{\mu}=|\mathbf{q}|, |\bm{\mu}| = 1} \left[\frac{1}{2}|\mathbf{q}|^2 - \mathbf{q}\cdot \mathbf{x} - \gamma\bm{\mu}\cdot\mathbf{y}\right].
\end{equation}

Let us denote $|\mathbf{q}|$ by $\theta$; since $\bm{\mu} = \mathbf{q}/|\mathbf{q}|$, the above relations imply that
\begin{equation}\label{eq:imply}
\mathbf{q} = \theta \bm{\mu}, ~~\theta > 0.
\end{equation}

Relation (\ref{eq:imply}) allows us to replace problem (\ref{eq:sigma1_2}) by the following constrained minimization problem in $\mathbf{R^3}$
\begin{equation}\label{eq:sigma1_3}
\inf_{(\theta,\bm{\mu})\in \mathbf{R}\times \mathbf{R}^2, \theta>0, |\bm{\mu}| = 1} \left[\frac{1}{2}\theta^2 - \theta\bm{\mu}\cdot \mathbf{x} - \gamma\bm{\mu}\cdot\mathbf{y}\right].
\end{equation}

In order to solve problem (\ref{eq:sigma1_3}), we observe that the above problem is equivalent to
\begin{equation}\label{eq:sigma1_4}
\inf_{\theta>0} \min_{\bm{\mu}\in \mathbf{R}^2, |\bm{\mu}| = 1}\left[\frac{1}{2}\theta^2 - \theta\bm{\mu}\cdot \mathbf{x} - \gamma\bm{\mu}\cdot\mathbf{y}\right].
\end{equation}

In order to minimize on a closed set of $\mathbf{R}^3$, the problem that we finally consider is the following variant of problem (\ref{eq:sigma1_4})
\begin{equation}\label{eq:sigma1_5}
\min_{\theta\geq 0} \min_{\bm{\mu}\in \mathbf{R}^2, |\bm{\mu}| = 1}\left[\frac{1}{2}\theta^2 - \theta\bm{\mu}\cdot \mathbf{x} - \gamma\bm{\mu}\cdot\mathbf{y}\right].
\end{equation}

The parameter $\theta$ being fixed, the solution $\bm{\mu}^{*}(\theta)$ of problem
\[
\min_{\bm{\mu}\in \mathbf{R}^2, |\bm{\mu}| = 1}\left[\frac{1}{2}\theta^2 - \theta\bm{\mu}\cdot \mathbf{x} - \gamma\bm{\mu}\cdot\mathbf{y}\right].
\]
is given by $\bm{\mu}^{*}(\theta) = \frac{\theta\mathbf{x}+\gamma\mathbf{y}}{|\theta\mathbf{x}+\gamma\mathbf{y}|}$, implying that problem (\ref{eq:sigma1_5}) reduces to
\begin{equation}\label{eq:sigma1_6}
\min_{\theta\geq 0} \left[\frac{1}{2}\theta^2 - |\theta\mathbf{x} + \gamma\mathbf{y}|\right].
\end{equation}

There are many ways to solve problem (\ref{eq:sigma1_6}), such as Newton's method, bisection or golden section methods, and a variety of fixed point methods (\cite{burden1985}). The method we have chosen is a fixed point one and has shown fast convergence properties. Let us denote by $E$ the function defined by
\[
E(\theta) = \frac{1}{2}\theta^2 - |\theta\mathbf{x} + \gamma\mathbf{y}|.
\]

We clearly have
\[
\frac{dE}{d\theta}(\theta) = \theta - \frac{\mathbf{x}\cdot(\theta\mathbf{x}+\gamma\mathbf{y})}{|\theta\mathbf{x}+\gamma\mathbf{y}|}.
\]

In order to solve equation $\frac{dE}{d\theta}(\theta) = 0$, we advocate the simple following \textit{fixed point} method:

\begin{equation} \label{eq:fixpoint}
\left\{
\begin{aligned}
&\theta^{0} = |\mathbf{x}|\\
&\text{For}~ k\geq 0, ~\theta^{(k)} \rightarrow \theta^{(k+1)} ~\text{as~follows} \\
&\theta^{(k+1)} = \max(0, \frac{\mathbf{x}\cdot(\theta^{(k)}\mathbf{x}+\gamma\mathbf{y})}{|\theta^{(k)}\mathbf{x}+\gamma\mathbf{y}|}).
\end{aligned}
\right.
\end{equation}

A more detailed presentation of our fixed point method reads as:\\

\noindent
\begin{tabular}{l}
	\toprule[1.5pt]
	\noindent\textbf{Algorithm 1:} Fixed point solution of problem (\ref{eq:sigma1_6})~~~~~~~~~~~~~~~~~~~~~~~~~~~~~~~~~~~~~~~~~~~~~~~\\
	\midrule[1.5pt]
	\noindent
	\textbf{Input:} ~~~$\mathbf{x}$, $\mathbf{y}$, $\gamma$\\
	\textbf{Output:} $\theta^{*}$ \\
	\textbf{Initialization:} $\theta^{(0)}$ = $|\mathbf{x}|$, $k = 0$\\
	
	\textbf{While:} $|\theta^{(k+1)} - \theta^{(k)}| > tol$ ~\text{and} ~$k < \text{M}_{it}$\\
	\indent 1) compute $\theta^{(k+1)}$ by \\
	\indent ~~~~~~~~~~~~ $\theta^{(k+1)} =  \max\left(0, \frac{\mathbf{x}\cdot(\theta^{(k)}\mathbf{x}+\gamma\mathbf{y})}{|\theta^{(k)}\mathbf{x}+\gamma\mathbf{y}|}\right)$\\
	\indent 2) $k = k + 1$\\
	\textbf{Endwhile.} \\
	
	\indent 3) One gets the final $\theta^{*}$ when iterations stop\\
	\bottomrule[1.5pt]
\end{tabular}\\

In Algorithm 1, $tol$ and $M_{it}$ denote a positive tolerance value and the maximum number of iterations, respectively. Actually, Algorithm 1 is not sensitive to these values. For all the experiments reported in this article we took $tol = 10^{-3}$ and $M_{it} = 100$.

Once $\theta^{*}$ is known, we obtain the vectors $\bm{\lambda}_{1}^{n+2/3}(x)$ and
$\mathbf{p}_{1}^{n+2/3}(x)$ (we defined them in Section \ref{sec:3.5.1}) via the following relations:
\begin{equation} \label{eq:lambda23_1}
\bm{\lambda}_{1}^{n+2/3}(x) = \frac{\theta^{*}\mathbf{p}^{n+1/3}(x)+\gamma\bm{\lambda}^{n+1/3}(x)}{|\theta^{*}\mathbf{p}^{n+1/3}(x)+\gamma\bm{\lambda}^{n+1/3}(x)|},
\end{equation}
and
\begin{equation} \label{eq:p23_1}
\mathbf{p}_{1}^{n+2/3}(x) = \theta^{*}\bm{\lambda}_{1}^{n+2/3}(x),
\end{equation}
respectively. A more rigorous notation would have been to use $\theta_{n}^{*}(x)$  instead of $\theta^{*}$, since the solution of problem (\ref{eq:sigma1_6}) varies with $x$ and $n$ (we recall that, in (\ref{eq:sigma1_6}), $\mathbf{x} = \mathbf{p}^{n+1/3}(x)$ and  $\mathbf{y} = \bm{\lambda}^{n+1/3}(x)$).

Once we have gotten $\left(\mathbf{p}_{1}^{n+2/3}(x), \bm{\lambda}_{1}^{n+2/3}(x) \right)$ from (\ref{eq:selectgamma_sol})--(\ref{eq:p013})  and $\left(\mathbf{p}_{1}^{n+2/3}(x), \bm{\lambda}_{1}^{n+2/3}(x) \right) $ from (\ref{eq:lambda23_1})--(\ref{eq:p23_1}), we obtain the minimizer of (\ref{eq:q_n23pro1}) through (\ref{eq:q_n23pro1-a}).

\subsection{Computing  $\mathbf{p}^{n+1}$ and $\lambda^{n+1}$ from (\ref{eq:split6})}\label{sec:3.6}
We clearly have
\begin{equation} \label{eq:lambda33}
\bm{\lambda}^{n+1} = \bm{\lambda}^{n+2/3}.
\end{equation}

On the other hand, the multi-valued equation verified by $\mathbf{p}^{n+1}$  in (\ref{eq:split6}) is the Euler-Lagrange equation of the following minimization problem
\begin{equation}\label{eq:p33_1}
\mathbf{p}^{n+1} = \arg \min_{\mathbf{q} \in \Sigma_f } ~ \left[\frac{1}{2}\int_{\Omega}\left|\mathbf{q} - \mathbf{p}^{n+2/3}\right|^2dx + \frac{\tau}{2}\int_{\Omega}\left|v_{\mathbf{q}} - f\right|^2dx\right],
\end{equation}
the function $v_{\mathbf{q}}$  being defined by (\ref{eq:eqpro}).

As we have mentioned in Section \ref{sec:3_4}, we shall use periodic boundary condition for all the subproblems.
%
%\noindent\textbf{Remark 3}.
Suppose that $\Omega$ is the rectangle $(0, L)\times (0, H)$. Next, define $\mathcal{H}_{p}^{1}(\Omega)$, a space of doubly periodic functions, by
\[
\mathcal{H}_{p}^{1}(\Omega) = \{v\in \mathcal{H}^{1}(\Omega); v(0, x_{2}) = v(L, x_{2}),~ \text{a.e. on} ~(0, H); v(x_{1}, 0) = v(x_{1}, H) ~ \text{a.e. on}~ (0, L)\}
\]
From the definition of $\Sigma_f$ (see Section \ref{sec:1}), problem (\ref{eq:p33_1}) is equivalent to

\begin{equation} \label{eq:p33_2}
\left\{
\begin{aligned}
&\mathbf{p}^{n+1} = \nabla u^{n+1} ~~~~\text{with}\\
&u^{n+1} = \arg\min_{v\in \mathcal{H}_p^{1}(\Omega)} \left[\frac{1}{2}\int_{\Omega}\left|\nabla v\right|^2dx + \frac \tau 2\int_{\Omega}\left|v - f\right|^2dx - \int_{\Omega}\mathbf{p}^{n+2/3}\cdot \nabla vdx  \right].
\end{aligned}
\right.
\end{equation}

Function $u^{n+1}$ is the unique solution of the following well-posed linear variational problem in $\mathcal{H}_p^{1}(\Omega)$:
\begin{equation} \label{eq:p33_3}
\left\{
\begin{aligned}
& u^{n+1}\in \mathcal{H}_p^{1}(\Omega),\\
& \int_{\Omega}\nabla u^{n+1}\cdot \nabla vdx + \tau\int_{\Omega}u^{n+1}vdx = \int_{\Omega}\mathbf{p}^{n+2/3}\cdot \nabla vdx + \tau\int_{\Omega}fvdx, ~~\forall v\in \mathcal{H}^{1}(\Omega)
\end{aligned}
\right.
\end{equation}

Problem (\ref{eq:p33_2}), (\ref{eq:p33_3}) has a unique solution which is the weak solution of the following problem:
\begin{equation} \label{eq:p33_5}
\left\{
\begin{aligned}
& -\nabla^{2} u^{n+1} + \tau u^{n+1} = -\nabla\cdot \mathbf{p}^{n+2/3} + \tau f, ~~ \text{in}~\Omega\\
& u^{n+1}(0, x_{2}) = u^{n+1}(L, x_2)~ \text{a.e. on} ~(0, H); ~u^{n+1}(x_1, 0) = u^{n+1}(x_1, H)~ \text{a.e. on} ~(0, L),\\
&\frac{\partial u^{n+1}}{\partial x_1}(0, x_2) = \frac{\partial u^{n+1}}{\partial x_1}(L, x_2)~ \text{a.e.  on} ~(0, H);  \frac{\partial u^{n+1}}{\partial x_2}(x_1, 0) = \frac{\partial u^{n+1}}{\partial x_2}(x_1, H)~ \text{a.e.  on} ~(0, L).
\end{aligned}
\right.
\end{equation}

In the particular case of rectangular domains $\Omega$, many efficient solvers are available for the solution of the finite dimensional analogues of problem (\ref{eq:p33_5}) obtained by finite difference discretization. Among these fast solvers let us mention those based on cyclic reduction and FFT.

\subsection{Summary}\label{sec:3.7}
The subproblems (\ref{eq:split4}), (\ref{eq:split5}) and (\ref{eq:split6}) encountered in our splitting method aim at minimizing consecutively the various components of the elastica cost functional. Our proposed algorithm is summarized in Algorithm 2

\noindent
\begin{tabular}{l}
	\toprule[1.5pt]
	\noindent\textbf{Algorithm 2:} A schematic description of the algorithm solving problem (\ref{eq:elastica44})\\
	\midrule[1.5pt]
	\noindent
	\textbf{Input:} ~The inputted image $f$, the parameters $a$, $b$  and $\tau$. \\
	\textbf{Output:} The computed image $u^{*}$.\\
	\textbf{Initialization:} $n = 0$, $u^{0} = f$, $\mathbf{p}^{0} = \nabla f$,
	$\bm{\lambda}^{0}(x) = \left\{
	\begin{aligned}
	& \mathbf{p}^{0}(x)/|\mathbf{p}^{0}(x)|, ~\text{if}~\mathbf{p}^{0}(x)\neq 0, \\
	& \mathbf{0}, ~~\text{otherwise}.
	\end{aligned}
	\right. x\in \Omega$.\\
	
	\textbf{While:} $\|u^{n+1} - u^{n}\|/\|u^{n+1}\| > tol$ ~\text{and} ~$n < \text{M}_{iter}$\\
	\indent 1. Using the methods discussed in Sections \ref{sec:2.2} and \ref{sec:3_4}, solve system (\ref{eq:split4}) to  \\
	\indent ~~~ obtain $\left(\mathbf{p}^{n+1/3}, \bm{\lambda}^{n+1/3}\right)$.\\
	\indent 2. Use the method discussed in Section \ref{sec:3.5}  to obtain $\left(\mathbf{p}^{n+2/3}, \bm{\lambda}^{n+2/3}\right)$ from (\ref{eq:split5}).\\
	\indent 3. Use the method discussed in Section \ref{sec:3.6} to obtain $\left(u^{n+1}, \mathbf{p}^{n+1}, \bm{\lambda}^{n+1}\right)$ from (\ref{eq:split6}).\\	
	\indent 4. Check convergence and go to the next iteration or stop.\\		
	\textbf{End While.} \\
	\indent If iterations stop, take $u^{*}= u^{n+1}$. \\
	\bottomrule[1.5pt]
\end{tabular}\\

In Algorithm 2, $tol$ is the stopping criterion tolerance, $M_{iter}$ is the maximum of iterations and the norm $||\cdot||$ is $L_2$ norm. All the subproblems encountered when using Algorithm 2 have either closed form solutions or can be solved by dedicated fast solvers. Due to the semi-implicit nature of the operator-splitting scheme, we can use (relatively) large values of $\tau$ and our numerical experiments show that the overall iteration number is (relatively) low. The model parameters $a$ and $b$ have to be given. Finally, the time-discretization step $\tau$ also needs to be provided. We want to say that
$\tau$ is easy to tune. The selection of $\gamma$ was addressed in Section \ref{sect:S0a}; further information will be provided in Section 6 about the choices of all these parameters.

\section{Numerical Discretization}\label{sec:2}
\subsection{Synopsis}\label{sec:syn}

\begin{figure}[!htp]
	\begin{center}
		
		\includegraphics[width=4.5in,height=2.5in]{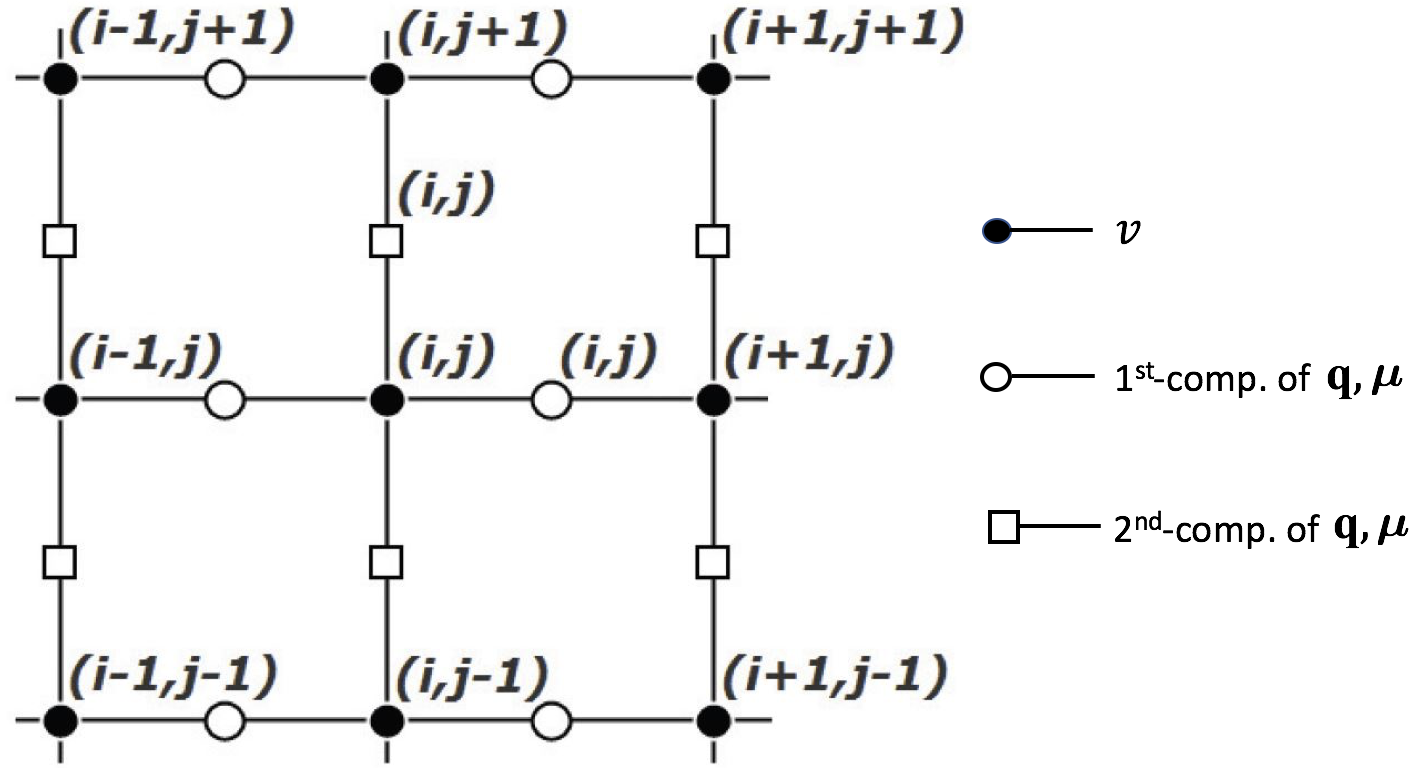}		
		
	\end{center}
	\caption{Indexation of the discrete analogues of the unknown functions $v$ (at the $\bullet$-nodes) and of the first (at the $\circ$-nodes) and second (at the $\square$-nodes) components of the vector-valued functions $\mathbf{q}$ and $\bm{\mu}$. }
	\label{fig:grid}
\end{figure}

As in \cite{tai2011}, we will assume that $\Omega$ is a rectangle. We assume that all functions are periodic in both the $x_1$ ans $x_2$ directions. To discretize the Euler elastica variational problem, we will use staggered grids as visualized in Fig. \ref{fig:grid}. In Fig. \ref{fig:grid}, the unknown function $v$ is discretized at the $\bullet$-nodes, while the first (resp., second) components of $\mathbf{q}$ and $\bm{\mu}$  are discretized at the $\circ$-nodes (resp., $\square$-nodes). Useful  notation will be introduced in Section \ref{sec:usefulnote}. The solution of the discrete sub-problems will be discussed in Sections \ref{sec:dis_pn13}--\ref{sec:dis_plambdan33}.

\subsection{ Some useful discrete operators}\label{sec:usefulnote}
After discretization, we denote by $\Omega_h$ the discrete image domain $\Omega_h = [1:M_1] h \times [1:N_1 ] h$ where $h = L/M_{1} = H/N_{1}$. Note that $\Omega_h$ is a set of $M_{1}N_{1}$ points in $\mathbf{R}^2$. Taking periodicity into account, we define the backward (--) and forward (+) discrete analogues of $\frac{\partial v}{\partial x_1}$ and $\frac{\partial v}{\partial x_2}$ by

\[
\partial_{1}^{-}v(i, j) = \left\{
\begin{aligned}
&(v(i, j) - v(i-1, j))/h, ~~ 1< i \leq M_{1} \\
&(v(1, j) - v(M_1, j))/h, ~~~~ i = 1,
\end{aligned}
\right.
\]

\[
\partial_{2}^{-}v(i, j) = \left\{
\begin{aligned}
&(v(i, j) - v(i, j-1))/h, ~~ 1< j \leq N_{1} \\
&(v(i, 1) - v(i, N_1))/h, ~~~~~ j = 1,
\end{aligned}
\right.
\]

\[
\partial_{1}^{+}v(i, j) = \left\{
\begin{aligned}
&(v(i+1, j) - v(i, j))/h, ~~ 1\leq i < M_{1} \\
&(v(1, j) - v(M_1, j))/h, ~~~~ i = M_1,
\end{aligned}
\right.
\]

\[
\partial_{2}^{+}v(i, j) = \left\{
\begin{aligned}
&(v(i, j+1) - v(i, j))/h, ~~ 1\leq j < N_{1} \\
&(v(i, 1) - v(i, N_1))/h, ~~~~~ j = N_1.
\end{aligned}
\right.
\]

With obvious notation, the discrete forward (+) and backward (--) gradient operators $\nabla^{+}$ and      $\nabla^{-}$  are defined by

\[
\nabla^{\pm}v(i, j) = (\partial_{1}^{\pm}v(i, j), \partial_{2}^{\pm}v(i, j)).
\]

The associated discrete forward (+) and backward (--) divergence operators $\text{div}^{+}$  and $\text{div}^{-}$  are defined (again with obvious notation) by

\[
\text{div}^{\pm}\mathbf{q}(i,j) = \partial_{1}^{\pm}q_{1}(i, j) + \partial_{2}^{\pm}q_{2}(i, j).
\]

If, in particular, a variable defined at the  $\circ$-nodes (resp., $\square$-nodes) needs to be evaluated at the $\square$-node (resp., $\circ$-node) $(i, j)$, it will be done using the following averaging operator:

\begin{equation} \label{eq:www1}
\mathcal{A}_{i,j}^{\square}(\mu_1) = \frac{\mu_{1}(i,j+1) + \mu_{1}(i-1,j+1) + \mu_{1}(i,j) + \mu_{1}(i-1,j)}{4},
\end{equation}
(resp.
\begin{equation} \label{eq:www11}
\mathcal{A}_{i,j}^{\circ}(\mu_{2}) = \frac{\mu_{2}(i+1,j) + \mu_{2}(i,j) + \mu_{2}(i+1,j-1) + \mu_{2}(i,j-1)}{4},
\end{equation}
where $\mu_1$ (resp., $\mu_2$) is defined at the $\circ$-nodes (resp., $\square$-nodes). In order to evaluate the magnitude of $\mathbf{q} = (q_1, q_2)$ at the $\bullet$-node $(i, j)$ we will use an additional averaging operator, namely

\begin{equation} \label{eq:www2}
\begin{aligned}
&|\mathcal{A}|_{i,j}^{\bullet}(\mathbf{q}) = \sqrt{ \left( \frac{q_1(i,j) + q_1(i-1,j)}{2}\right)^2 +  \left( \frac{q_2(i,j) + q_2(i,j-1)}{2}\right)^2 },\\
\end{aligned}
\end{equation}
where $q_1$ and $q_2$ are defined on $\circ$-nodes and $\square$-nodes, respectively. Similarly, the discrete divergence $\text{div}_{i,j}^{\bullet}(\bm{\mu})$ of $\bm{\mu} = (\mu_1, \mu_2)$ at the $\bullet$-node $(i, j)$ is defined by

\begin{equation} \label{eq:www3}
\begin{aligned}
&\text{div}_{i,j}^{\bullet}(\bm{\mu}) = [\mu_1(i,j) - \mu_1(i-1,j) + \mu_2(i,j) - \mu_2(i,j-1)]/h,\\
\end{aligned}
\end{equation}
where $\mu_1$ (resp. $\mu_2$) is defined at the $\circ$-nodes (resp. $\square$-nodes). Finally, we define shifting and identity operators by
\begin{equation} \label{eq:www4}
\begin{aligned}
&\mathcal{S}_{1}^{\pm}\varphi(i,j) = \varphi(i\pm1,j),  ~~\mathcal{S}_{2}^{\pm}\varphi(i,j) = \varphi(i,j\pm1)  ~~\text{and}~~\mathcal{I}\varphi(i,j) = \varphi(i,j).
\end{aligned}
\end{equation}

\subsection{Computation of the discrete analogue of  $\mathbf{p}^{n+1/3}$ in (\ref{eq:q_n13pro_solut})}\label{sec:dis_pn13}
Let us recall that from (\ref{eq:q_n13pro_solut}) one has

\begin{equation} \label{eq:dis_p_n13pro}
\mathbf{p}^{n+1/3} = \text{max}\bigg\{0, 1 - \frac{c}{|\mathbf{p}^{n}|}\bigg\} \mathbf{p}^{n},
\end{equation}
where $c = \tau a + \tau b |\nabla \cdot \bm{\lambda}^n|^2$. In the discrete setting, the first (resp. second) component of $\bf{p}^{n}$ and $\bm{\lambda}^{n}$ is defined at $\circ$-nodes (resp. $\square$-nodes), we need to discuss the two situations we will encounter when discretizing (\ref{eq:dis_p_n13pro}) (for simplicity, we will denote $\bm{\lambda}^n$ by $\bm{\lambda}$ and $\mathbf{p}^n$ by $\mathbf{p}$).

1) If $(i,j)$ is a $\circ$-node, the corresponding discretization of $\bf{p}$ and $c$ is given as follows:
\begin{equation}\label{eq:pp1}
p_{1}^{(1)}(i,j) = p_{1}(i,j); ~~ p_{2}^{(1)}(i,j) = \mathcal{A}_{i,j}^{\circ}(p_2),
\end{equation}

\footnotesize{
	\begin{equation}\label{eq:p1}
	\begin{aligned}
	c^{(1)}(i,j)& = \tau \left[ a + b |\partial_1 \lambda_{1}(i,j) + \partial_2 \lambda_{2}(i,j)|^2 \right] \\
	& = \tau \left[a + b \left| \frac{\lambda_1(i+1,j) - \lambda_1(i-1,j)}{2h} + \frac{\lambda_2(i+1,j) + \lambda_2(i,j)}{2h} - \frac{\lambda_2(i,j-1) + \lambda_2(i+1,j-1)}{2h}\right|^2\right].
	\end{aligned}
	\end{equation}}

\normalsize
2)If $(i,j)$ is a $\square$-node, the corresponding discretization of $\bf{p}$ and $c$ is given as follows:

\begin{equation}\label{eq:pp2}
p_{1}^{(2)}(i,j) = \mathcal{A}_{i,j}^{\square}(p_1); ~~ p_{2}^{(2)}(i,j) = p_{2}(i,j),
\end{equation}

\footnotesize{
	\begin{equation}\label{eq:p2}
	\begin{aligned}
	c^{(2)}(i,j)& = \tau \left[a +  b \left| \partial_1 \lambda_{1}(i,j) + \partial_2 \lambda_{2}(i,j)\right|^2 \right] \\
	& = \tau \left[a + b \left| \frac{\lambda_1(i,j) + \lambda_1(i,j+1)}{2h} - \frac{\lambda_1(i-1,j) + \lambda_2(i-1,j+1)}{2h} + \frac{\lambda_2(i,j+1) - \lambda_2(i,j-1)}{2h}\right|^2\right]. \\
	\end{aligned}
	\end{equation}}

\normalsize
Finally,

\begin{equation}\label{eq:final_p13}
p^{n+1/3}_{\alpha}(i,j) = \text{max}\left\{0, 1 - \frac{c^{(\alpha)}(i,j)}{\sqrt{|p_{1}^{(\alpha)}(i,j)|^2 + |p_{2}^{(\alpha)}(i,j)|^2}}\right\}p_{\alpha}^{(\alpha)}(i,j), ~~ \alpha=\{1, 2\}.
\end{equation}

\subsection{Computation of the discrete analogue of  $\lambda^{n+1/3}$ in (\ref{eq:u_n13pro_solut1})}\label{sec:dis_lambdan13}
We recall that (\ref{eq:u_n13pro_solut1}) reads as

\begin{equation} \label{eq:u_n13pro_solut_new}
\gamma \bm{\lambda}^{n+1/3} - \tau \nabla(2b|\mathbf{p}^{n+1/3}|\nabla\cdot{\bm{\lambda}}^{n+1/3}) = \gamma  \bm{\lambda}^{n}, \mbox{  in  }{\Omega},
\end{equation}
It is completed by periodic boundary conditions. For simplicity, we denote the (known) vector $(\mathbf{p}^{n+1/3}, \bm{\lambda}^n)$ by $(\widetilde{\mathbf{p}}, \widetilde{\bm{\lambda}})$ and $\bm{\lambda}^{n+1/3}$ (an unknown one) by $\bm{\lambda}$. Following \cite{tai2011}, we discretize (\ref{eq:u_n13pro_solut_new}) as follows

\begin{equation}\label{eq:final_lambda13}
\gamma \bm{\lambda} -\tau \nabla^{+} (2b|\widetilde{\mathbf{p}}|\text{div}^{-}\bm{\lambda}) =  \gamma \widetilde{\bm{\lambda}}.
\end{equation}
%where the discrete gradient and divergence operators $\nabla^{+}$ and $\text{div}^{-}$ have been defined in Section \ref{sec:usefulnote}.

To solve (\ref{eq:final_lambda13}), we will employ (as in \cite{tai2011}) \textit{a frozen coefficient} approach where instead of solving (\ref{eq:final_lambda13}) we solve

\begin{equation} \label{eq:final_lambda13_1}
\gamma \bm{\lambda} - c^* \nabla^{+} (\text{div}^{-}\bm{\lambda})  =  \gamma \widetilde{\bm{\lambda}}
- \nabla^{+} \left[(c^* - 2\tau b|\widetilde{\mathbf{p}}|)\text{div}^{-}\widetilde{\bm{\lambda}}\right],
\end{equation}
with $c^*$ properly chosen. Following \cite{tai2011}, we advocate taking $c^* = \max_{\bullet\text{-nodes}(i,j)} 2\tau b |\mathcal{A}|_{i,j}^{\bullet}(\widetilde{\mathbf{p}})$.

In matrix form, (\ref{eq:final_lambda13_1}) can be written as, in $\Omega_h$,
\begin{equation} \label{eq:final_lambda13_2}
\gamma h^2 \left(\begin{array}{c}
\lambda_{1}\\
\lambda_{2}
\end{array}
\right)
-c^* \left(\begin{array}{c}
\partial_{1}^{+}\\
\partial_{2}^{+}
\end{array}
\right)
\left(\begin{array}{cc}
\partial_{1}^{-} & \partial_{2}^{-}
\end{array}
\right)
\left(\begin{array}{c}
\lambda_{1}\\
\lambda_{2}
\end{array}
\right)
=
\gamma h^2 \left(\begin{array}{c}
\widetilde{\lambda}_{1}\\
\widetilde{\lambda}_{2}
\end{array}
\right)
-
\left(\begin{array}{c}
\partial_{1}^{+}\\
\partial_{2}^{+}
\end{array}
\right)
(c^*h - 2\tau bh|\widetilde{\mathbf{p}}|)\text{div}^{-}\widetilde{\bm{\lambda}},
\end{equation}
or, equivalently,
\begin{equation} \label{eq:final_lambda13_3}
\left\{
\begin{array}{*{20}c}
\left(\gamma h^2 -c^* \partial_{1}^{+}\partial_{1}^{-}\right)\lambda_1
- c^* \partial_{1}^{+}\partial_{2}^{-} \lambda_2
=
\gamma h^2 \widetilde{\lambda}_1
- \partial_{1}^{+}(c^*h - 2\tau b h|\widetilde{\mathbf{p}}|)\text{div}^{-}\widetilde{\bm{\lambda}},\\
-c^* \partial_{2}^{+}\partial_{1}^{-} \lambda_1
+ \left(\gamma h^2  - c^* \partial_{2}^{+}\partial_{2}^{-} \right)\lambda_2
=
\gamma h^2 \widetilde{\lambda}_2
- \partial_{2}^{+}(c^*h - 2\tau b h|\widetilde{\mathbf{p}}|)\text{div}^{-}\widetilde{\bm{\lambda}}.

\end{array}
\right.
\end{equation}

Using the shifting and identity operator defined in Section \ref{sec:usefulnote}, for each pair $(i, j)$ the first equation in (\ref{eq:final_lambda13_3}) reads as
\begin{equation} \label{eq:final_lambda13_4}
\left[\gamma h^2 \mathcal{I} + c^* (\mathcal{I} - \mathcal{S}_{1}^{+})(\mathcal{I} - \mathcal{S}_{1}^{-}) \right]\lambda_{1}(i,j) + c^* (\mathcal{I} - \mathcal{S}_{1}^{+})(\mathcal{I} - \mathcal{S}_{2}^{-})\lambda_2(i,j)
= g_1(i,j),
\end{equation}
where
\[
g_1(i,j) = \gamma h^2 \widetilde{\lambda}_1(i,j) - \left[\left(c^* h - 2\tau b h|\mathcal{A}|_{i+1,j}^{\bullet}(\widetilde{\mathbf{p}})\right)\text{div}_{i+1,j}^{\bullet}\widetilde{\bm{\lambda}} - \left(c^*h - 2\tau b h|\mathcal{A}|_{i,j}^{\bullet}(\widetilde{\mathbf{p}})\right)\text{div}_{i,j}^{\bullet}\widetilde{\bm{\lambda}}\right].
\]

Similarly, the second equation of (\ref{eq:final_lambda13_3}) reads as
\begin{equation} \label{eq:final_lambda13_5}
c^* (\mathcal{I} - \mathcal{S}_{2}^{+})(\mathcal{I} - \mathcal{S}_{1}^{-}) \lambda_1(i,j)
+ \left[\gamma h^2 \mathcal{I} + c^* (\mathcal{I} - \mathcal{S}_{2}^{+})(\mathcal{I} - \mathcal{S}_{2}^{-})\right]\lambda_2(i,j)
= g_2(i,j),
\end{equation}
where
\[g_2(i,j) = \gamma h^2 \widetilde{\lambda}_2(i,j) - \left[\left(c^*h - 2\tau b h|\mathcal{A}|_{i,j+1}^{\bullet}(\widetilde{\mathbf{p}})\right)\text{div}_{i,j+1}^{\bullet}(\widetilde{\bm{\lambda}}) - \left(c^*h - 2\tau b h|\mathcal{A}|_{i,j}^{\bullet}(\widetilde{\mathbf{p}})\right)\text{div}_{i,j}^{\bullet}\widetilde{\bm{\lambda}}\right].
\]

The boundary conditions we consider being the periodic ones, we may apply the \textit{discrete Fourier transform} $\mathcal{F}$ to equations (\ref{eq:final_lambda13_4}), (\ref{eq:final_lambda13_5}). We obtain then

\begin{equation} \label{eq:final_lambda13_6}
\left(\begin{array}{cc}
a_{11} & a_{12}\\
a_{21} & a_{22}
\end{array}
\right)
\mathcal{F}
\left(\begin{array}{c}
\lambda_{1}(y_i, y_j)  \\
\lambda_{2}(y_i, y_j)
\end{array}
\right)
=
\mathcal{F}
\left(\begin{array}{c}
g_{1}(y_i, y_j)  \\
g_{2}(y_i, y_j)
\end{array}
\right),
\end{equation}
where in (\ref{eq:final_lambda13_6}), one has:
\[
a_{11} = \gamma h^2 - 2c^*(\cos z_i -1), ~a_{12} =  c^*(\cos z_i -1 + \sqrt{-1}\sin z_i )(\cos z_j - 1 - \sqrt{-1}\sin z_j),
\]
\[
a_{21} = c^*(\cos z_j -1 + \sqrt{-1}\sin z_j)(\cos z_i - 1 - \sqrt{-1}\sin z_i), ~a_{22} = \gamma h^2 - 2c^*(\cos z_j -1),
\]
with
\begin{equation} \label{eq:wahaha}
z_i = \frac{2 \pi}{M_1}y_i, ~~y_i = 1, 2, \cdots, M_1, ~~\text{and}~~ z_j = \frac{2 \pi}{N_1}y_j, ~~y_j = 1, 2, \cdots, N_1.
\end{equation}

The determinant $D(i,j)$ of the coefficient matrix in (\ref{eq:final_lambda13_6}) is given by
\[
D(i,j) = \gamma^2 h^4 + 2\gamma h^2 c^*(2 - \cos z_i - \cos z_j),
\]
implying $D(i,j) > 0$ if $\gamma > 0$. It follows then from (\ref{eq:final_lambda13_6}), that (with obvious notation) the solution $\bm{\lambda}$ of problem (\ref{eq:final_lambda13_1}) (that is the discrete analogue of $\bm{\lambda}^{n+1/3}$ in (\ref{eq:u_n13pro_solut1})) is given by

\begin{equation} \label{eq:final_lambda13_62}
\left\{
\begin{aligned}
&\lambda_1 = Real\left[\mathcal{F}^{-1} \left( \frac{a_{22}\mathcal{F}(g_1) - a_{12}\mathcal{F}(g_2)}{D}\right)\right],\\
&\lambda_2 = Real\left[\mathcal{F}^{-1} \left( \frac{- a_{21}\mathcal{F}(g_1) + a_{11}\mathcal{F}(g_2)}{D}\right)\right],
\end{aligned}
\right.
\end{equation}
where $Real(x + \sqrt{-1}y) = x$, and $\bm{\lambda} = (\lambda_1, \lambda_2)$.

\subsection{Computation of the discrete analogue of  $(\mathbf{p}^{n+2/3}, \lambda^{n+2/3})$ in (\ref{eq:split5})}\label{sec:dis_plambdan23}
We need to solve problem (\ref{eq:q_n23pro}) to get the solutions. In the following, we give the details of its discretization.

\subsubsection{Solution of (\ref{eq:selectgamma_sol})}\label{sec:dis_plambdan23_s0}
From Section \ref{sect:S0}, we see that the minimizer of the functional in (\ref{eq:q_n23pro1}) over $\sigma_0$ is given by:
\begin{equation}\label{eq:dis_s0}
\left(\mathbf{p}_{0}^{n+2/3}(x), \bm{\lambda}_{0}^{n+2/3}(x)\right) = \left(\mathbf{0}, \frac{\bm{\lambda}^{n+1/3}(x)}{\text{max}[1, |\bm{\lambda}^{n+1/3}(x)|]}\right).
\end{equation}
The discrete analogue of (\ref{eq:dis_s0}) reads as

\begin{equation}\label{eq:dis_s0_result}
\left(\mathbf{p}_{0}^{n+2/3}(i, j), \bm{\lambda}_{0}^{n+2/3}(i, j)\right) = \left(\mathbf{0}, \frac{\bm{\lambda}^{n+1/3}(i, j)}{\text{max}\left[1, \sqrt{|\lambda_{1}^{n+1/3}(i, j)|^2 + |\lambda_{2}^{n+1/3}(i, j)|^2}\right]}\right),
\end{equation}
with $\bm{\lambda}^{n+1/3}(i, j) = \left(\lambda_{1}^{n+1/3}(i, j), \lambda_{2}^{n+1/3}(i, j)\right)$.

\subsubsection{Discretization of problem (\ref{eq:sigma1})}\label{sec:dis_plambdan23_s1}
Section \ref{sect:S1} was dedicated to the solution of problem  (\ref{eq:sigma1}), a constrained minimization problem in $\mathbf{R}^4$ defined by

\begin{equation}\label{eq:sigma1_dis}
\inf_{(\mathbf{q},\bm{\mu})\in \mathbf{R}^2\times \mathbf{R}^2, \mathbf{q} \neq \mathbf{0}, \mathbf{q}\cdot \bm{\mu}=|\mathbf{q}|, |\bm{\mu}| = 1} \left[|\mathbf{q} - \mathbf{p}^{n+1/3}(x)|^2 + \gamma|\bm{\mu} - \bm{\lambda}^{n+1/3}(x)|^2 \right].
\end{equation}

Proceeding as in Section  \ref{sect:S1}, we define $\mathbf{x}_{i, j}$ and $\mathbf{y}_{i, j}$ by
\[
\mathbf{x}_{i,j} = \left(\mathbf{x}_{i,j}^{(1)}, \mathbf{x}_{i,j}^{(2)}\right) = \left(p_{1}^{n+1/3}(i,j),  p_{2}^{n+1/3}(i,j)\right),
\]
\[\mathbf{y}_{i,j} = \left(\mathbf{y}_{i,j}^{(1)}, \mathbf{y}_{i,j}^{(2)}\right) = \left(\lambda_{1}^{n+1/3}(i,j),  \lambda_{2}^{n+1/3}(i,j)\right).
\]
Then, we use the following discrete variant of algorithm (\ref{eq:fixpoint}) to compute $\theta_{i,j}^{*}$
\begin{equation} \label{eq:fixpoint_dis}
\left\{
\begin{aligned}
&\theta_{i, j}^{(0)} = |\mathbf{x}_{i, j}|\\
&\text{For}~ k\geq 0, ~\theta_{i, j}^{(k)} \rightarrow \theta_{i, j}^{(k+1)} ~\text{as~follows} \\
&\theta_{i, j}^{(k+1)} = \max\left[0, \frac{\mathbf{x}_{i, j}\cdot(\theta_{i,j}^{(k)}\mathbf{x}_{i,j}+\gamma\mathbf{y}_{i,j})}{|\theta_{i,j}^{(k)}\mathbf{x}_{i,j}+\gamma\mathbf{y}_{i,j}|})\right].
\end{aligned}
\right.
\end{equation}
Once $\theta_{i,j}^{*}$ is computed we obtain the discrete analogues of $\left(\mathbf{p}_{1}^{n+2/3}(x), \bm{\lambda}_{1}^{n+2/3}(x)\right)$ from the following formula which is the discrete analogue of  (\ref{eq:lambda23_1}), (\ref{eq:p23_1}):
\begin{equation} \label{eq:thetapro_final}
\left\{
\begin{aligned}
&\bm{\lambda}^{n+2/3}(i,j) = \frac{\theta_{i,j}^{*} \mathbf{x}_{i,j} + \gamma \mathbf{y}_{i,j}}{\sqrt{\left|\theta_{i,j}^{*} x_{i,j}^{(1)} + \gamma y_{i,j}^{(1)}\right|^{2} + \left|\theta_{i,j}^{*} x_{i,j}^{(2)} + \gamma y_{i,j}^{(2)}\right|^2}},\\
&\mathbf{p}^{n+2/3}(i,j) = \theta_{i,j}^{*} \bm{\lambda}^{n+2/3}(i,j).
\end{aligned}
\right.
\end{equation}

\subsection{Discretization of problem (\ref{eq:p33_5}) }\label{sec:dis_plambdan33}
From Section \ref{sec:3.6}, we have $\bm{\lambda}^{n+1/3} = \bm{\lambda}^{n+2/3}$ and $\mathbf{p}^{n+1} = \nabla u^{n+1}$ where $u^{n+1}$ is the solution of the following linear elliptic problem

\begin{equation} \label{eq:1upro_solut_1}
- \nabla^{2} u^{n+1} + \tau u^{n+1} = -\nabla\cdot \mathbf{p}^{n+2/3} + \tau f, ~~\text{in}~{\Omega},
\end{equation}
completed by periodic boundary conditions. We need to discretize this problem.  Denoting $\mathbf{p}^{n+2/3}$ by $\widetilde{\mathbf{p}}$, we employ the following finite difference scheme to approximate (\ref{eq:1upro_solut_1}):

\begin{equation} \label{eq:1upro_solut_2}
\left(\begin{array}{cc}
\partial_{1}^{-} & \partial_{2}^{-}
\end{array}
\right)
\left[
\left(\begin{array}{c}
\partial_{1}^{+}\\
\partial_{2}^{+}
\end{array}
\right)u^{n+1}
-
h\left(\begin{array}{c}
\widetilde{p}_{1}\\
\widetilde{p}_{2}
\end{array}
\right)
\right]
+
\tau h^2
\left(f - u^{n+1}
\right)
= 0, ~ \text{in}~~ \Omega_h
\end{equation}

Problem (\ref{eq:1upro_solut_2}) is equivalent to
\begin{equation} \label{eq:1upro_solut_3}
\left(\partial_{1}^{-}\partial_{1}^{+} + \partial_{2}^{-}\partial_{2}^{+} - \tau h^2\right)u^{n+1}
=
h(\partial_{1}^{-}\widetilde{p}_{1} + \partial_{2}^{-}\widetilde{p}_{2}) - \tau h^2 f,
\end{equation}

Relation (\ref{eq:1upro_solut_3}) can be written also as
\begin{equation} \label{eq:1upro_solut_4}
\left[(\mathcal{I} - \mathcal{S}_{1}^{-})(\mathcal{S}_{1}^{+} - \mathcal{I}) + (\mathcal{I} - \mathcal{S}_{2}^{-})(\mathcal{S}_{2}^{+} - \mathcal{I}) - \tau h^2\mathcal{I} \right]u^{n+1}(i,j)
=
g(i,j),
\end{equation}
where $g(i,j) = h(\partial_{1}^{-}\widetilde{p}_{1}(i,j) + \partial_{2}^{-}\widetilde{p}_{2}(i,j)) - \tau h^2 f(i,j)$. From the  periodicity of the boundary conditions, it makes sense to use FFT to solve problem (\ref{eq:1upro_solut_4}). We obtain then

\begin{equation} \label{eq:1upro_solut_5}
w_{i,j}\mathcal{F}(u^{n+1}(i,j))
=
\mathcal{F}(g(i,j)),
\end{equation}
where $w(i,j) = \left[(1 - e^{-\sqrt{-1}z_i})(e^{\sqrt{-1}z_i} - 1) + (1 - e^{-\sqrt{-1}z_j})(e^{\sqrt{-1}z_j} - 1) - \tau h^2 \right]$, with $z_i$ and $z_j$ as in (\ref{eq:wahaha}). From (\ref{eq:1upro_solut_5}), we obtain (with obvious notation)

\begin{equation} \label{eq:1upro_solut_6}
u^{n+1}
=
Real\left[\mathcal{F}^{-1} \left( \frac{\mathcal{F}(g)}{w}\right)\right],
\end{equation}
with $Real(\cdot)$ as in Section \ref{sec:dis_lambdan13}. Once $u^{n+1}$ is known we compute $\mathbf{p}^{n+1}$ by

\begin{equation} \label{eq:pn1solution}
\mathbf{p}^{n+1}=\nabla^{+} u^{n+1} =
\left(\begin{array}{c}
\partial_{1}^{+}u^{n+1}\\
\partial_{2}^{+}u^{n+1}
\end{array}
\right),
\end{equation}
(operators have been defined in Section \ref{sec:usefulnote}). Finally, the discrete analogue of $\bm{\lambda}^{n+1}(x)$, for a.e. $x\in \Omega$, is given by
\begin{equation} \label{eq:thetapro666}
\left\{
\begin{aligned}
&\lambda_{1}^{n+1}(i,j) = \lambda_{1}^{n+2/3}(i,j),\\
&\lambda_{2}^{n+1}(i,j) = \lambda_{2}^{n+2/3}(i,j).
\end{aligned}
\right.
\end{equation}

\subsection{Further comments}\label{sec:4.7}
In Sections \ref{sec:dis_pn13} to \ref{sec:dis_plambdan33}, we have supplied the details for the discretization for the sub-problems associated with the operator-splitting scheme (\ref{eq:initial2})-(\ref{eq:split6}). In Section \ref{sec:results}, we will apply the above methodology to the solution of image smoothing problems. It will allow us to demonstrate that with our approach, one can handle the elastica energy functional efficiently and accurately. In addition, we will use further experiments to show the good properties of the proposed method.  These include modularity,  good stability and the low cost of the algorithm.

\section{Numerical Results}\label{sec:results}
In this section, the proposed method is applied to image smoothing to test its effectiveness. All experiments are implemented in MATLAB(R2016a) on a laptop of 8Gb RAM and Intel(R) Core(TM) i7-7500 CPU: @2.70 GHz 2.90GHz. Note that the intensities of all images are in the range of $[0, 1]$. For simplicity, we also use mesh size $h=1$.

In our experiments, it is reasonable to stop the iteration if the following defined relative error (ReErr) of the solution is smaller than the predefined tolerance $tol$, i.e.,
\begin{equation} \label{eq:rte}
\begin{aligned}
\text{ReErr} = \frac{\|u^{n+1} - u^{n}\|_2}{\|u^{n+1}\|_2} < tol,
\end{aligned}
\end{equation}
where $tol$ is a pre-defined positive value. In particular, a bigger $tol$ may result in a faster stopping of the proposed iterative method.

One of the main advantage of the new method is that it only involves the time step $\tau$ as free algorithm parameter to be chosen.
The fast speed and robustness of the proposed method are also verified in this section by some specially designed experiments.

In what follows, we apply, in Section \ref{sec:denoise}, the proposed method to image smoothing. Then, in Section \ref{sec:4.3}, we compare the speed of convergence and stability properties of this method with those of the THC algorithm (\cite{tai2011}). In Section \ref{sec:4.3} we further discuss various aspects of the new method and draw some conclusions concerning its ability at solving smoothing problems.\\

%\begin{remark}
 \textbf{Remark 5.1}:	\textit{In some earlier works, c.f. \cite{tai2011,duanwang2013,zhang2017}, the Euler elastica model was applied to image denoising. We found, however, that  ``edge-preserving smoothing'' describes better than ``denoising'' the properties of the proposed method. Indeed, minimizing the Euler elastica energy functional is actually a way to enforce the curvature of an image to be small, a property leading to the smoothing of image details in non-edge regions, while preserving and smoothing the edges. The ``denoising'' effect is just an intermediate result, ``smoothing'' being actually the final result of the elastica energy functional minimization. Therefore, in this article, we will use “smoothing” instead of ``denoising'', a departure from the terminology we used in previous works. }
%\end{remark}

\subsection{Image smoothing}\label{sec:denoise}
In this section, we first apply (in Section  \ref{sec:4.1.1}) the proposed method to the ROF  model (i.e., $b = 0$) and then show, in Section \ref{sec:4.1.2}, some results of image smoothing with the Euler elastica model.

\subsubsection{The proposed method for the ROF model}\label{sec:4.1.1}

\begin{figure*}[!htp]
	\begin{center}
		
		\includegraphics[width=1.9in,height=1.7in]{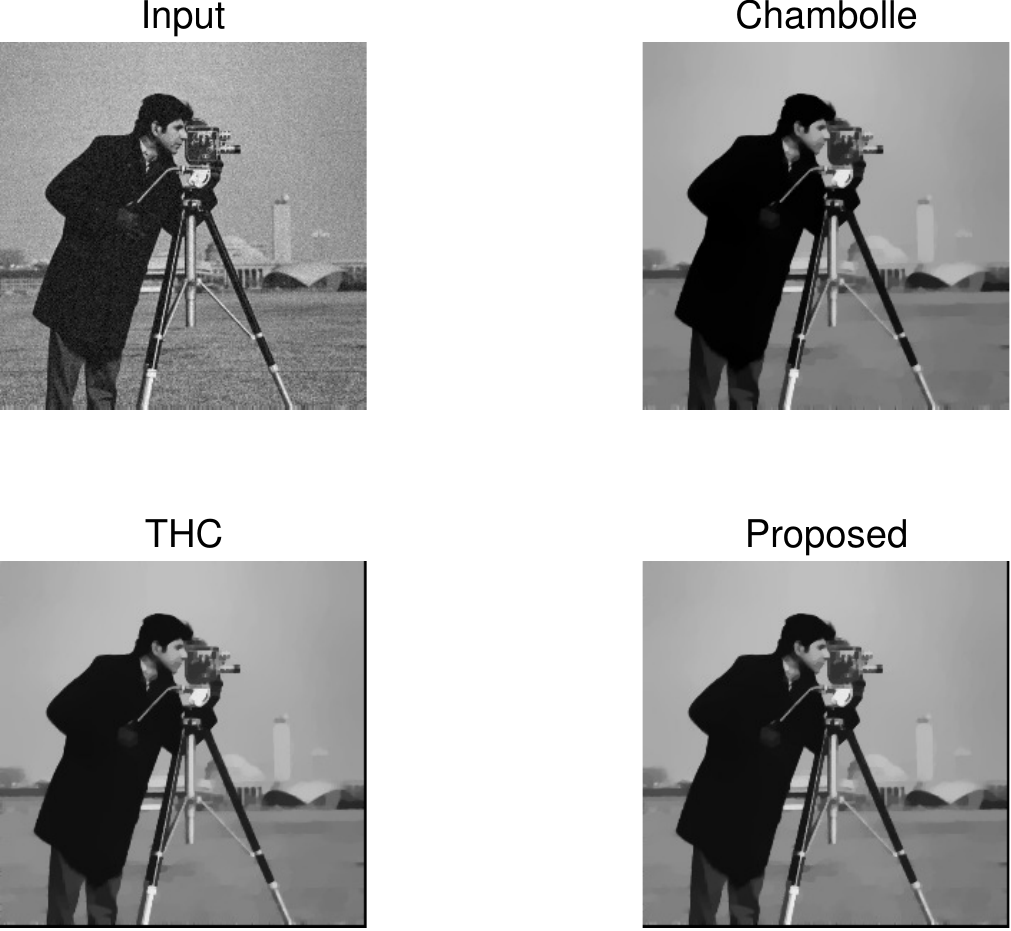}
		\includegraphics[width=2.1in,height=1.7in]{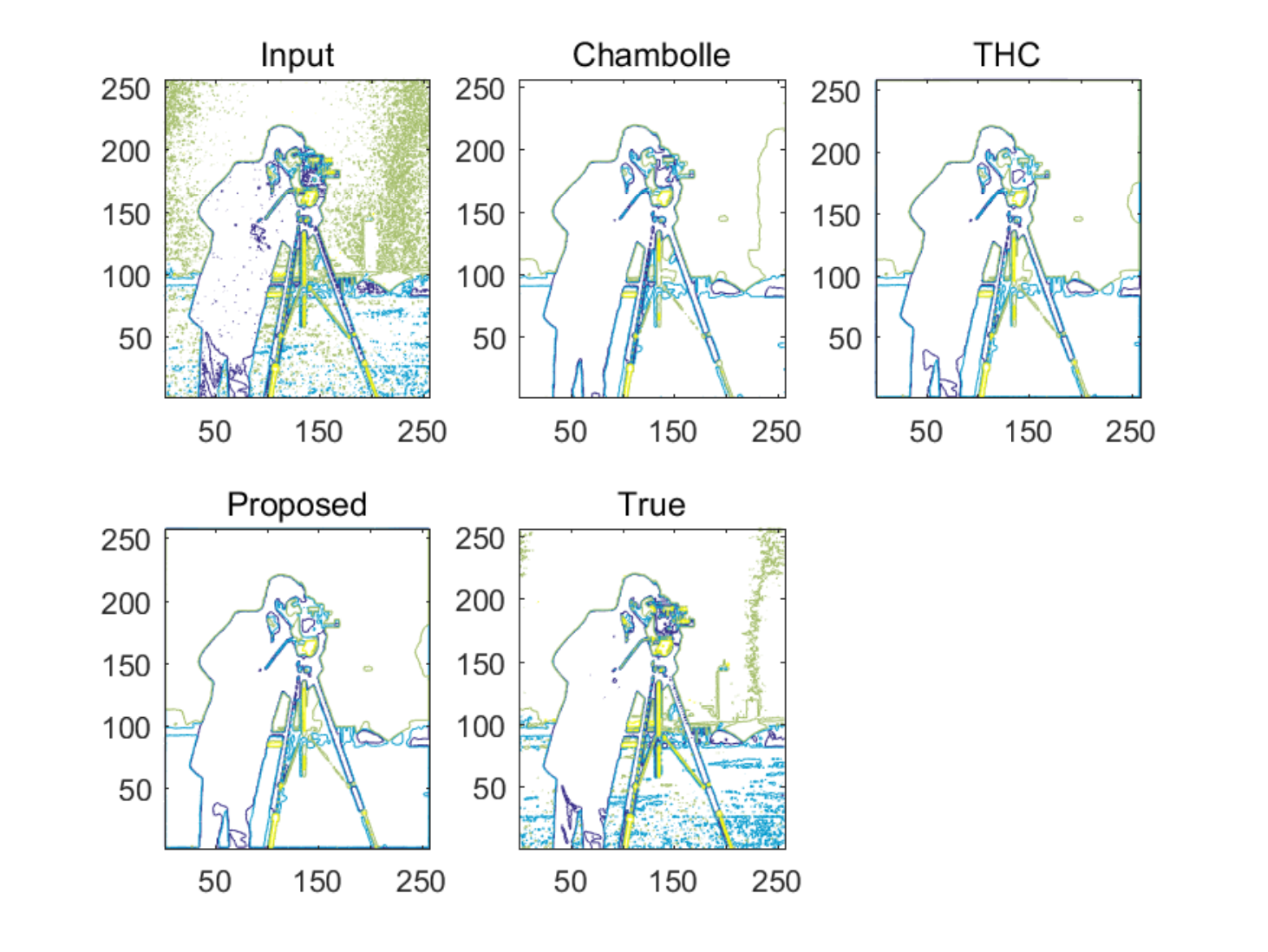}
		\includegraphics[width=1.9in,height=1.7in]{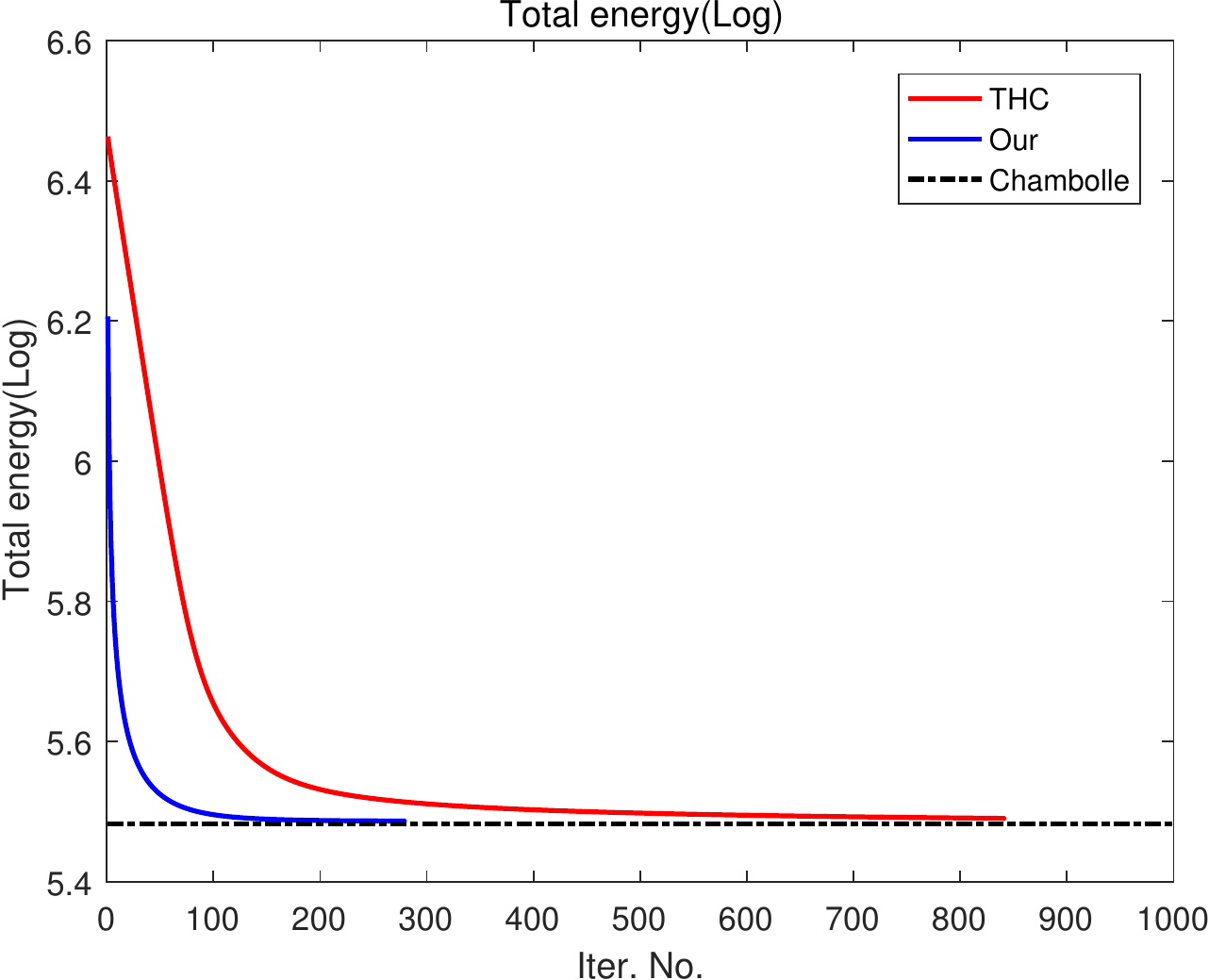}
		
		(a)
		
		\includegraphics[width=1.9in,height=1.7in]{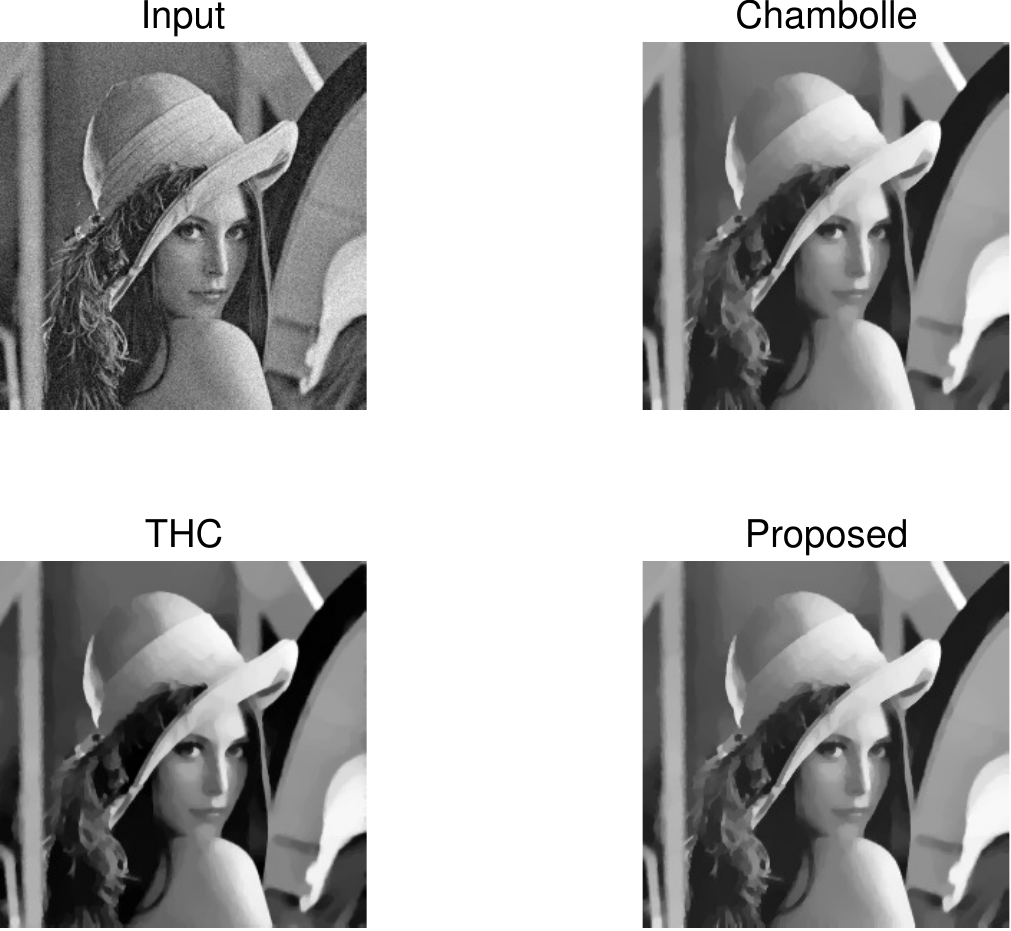}
		\includegraphics[width=2.1in,height=1.7in]{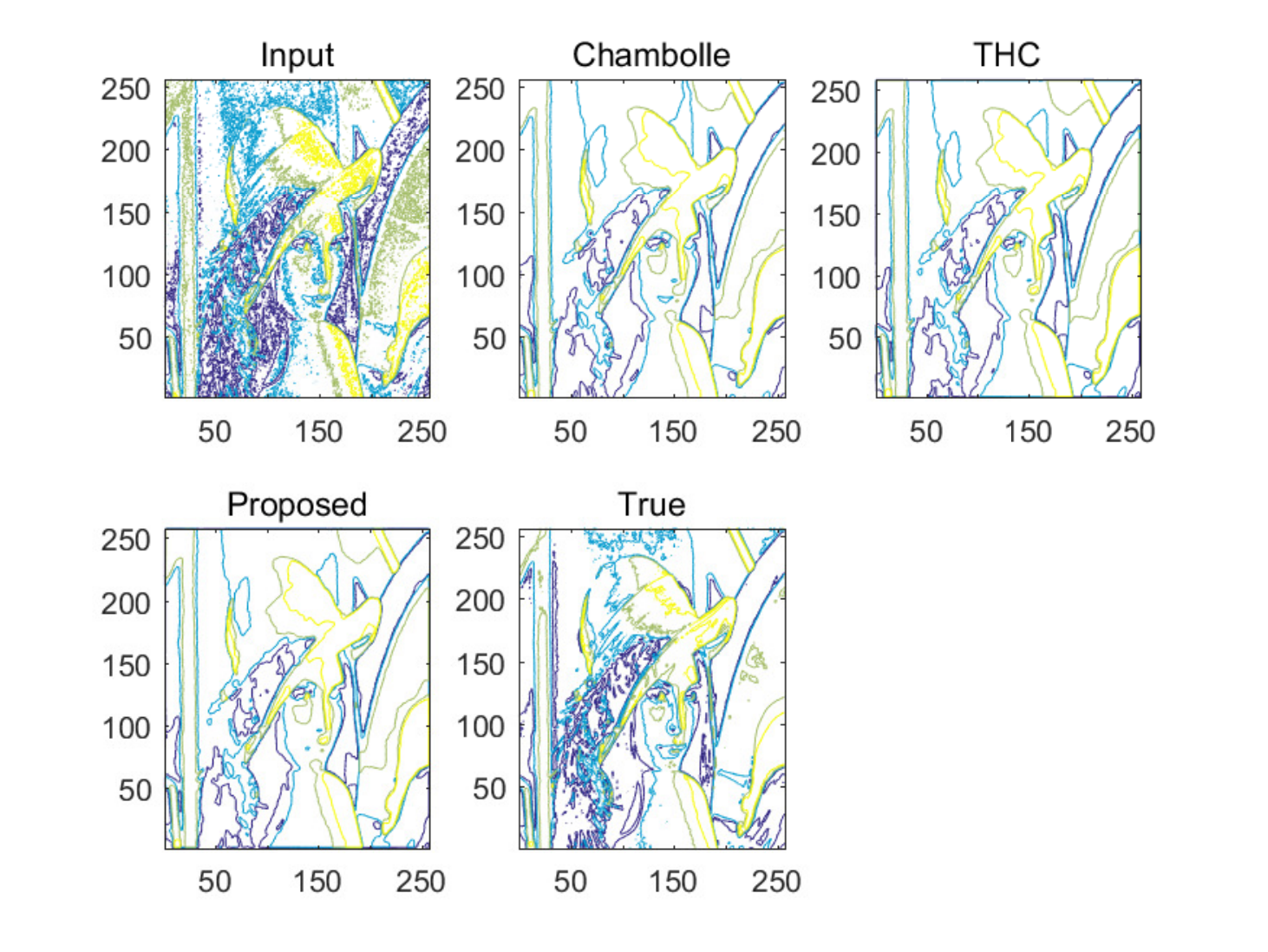}
		\includegraphics[width=1.9in,height=1.7in]{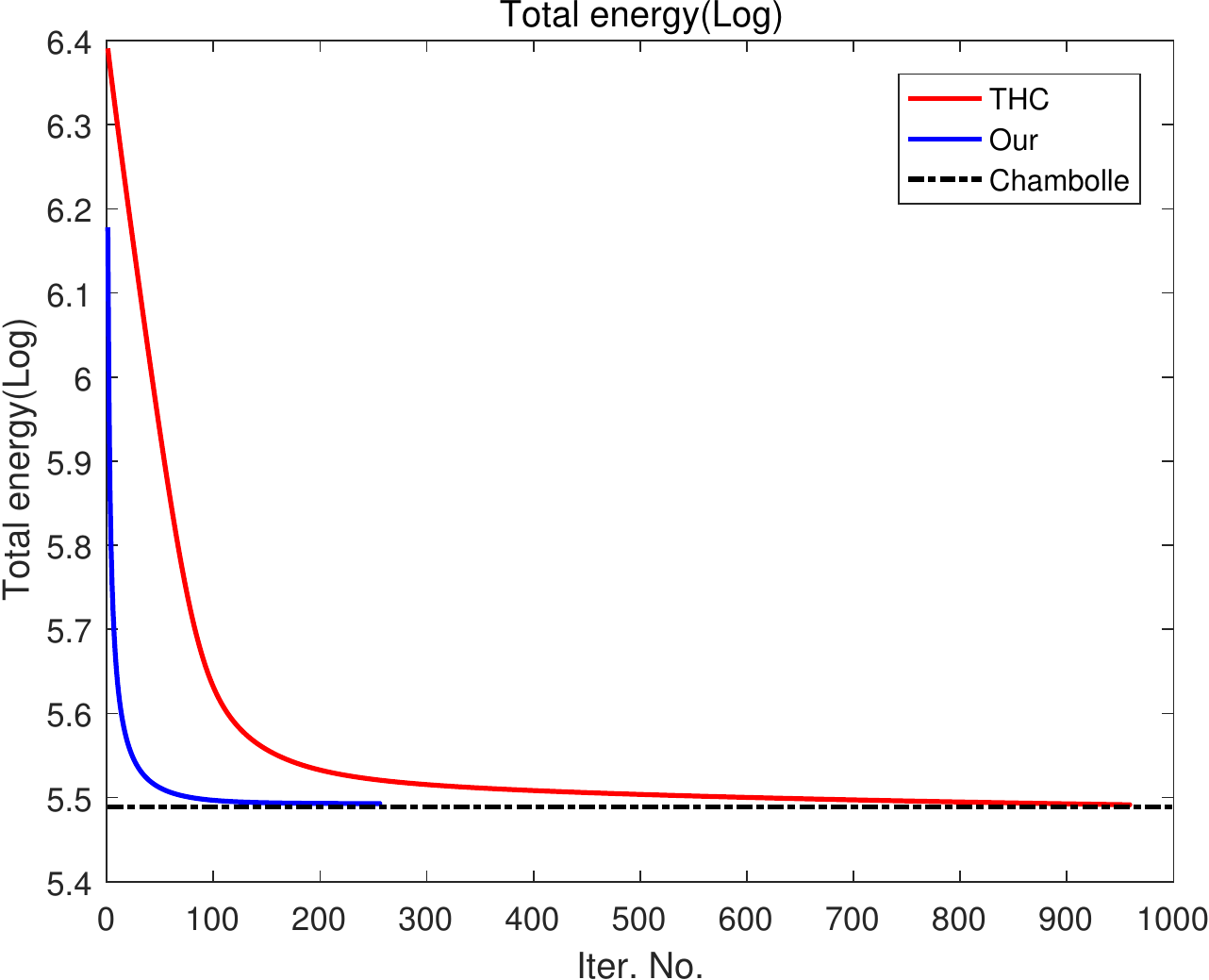}
		
		(b)
		
	\end{center}
	\caption{We use Chambolle's algorithm, the THC method and the proposed method to apply the ROF model when setting $b = 0$ to see if all methods converge to the same solution and same energy. Note that we set $a = 0.1$, and $tol = 1\times 10^{-5}$ for all methods, and $\tau = 0.1$ for the proposed method. The experiments are tested on the image ``cameraman'' (a) and ``Lena'' (b). The results (the first two columns), the contour plots (the middle column), and energy plots (the right column) by the three compared methods are shown. Note that the black dashed lines in the energy plots represent the final energy of the ROF model solved by the Chambolle's method. From this figure, we observe that the results of three compared approaches for the ROF model converge finally to the same energy value, which is a preliminary validation of the proposed method. }
	\label{fig:img_denoise}
\end{figure*}

We apply Chambolle's method \cite{chambolle2004}, the THC method \cite{tai2011} and the proposed method to the ROF model which is actually a special case of the Euler elastica energy when setting $b = 0$ in (\ref{eq:elastica}).
In Fig. \ref{fig:img_denoise}, we set $b = 0$ and fix $a = 0.1$ for the Euler elastica energy based image restoration problem (\ref{eq:elastica}), which is just the ROF model.  The results of Chambolle's method, of the THC method and of the proposed method for the ROF model are shown in Fig. \ref{fig:img_denoise}. In particular, we implemented our method with $\tau = 0.1$ and $\gamma^{n} = \text{max}\left(|\mathbf{p}^{n+1/3}|, \sqrt{\tau}\right)$, c.f.  (\ref{eq:gamma}). All three algorithms are solving the same ROF based problem and their energy converges to the same value. The restored images are also shown in Fig. \ref{fig:img_denoise}. We use this example to show that our algorithm works also for the ROF model.

\subsubsection{Application of the proposed method to image smoothing }\label{sec:4.1.2}

In what follows, we show the capability of the new method at image smoothing. In addition, we also demonstrate the superiority of the Euler elastica model when compared with the ROF model.

We report the results of image smoothing by the Euler elastica model solved by the proposed method, and by the ROF model solved by Chambolle's method \cite{chambolle2004} as well. The results demonstrate the competitive ability of edge-preserving image smoothing of the Euler elastica model.

\begin{figure*}
	\begin{center}
		
		\includegraphics[width=1.5in,height=1.5in]{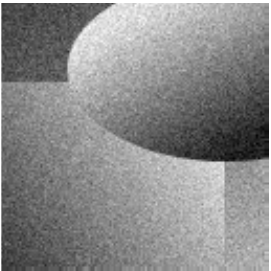}
		\includegraphics[width=1.5in,height=1.5in]{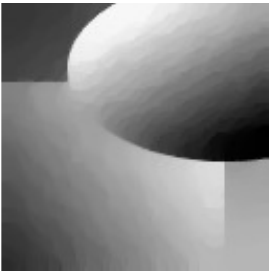}
		\includegraphics[width=1.5in,height=1.5in]{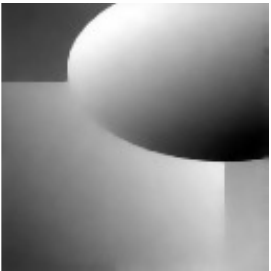}

		\includegraphics[width=1.5in,height=1.5in]{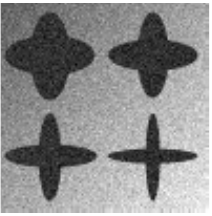}
		\includegraphics[width=1.5in,height=1.5in]{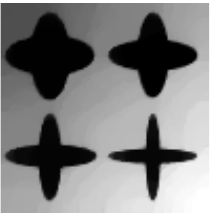}
		\includegraphics[width=1.5in,height=1.5in]{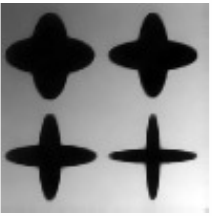}

		\includegraphics[width=1.5in,height=1.5in]{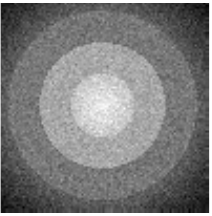}
		\includegraphics[width=1.5in,height=1.5in]{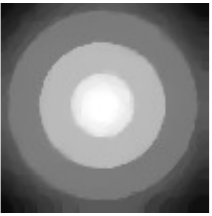}
		\includegraphics[width=1.5in,height=1.5in]{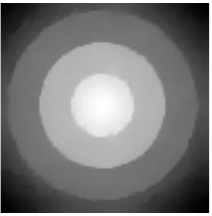}
		
		\includegraphics[width=1.5in,height=1.5in]{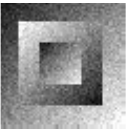}
		\includegraphics[width=1.5in,height=1.5in]{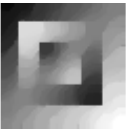}
		\includegraphics[width=1.5in,height=1.5in]{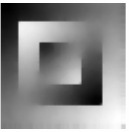}
		
		\includegraphics[width=1.5in,height=1.5in]{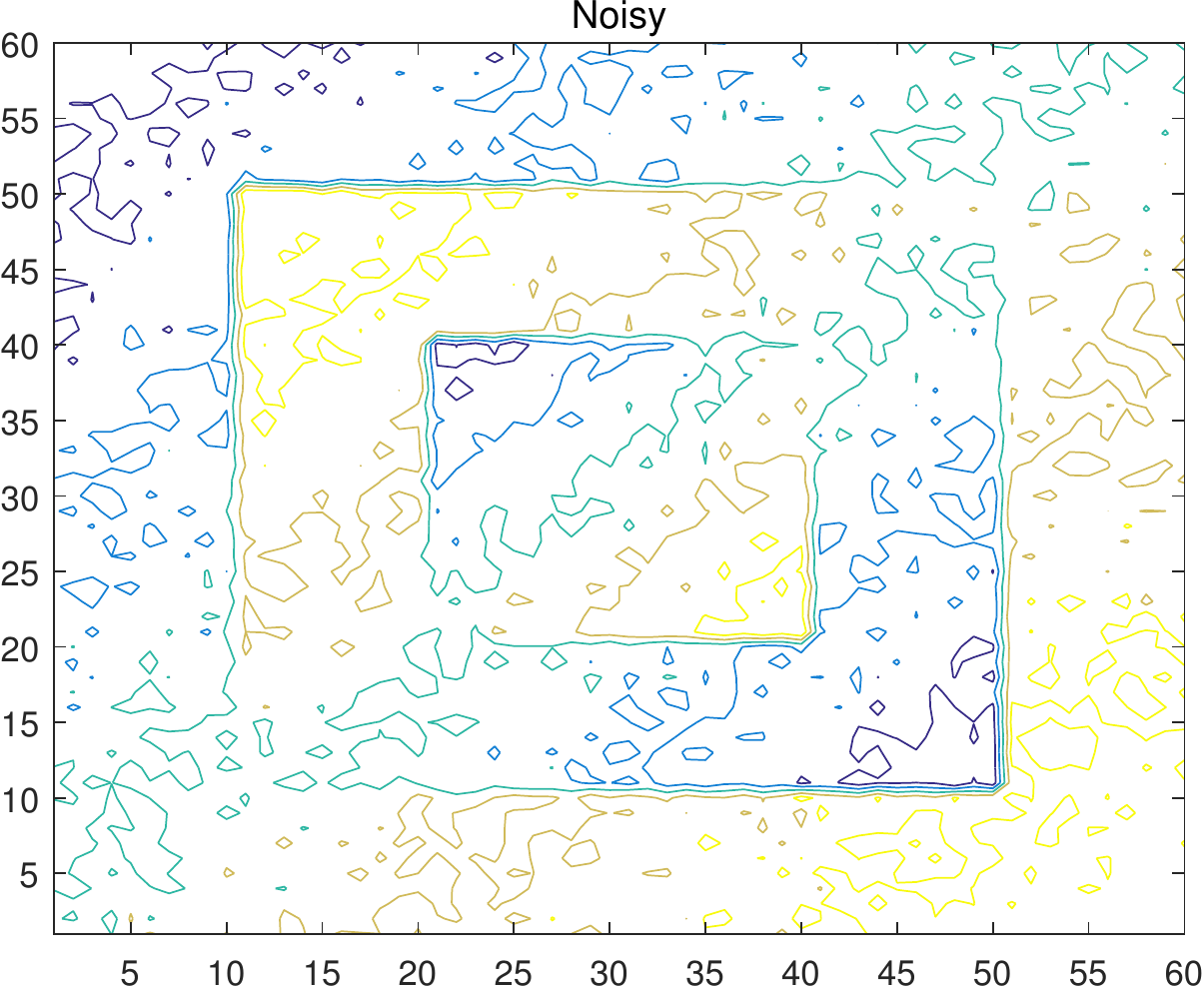}
		\includegraphics[width=1.5in,height=1.5in]{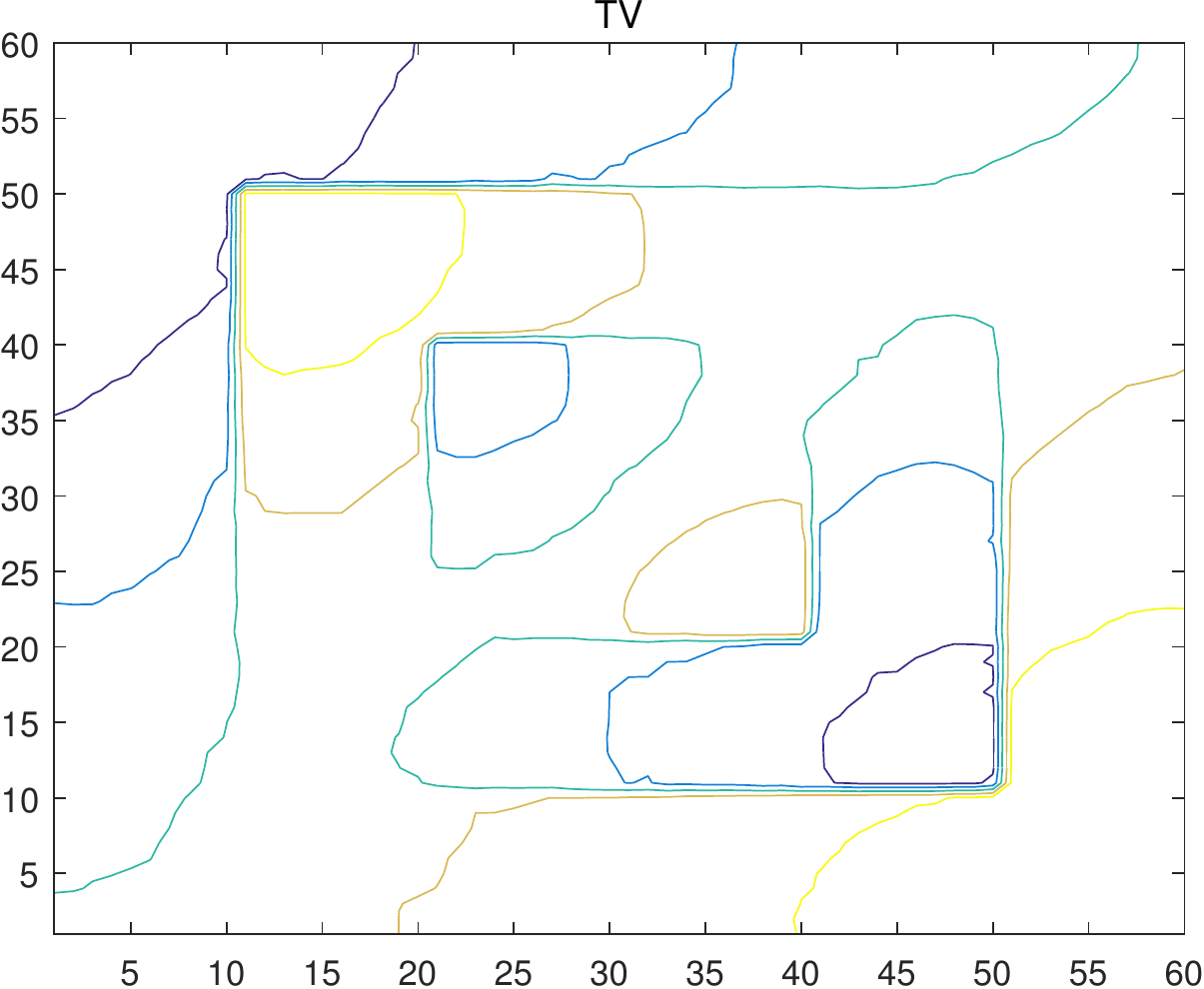}
		\includegraphics[width=1.5in,height=1.5in]{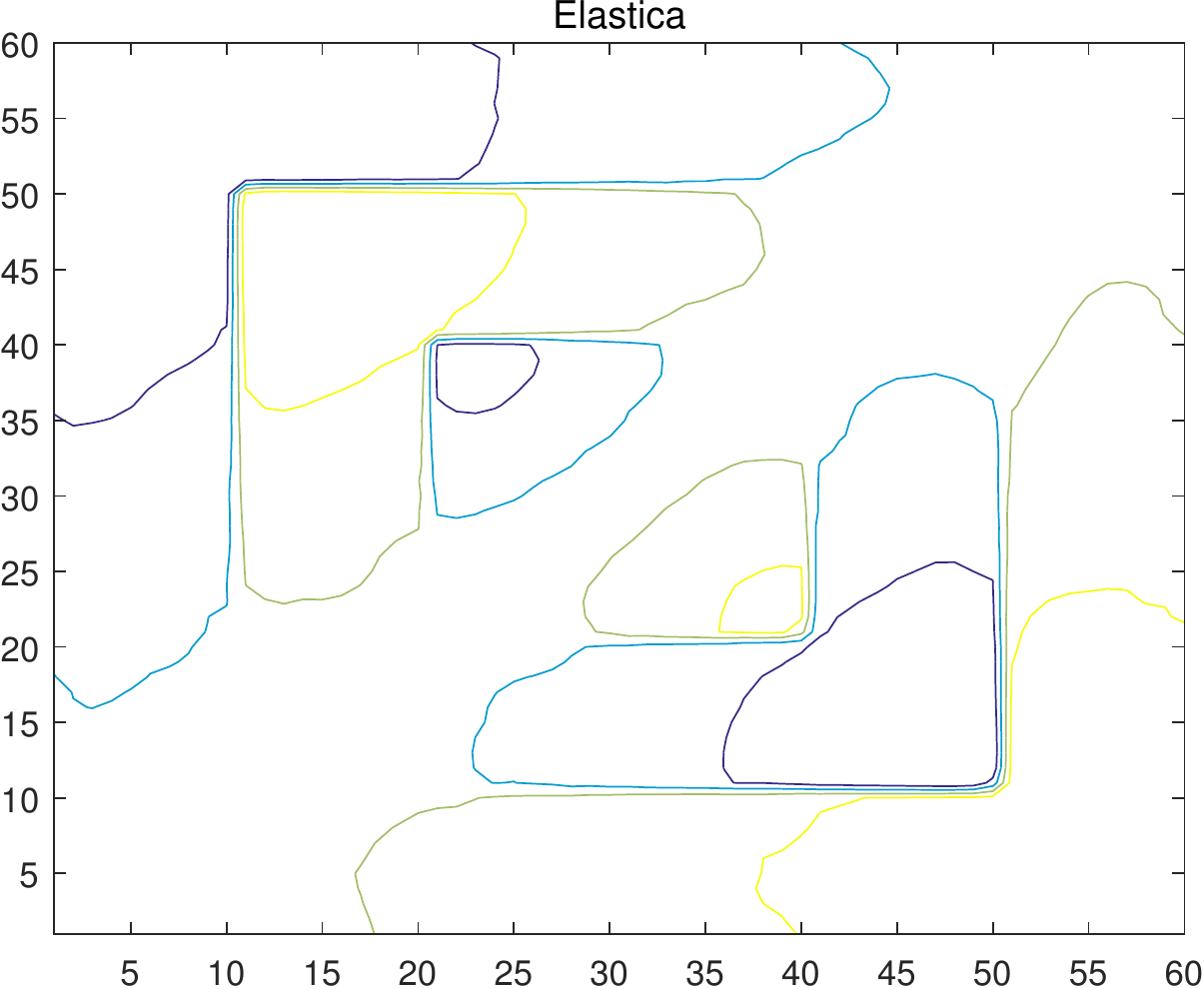} 			
		
		Noisy Image \hspace{0.8in} TV  \hspace{1in} Elastica
		
	\end{center}
	\caption{Image smoothing using the ROF and Euler elastica models. Left: Noisy images; Middle: ROF model treated by Chambolle's method; Right: the Euler elastica model treated by the proposed algorithm. The bottom row shows the contour of the last image. We can see that the elastica model gives images with smoother contours than the ROF model. From the figure, we see that the ROF model creates undesired staircase effects, while the elastica model overcomes it. Note that the parameters in our method for these four test examples are all set as $a = b = 0.1$, $\tau = 0.1$ and $\gamma^{n} = \text{max}\left(|\mathbf{p}^{n+1/3}|, \sqrt{\tau}\right)$. }
	\label{fig:img_elas_denoise}
\end{figure*}

Fig. \ref{fig:img_elas_denoise} shows the results of the proposed algorithm for Euler's elastica model and Chambolle's algorithm for ROF model on four synthetic images. The noisy images are shown in the left column, and the smoothed images by ROF model and the Euler elastica model are shown in the middle and right columns, respectively. Gaussian white noise with zero mean and a 20 standard deviation is used for the first three images, i.e., ``ball'', ``star'' and ``circle'', a 10 standard deviation being used for the fourth image, i.e., ``square''. We acknowledge that all test images in this figure are taken from \cite{tai2011}.

From Fig. \ref{fig:img_elas_denoise}, the ROF model is able to well preserve image discontinuous jumps, e.g., sharp edges, but it leads to some undesired artifacts, for example, the staircase effect in the smooth regions. The Euler elastica model applied via our method not only well preserves the jumps, but also removes the noise  without leading to undesired artifacts in the smooth regions.
In the last row of Fig. \ref{fig:img_elas_denoise}, we have visualized the contours of the image ``square'' (noisy on the right, after ROF smoothing in the center, after elastica smoothing on the right). The smoothest contours are the ones obtained by the elastica model via our method. An analysis of these properties can be found in \cite{zhu2012}.

Note that for all the experiments reported in Fig. \ref{fig:img_elas_denoise}, all the involved parameters have the same values, i.e.,   $a = b = 0.1$, $\tau = 0.1$, $\gamma^{n} = \text{max}\left(|\mathbf{p}^{n+1/3}|, \sqrt{\tau}\right)$, and $tol = 1\times 10^{-5}$. This shows that our method is stable with respect to parameter choice, this property being one of its main advantages.

\begin{figure*}[!htp]
	\begin{center}
		
		\includegraphics[width=5.6in,height=4in]{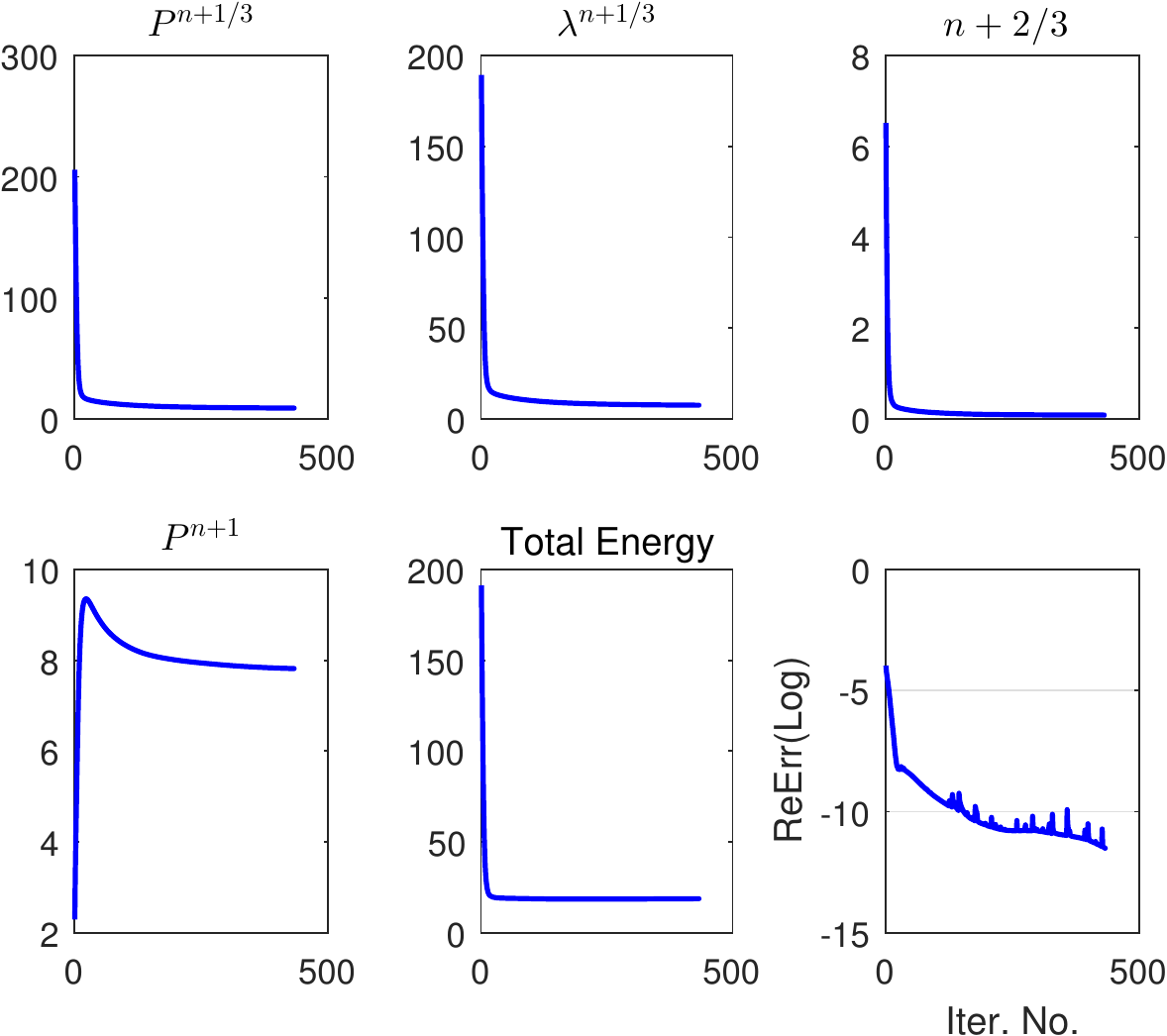}
		
	\end{center}
	\caption{We take the image ``square'' which is the last one in Fig. \ref{fig:img_elas_denoise} as an example to show the energy changes of each subproblem, i.e., the $\mathbf{p}^{n+1/3}$-subproblem (\ref{eq:q_n13pro}), the $\bm{\lambda}^{n+1/3}$-subproblem (\ref{eq:lambda_n13pro}), the $(n+2/3)$-subproblem (\ref{eq:q_n23pro}), the $\mathbf{p}^{n+1}$-subproblem (\ref{eq:p33_1}) and the total energy (\ref{eq:transelastica22}).   }
	\label{fig:energy_elas}
\end{figure*}

In Fig. \ref{fig:energy_elas}, we also monitor the energy changes of the subproblems and of the original problem (\ref{eq:transelastica22}). From this figure, it is clear that the energies of the $\mathbf{p}^{n+1/3}$ subproblem (\ref{eq:q_n13pro}), of the $\bm{\lambda}^{n+1/3}$ subproblem (\ref{eq:lambda_n13pro}), of $(n+2/3)$ subproblem (\ref{eq:q_n23pro}) (including $\mathbf{p}^{n+2/3}$ and $\bm{\lambda}^{n+2/3}$ subproblems) all decrease as $n$ increases, while the energy of the $\mathbf{p}^{n+1}$ subproblem (\ref{eq:p33_1}) increases to a stable value. This is because the role of the $(n +1/3)$ subproblem is to minimize the value of the elastica energy term, without taking the fidelity term into account.  Nevertheless, the total energy of the original problem (\ref{eq:transelastica22}) always decreases as $n$ increases.

%%%%%%%%%%%%%%%%%%%%%%%%%energy plot%%%%%%%%%%%%%%%%%%%%

\subsection{Advantages of the proposed method}\label{sec:4.3}	
In Section \ref{sec:denoise}, we applied our method to the Euler elastica model.
In what follows, some special experiments will be designed and implemented to illustrate the superiority of the proposed method compared with the THC method \cite{tai2011}.

The THC method proposed in \cite{tai2011} is an efficient approach to solve the Euler's elastica problem.
As shown by the results reported in \cite{tai2011}, the THC method can solve the Euler elastica problem hundredfold times faster than the Chan-Kang-Shen (CKS) method in \cite{chankang2002}. After the THC method, some promising approaches, see e.g., \cite{duanwang2013,duan2014,zhang2017,bae2017}, were proposed for the Euler elastica problem. In \cite{duanwang2013}, Duan et al. proposed another fast augmented Lagrangian method to solve the Euler elastica problem based on the framework of the THC method. Afterwards, Duan et al. in \cite{duan2014} applied the THC based method in to solve the Euler elastica regularized Mumford-Shah problem, aiming to deal with two-stage image segmentation. Also, Zhu et al. \cite{zhu2013} applied the THC method for the Euler elastica regularized Chan-Vese problem, which gets excellent segmentation results. In \cite{zhang2017}, Zhang et al. proposed a fast linearized augmented Lagrangian approach to solve the Euler elastica problem and applied it to image denoising.

However, the ALM method has some limitations. First, it needs three Lagrange multipliers and three augmentation parameters. For practical applications, it is rather difficult to tune these parameters. We have observed, as shown later in this section, that ALM has a fast convergence and produces very good results when these parameters are chosen correctly. However, if we just change these parameters slightly from their ``correct'' values, the algorithm will slow down dramatically. Moreover, these parameters are often image dependent and need to be chosen properly for different images.

The method proposed in this article is a simple and new operator splitting approach. It requires only the solution of few simple subproblems, moreover it requires the tuning of only one parameter, namely the time-discretization step $\tau$.

In what follows, we will design some numerical experiments to verify the above-mentioned advantages. In particular, the first advantage, i.e., fewer parameters, holds obviously. Thus we only need to verify the second and third advantages.

\subsubsection{Parameter sensibility: a numerical testing approach}\label{sec:4.3.1}

In order to assess the stability properties, with respect to parameter variations, of the method we introduced in this article, we will proceed as follow. First, we will fix the model parameters $a$ and $b$. Next, we will tune the augmentation parameters of the augmented Lagrangian in the THC method and the time-discretization step $\tau$ of our method. Finally, we will compare the results obtained by both methods.

\begin{figure*}[!t]
	\begin{center}
		\includegraphics[width=1.3in,height=1.3in]{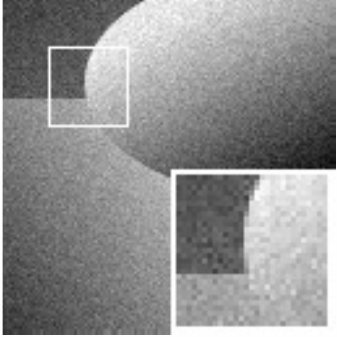}
		\includegraphics[width=1.3in,height=1.3in]{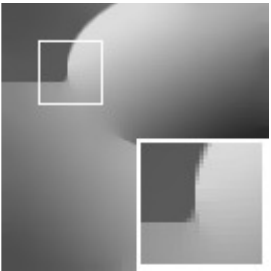}
		\includegraphics[width=1.3in,height=1.3in]{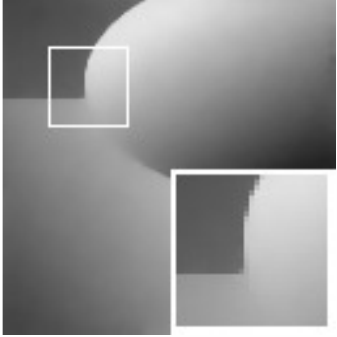}
		\includegraphics[width=1.7in,height=1.3in]{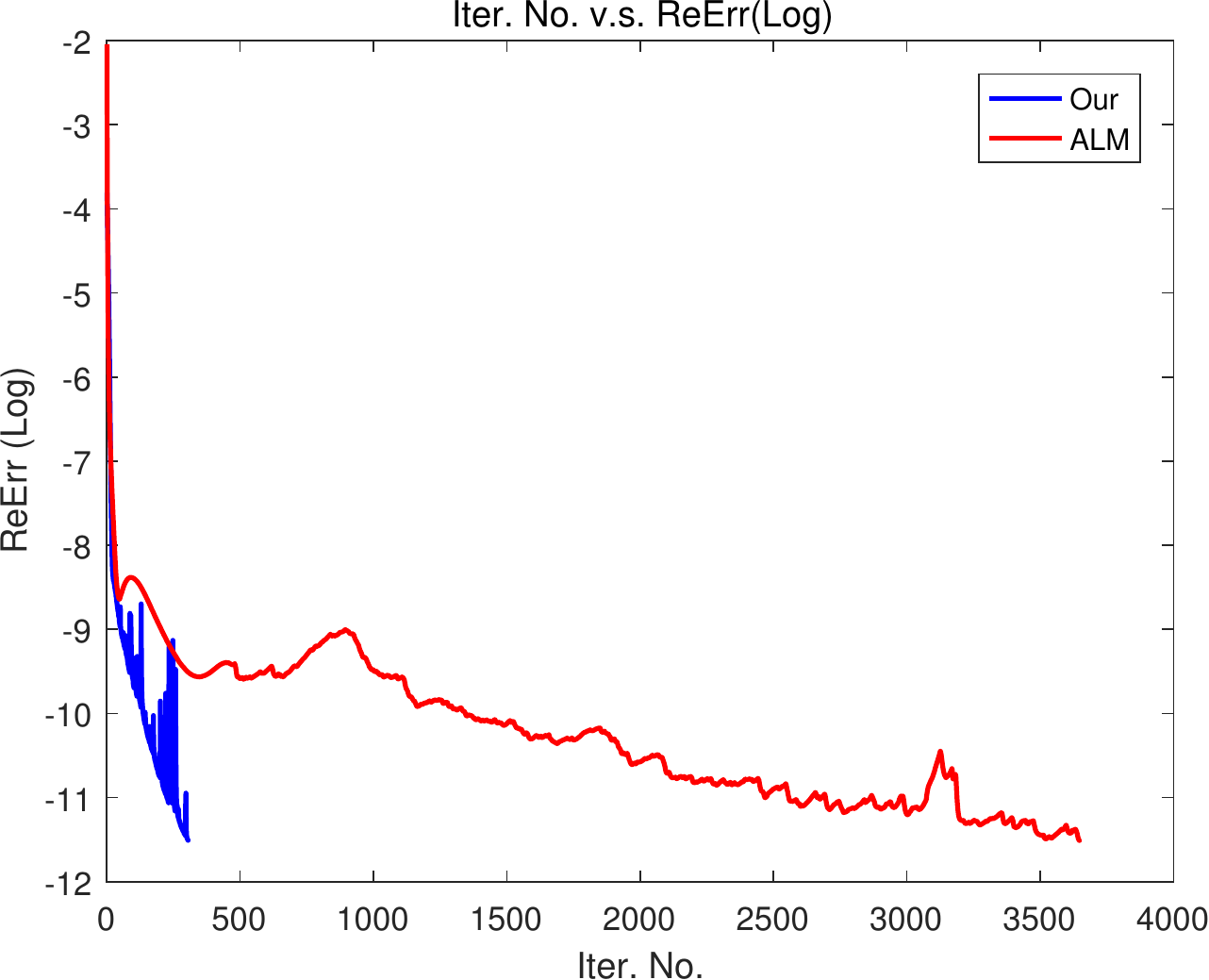}
		
		\includegraphics[width=1.3in,height=1.3in]{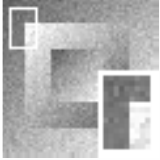}
		\includegraphics[width=1.3in,height=1.3in]{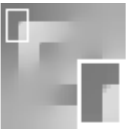}
		\includegraphics[width=1.3in,height=1.3in]{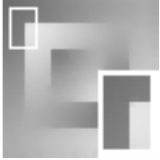}
		\includegraphics[width=1.7in,height=1.3in]{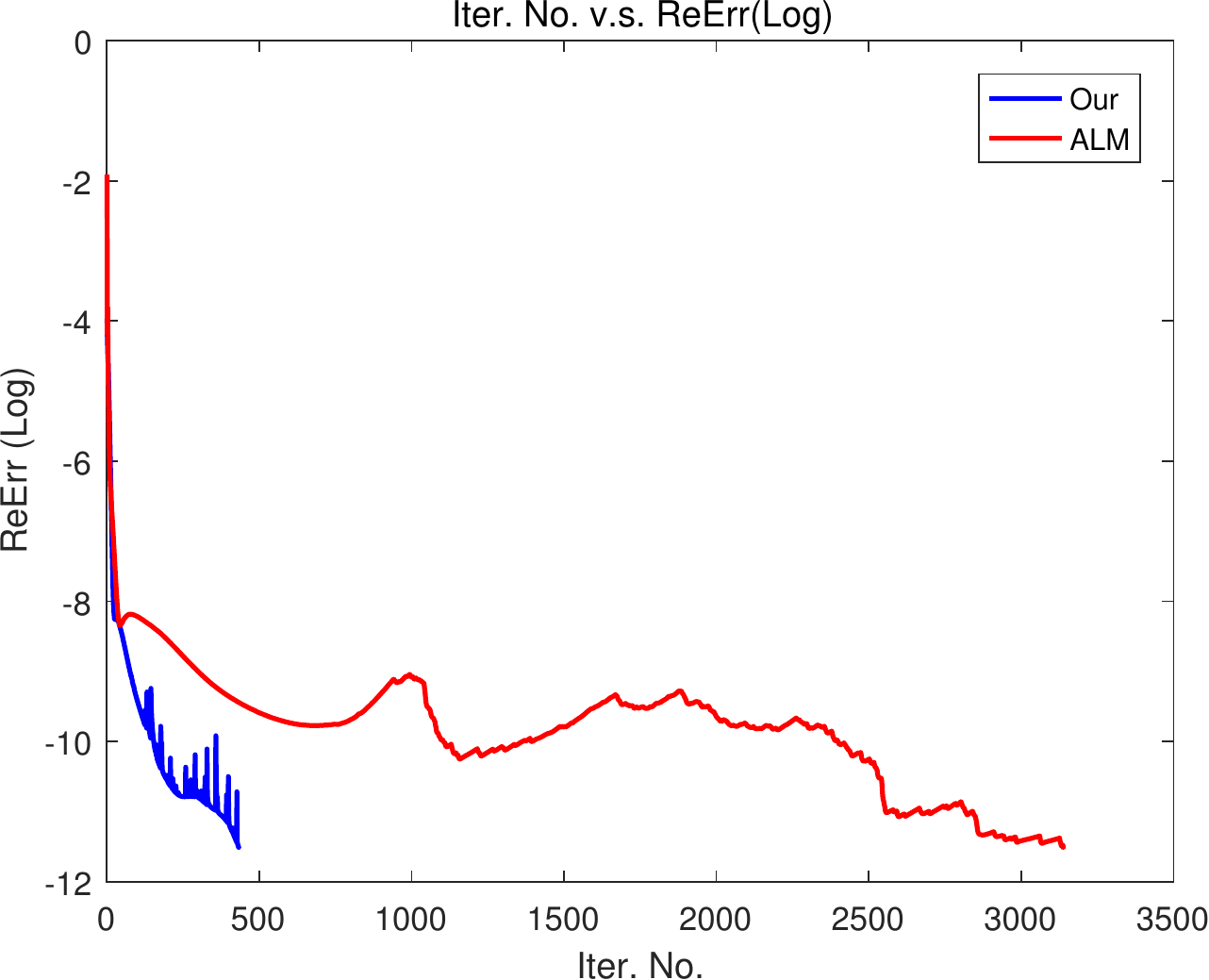}
		
		\includegraphics[width=1.3in,height=1.3in]{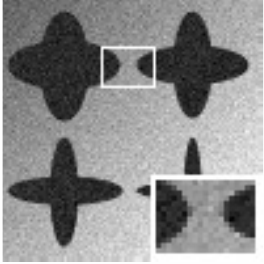}
		\includegraphics[width=1.3in,height=1.3in]{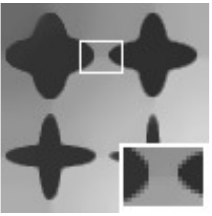}
		\includegraphics[width=1.3in,height=1.3in]{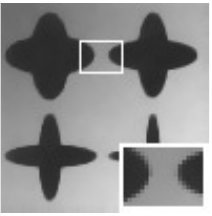}
		\includegraphics[width=1.7in,height=1.3in]{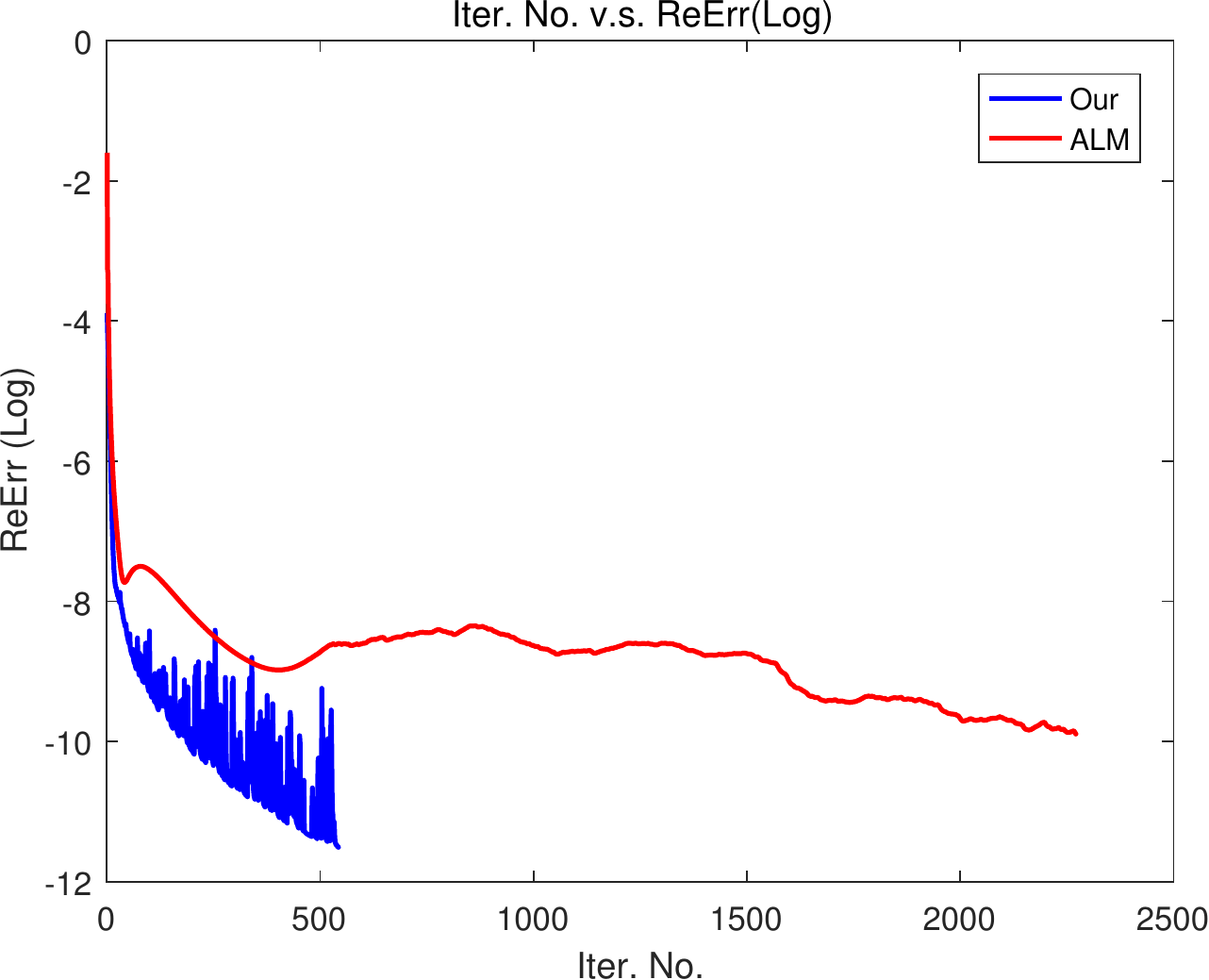}

		\includegraphics[width=1.3in,height=1.3in]{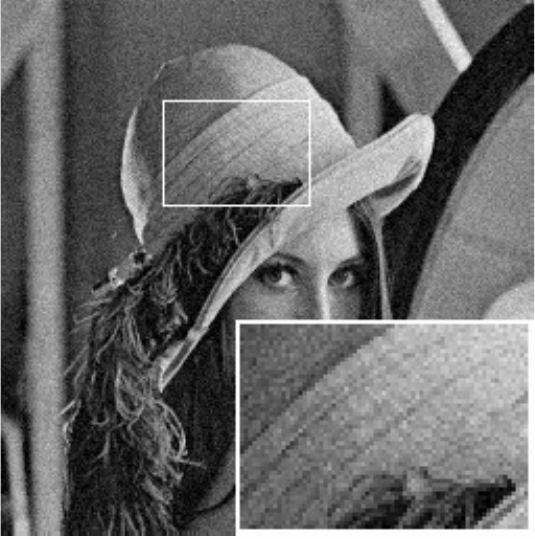}
		\includegraphics[width=1.3in,height=1.3in]{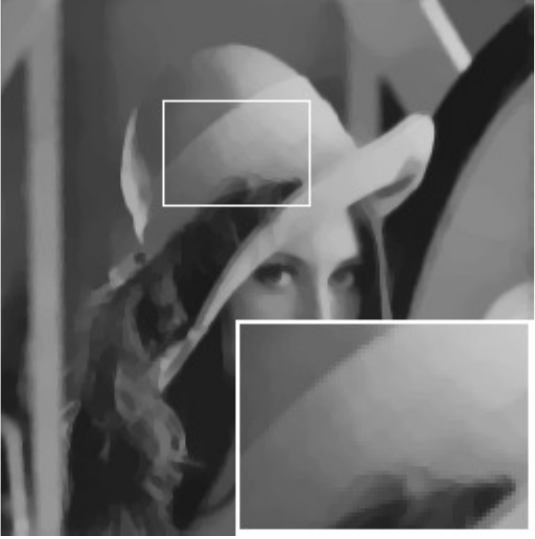}
		\includegraphics[width=1.3in,height=1.3in]{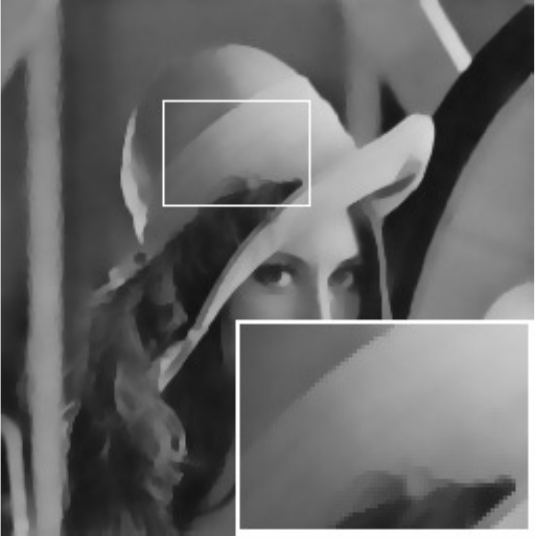}
		\includegraphics[width=1.7in,height=1.3in]{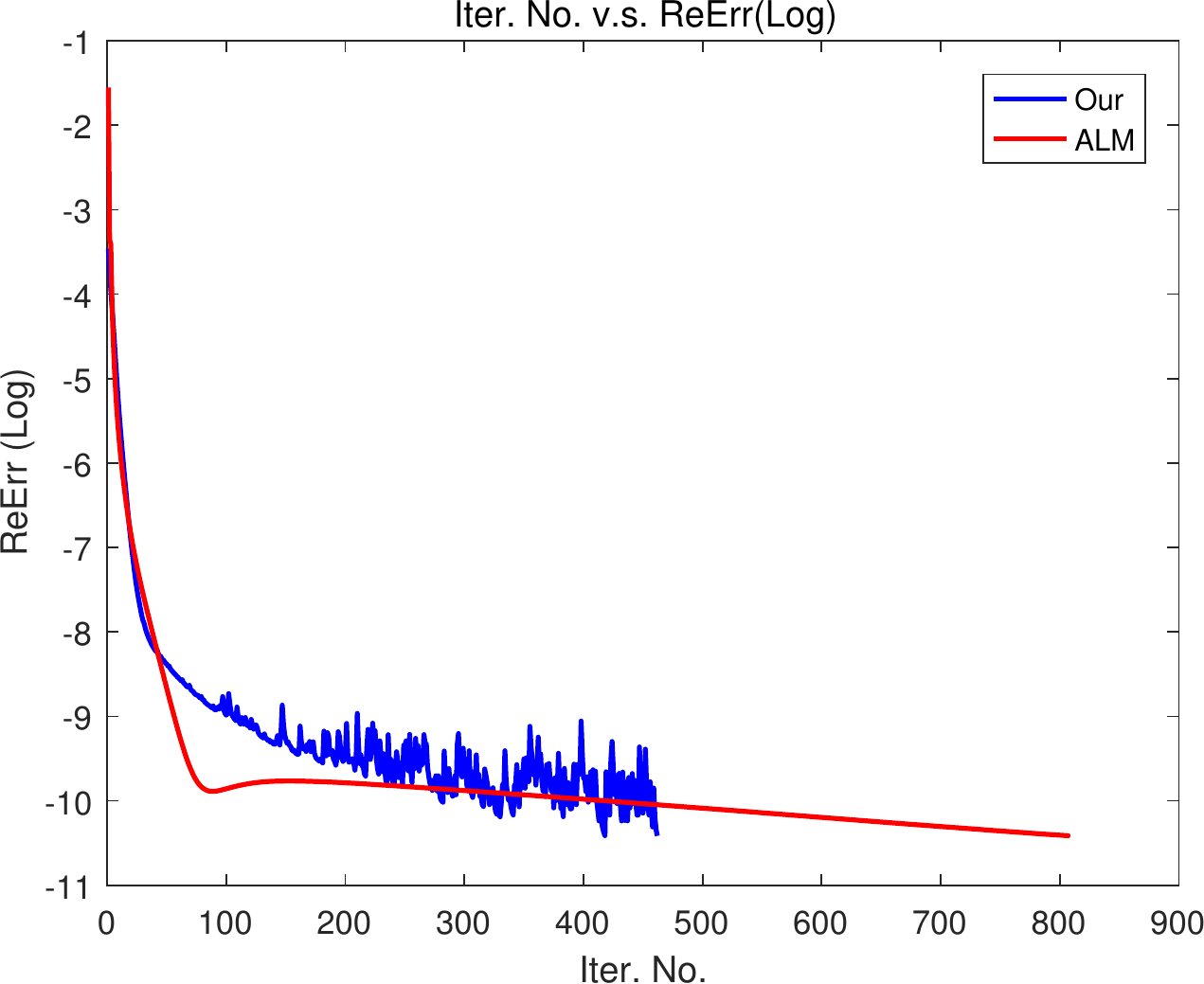}
		
		Noisy \hspace{0.8in} THC \hspace{0.8in} Proposed \hspace{1in} ReErr
		
	\end{center}
	\caption{First row: The figure shows that THC method and the method proposed in this article obtain both good results for $r_1 = 0.01$, $r_2 = 10$, $r_3 = 100$ and $\tau = 0.1$, respectively. For both methods we took $a = b = 0.1$ and $tol = 1\times 10^{-5}$ (resp., $3\times 10^{-5}$) for the first three examples (resp., for ``Lena'' our fourth example). The column on the right shows that the proposed method is significantly faster than the THC one (although its convergence is more oscillatory). }
	\label{fig:test1}
\end{figure*}

Our intention with the first experiment we performed was testing the sensitivity to one set of parameters for multiple images. For instance, $a$ and $b$ being fixed, we selected one specific image, then tuned the THC method (resp., the proposed method) augmentation parameters (resp., time-discretization step $\tau$) in order to obtain high quality image smoothing. Then, leaving the augmentation parameters and $\tau$ unchanged, we applied both methods to the smoothing of other images to see if one still obtains good results. In Fig. \ref{fig:test1} we have reported the results of the experiment described hereafter: (i) One considers four noisy images, namely ``ball'' ($128\times 128$), ``square'' ($60\times 60$), ``star'' ($100\times 100$) and ``Lena'' ($256\times 256$). (ii) We take $a = b = 0.1$ for both methods, to insure that they solve the same problem. (iii) The tolerance of the stopping criterion is set at $tol = 1\times 10^{-5}$ (resp. $3\times 10^{-5}$) for ``ball'', ``square'' and ``star'' (resp., ``Lena''). (iv) Taking ``ball'' as image of reference, we selected $r_1 = 0.01$, $r_2 = 10$ and $r_3 = 100$ for the THC method (resp., $\tau = 0.1$ and $\gamma$ given by (\ref{eq:gamma}) for the proposed method), these values producing high quality smoothing of the noisy ``ball'' image (see \cite{tai2011} for details about the THC method augmentation parameters $r_1$, $r_2$, $r_3$). (v) Finally, keeping the same values for the above parameters, we applied both methods to the other three images. The results reported in Fig. \ref{fig:test1} show that the method we propose in this article is still operational, unlike the THC method that leads to undesired image artifacts. Moreover, the right column of Fig. \ref{fig:test1} and Tab. \ref{tab:test1} show that the new method requires significantly less iterations than the THC one to verify the stopping criterion. Besides, the average computational time per iteration of the proposed method is also smaller than the one of the THC method, a possible explanation being that the method we propose in this article has fewer sub-problems, these sub-problems either having closed form solutions or being solvable by fast algorithms such as FFT.

\begin{table}[!t]
	\caption{In this table we have compared the following performances of the method introduced in this article and of the THC method: Third column: Number of iterations necessary to achieve convergence; Fourth column: Total computational time (in seconds); Fith column: Averaged computational time per iteration (in seconds). }\label{tab:test1}
	\scriptsize
	\begin{center}
		\begin{tabular}{|c|c|c|c|c|}
			\hline
			Image & Method & Iterations  & Time (s)      & Average Time (s)/per iteration        \\
			\hline\hline
			\multirow{2}{*}{\textbf{ball ($128\times 128$)}}&{\textbf{Proposed}}&{306}&{2.57}&{0.008}\\
			
			\cline{2-5}
			& \textbf{THC} &  3648 & 37.19  &  0.010   \\
			
			\hline
			\multirow{2}{*}{\textbf{square ($60\times 60$)}}&{\textbf{Proposed}}&{434}&{1.16}&{0.002}\\
			
			\cline{2-5}
			& \textbf{THC} &  3339 & 10.27  &  0.003   \\
			
			\hline
			\multirow{2}{*}{\textbf{star($100\times 100$)}}&{\textbf{Proposed}}&{562}&{3.70}&{0.006}\\
			
			\cline{2-5}
			& \textbf{THC} & 2234 & 17.58  &  0.007   \\
			
			\hline
			\multirow{2}{*}{\textbf{Lena ($256\times 256$)}}&{\textbf{Proposed}}&{462}&{15.21}&{0.033}\\
			
			\cline{2-5}
			& \textbf{THC} &  808 & 31.18  &  0.039   \\
			
			\hline
		\end{tabular}
	\end{center}
\end{table}

%%%%%%%%%%%%%%%%%%%%%% TEST12 %%%%%%%%%%%%%%%%%%%%%%%%%%%%%%%%%%%%%%%%%%%%%%%%%%%

\begin{figure*}
	\begin{center}
		\includegraphics[width=1.2in,height=1.2in]{test1_square_alm-eps-converted-to.pdf}
		\includegraphics[width=1.2in,height=1.2in]{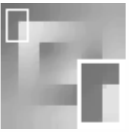}
		\includegraphics[width=1.2in,height=1.2in]{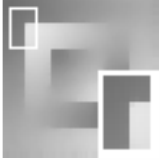}
		
		Iter = 3314 \hspace{0.25in} Iter  = 10000 \hspace{0.25in} Iter = 30000
		
		\includegraphics[width=1.2in,height=1.2in]{test1_star_alm-eps-converted-to.pdf}
		\includegraphics[width=1.2in,height=1.2in]{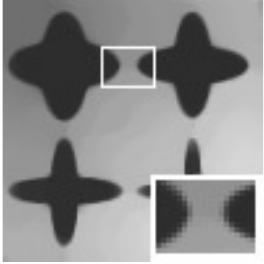}
		\includegraphics[width=1.2in,height=1.2in]{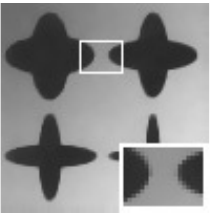}
		
		Iter = 1416 \hspace{0.25in} Iter  = 10000 \hspace{0.25in} Iter = 30000
		
		\includegraphics[width=1.2in,height=1.2in]{test1_lena_alm-eps-converted-to.pdf}
		\includegraphics[width=1.2in,height=1.2in]{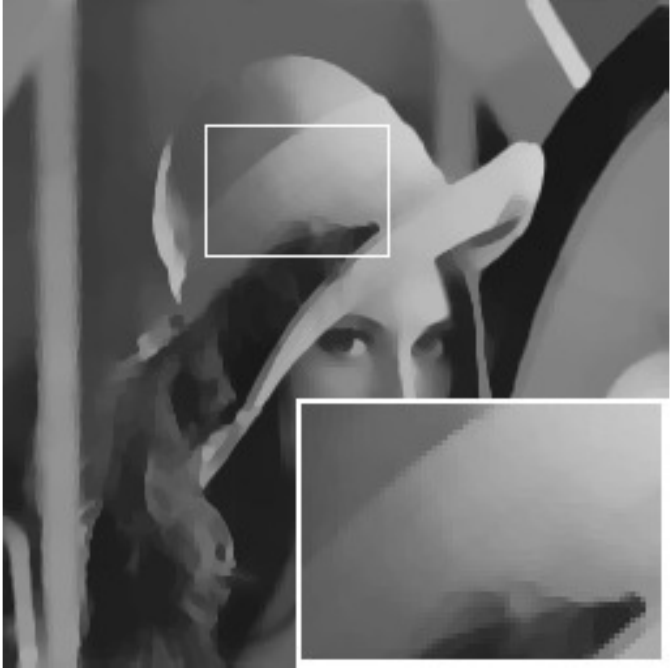}
		\includegraphics[width=1.2in,height=1.2in]{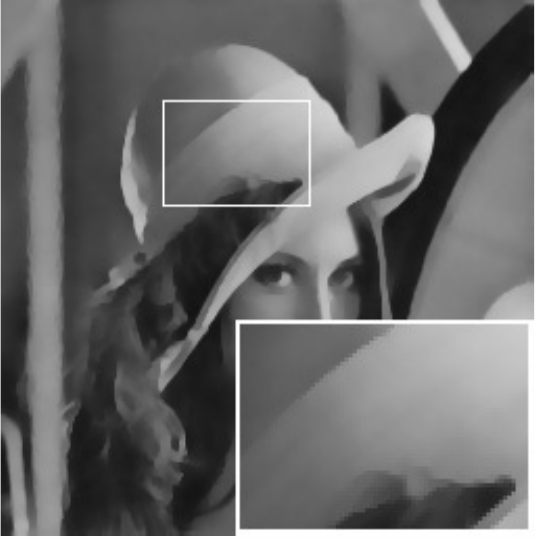}
		
		Iter = 3526 \hspace{0.25in} Iter  = 10000 \hspace{0.25in} Iter = 30000
		
	\end{center}
	\caption{This figure shows that one can improve the smoothing qualities of the THC method (for a non-optimal choice of the augmentation parameters $r_1$, $r_2$ and $r_3$) by requiring more iterations. We did the computations with $r_1 = 0.01$, $r_2 = 10$ and $r_3 = 100$, an augmentation parameter choice  which is optimal (or near optimal) for the ``ball'' image, but not for ``square'', ``star'' and ``Lena''.  Note that the left images are identical to the corresponding ones in the second column of Fig. \ref{fig:test1}. }
	\label{fig:test12}
\end{figure*}

\begin{figure*}
	\begin{center}
		\includegraphics[width=1.5in,height=1.5in]{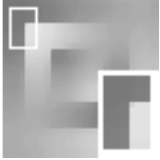}
		\includegraphics[width=1.5in,height=1.5in]{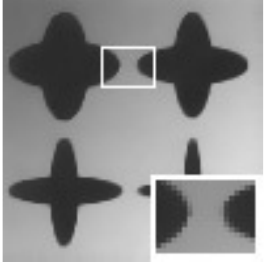}
		\includegraphics[width=1.5in,height=1.5in]{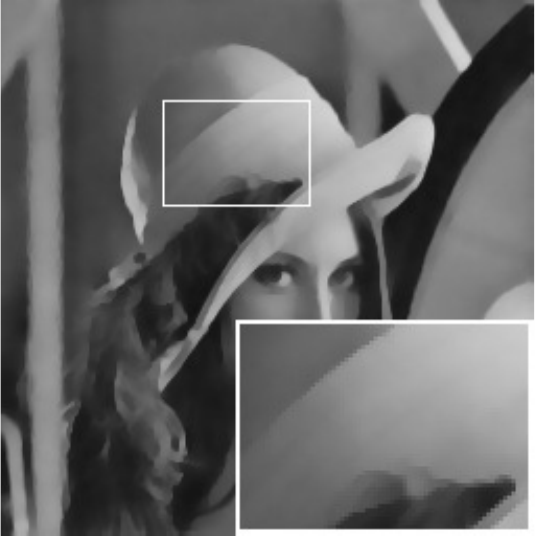}
		\includegraphics[width=1.5in,height=1.5in]{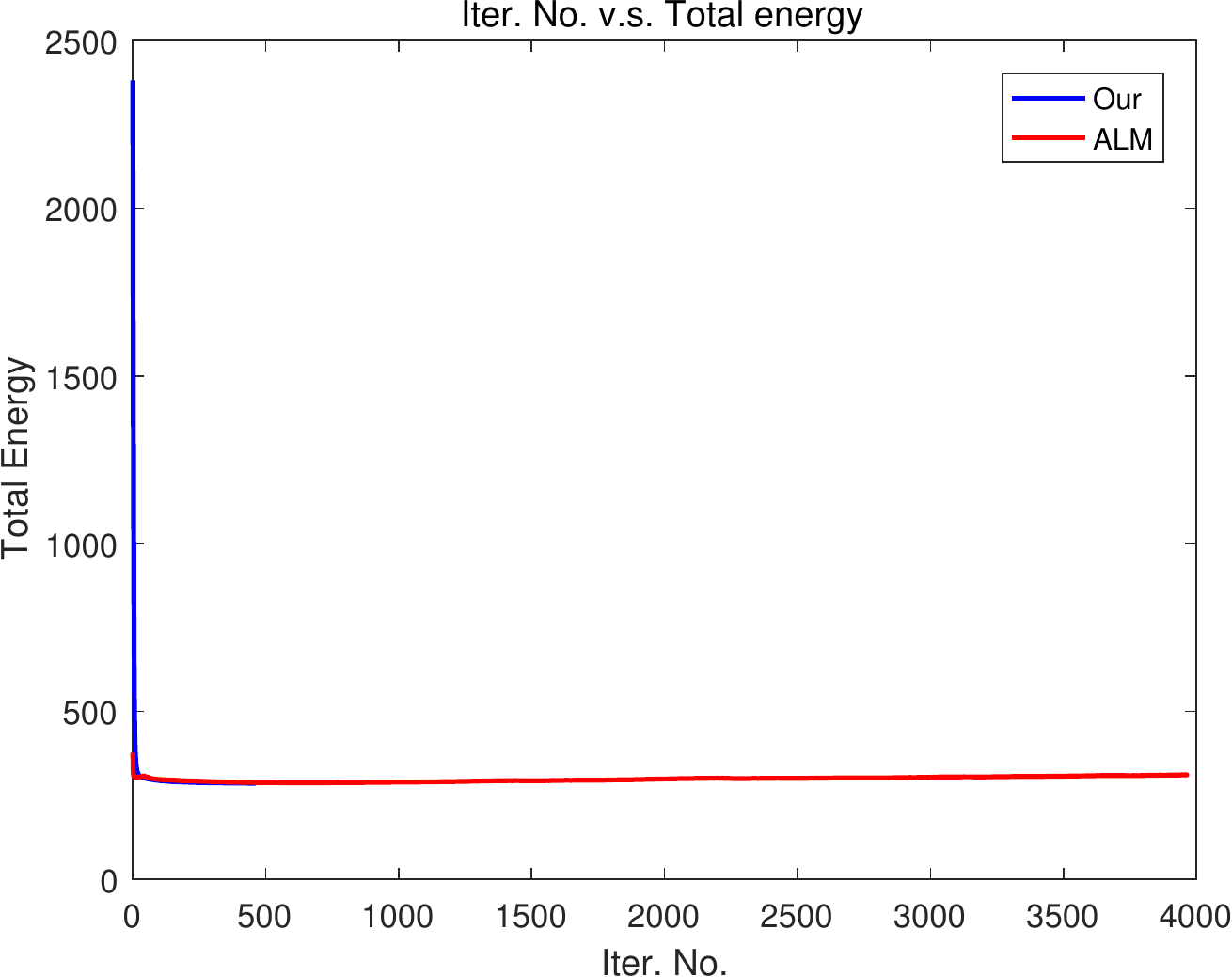}	
		
		(a) \hspace{1.3in} (b) \hspace{1.3in} (c) \hspace{1.3in} (d)
		
	\end{center}
	\caption{This figure shows that by a proper tuning of the augmentation parameters $r_1$, $r_2$ and $r_3$ one can significantly improve the smoothing properties of the THC method. Good choices are: (a)  $r_1 = 0.05$, $r_2 = 10$ and $r_3 = 100$ for ``square''. (b) $r_1 = 0.005$, $r_2 = 10$ and $r_3 = 100$ for ``star''. (c) $r_1 = 0.01$, $r_2 = 10$ and $r_3 = 300$ for ``Lena''. In Fig. \ref{fig:test2}(d) we have visualized for the ``Lena'' image, the variations of the elastica energy v.s. the iteration number for the proposed method (blue curve) and for the THC method (red curve): both methods reached the same limit. These results show that the THC method suffers from a strong image dependence concerning a good choice for the augmentation parameters $r_1$, $r_2$ and $r_3$.   }
	\label{fig:test2}
\end{figure*}

As already mentioned the THC method gives poor results when it reaches the stopping criteria given above. As expected, one can make the THC method operational again by either increasing the number of iterations (see Fig. \ref{fig:test12}) or modifying the augmentation parameters (see Fig. \ref{fig:test2}). We see in particular on Fig. \ref{fig:test2}(d) that, for the ``Lena'' image, the THC method with properly tuned augmentation parameters converges to a solution with the same energy than the one we obtain via the method proposed in this article.

\begin{figure*}[!t]
	\begin{center}
		\includegraphics[width=1.2in,height=1.2in]{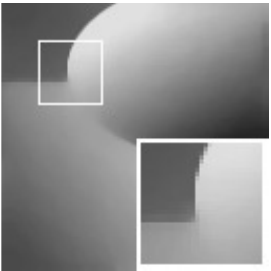}
		\includegraphics[width=1.2in,height=1.2in]{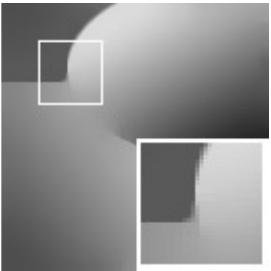}
		\includegraphics[width=1.2in,height=1.2in]{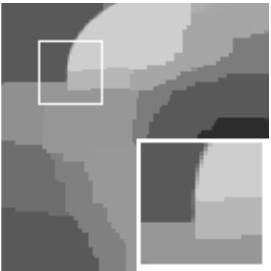}
		\includegraphics[width=1.2in,height=1.2in]{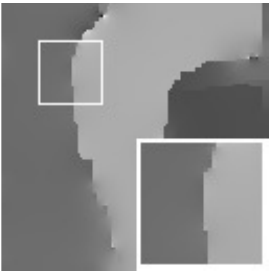}
		
		$r_1 = 0.001$ \hspace{0.6in} $r_1 = 0.01$ \hspace{0.6in} $r_1 = 0.1$ \hspace{0.6in} $r_1 = 1$
		
		\includegraphics[width=1.2in,height=1.2in]{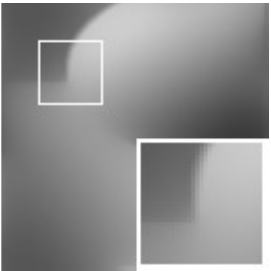}
		\includegraphics[width=1.2in,height=1.2in]{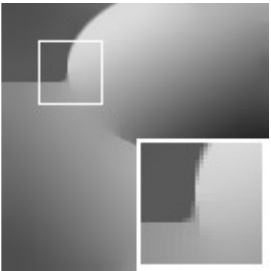}
		\includegraphics[width=1.2in,height=1.2in]{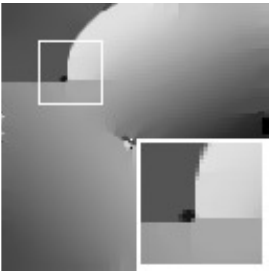}
		\includegraphics[width=1.2in,height=1.2in]{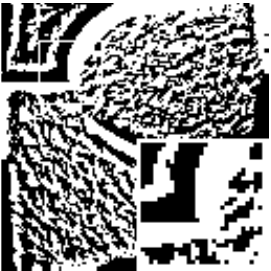}
		
		$r_2 = 1$ \hspace{0.5in} $r_2 = 10$ \hspace{0.5in} $r_2 = 20$ \hspace{0.5in} $r_2 = 30$
		
		\includegraphics[width=1.2in,height=1.2in]{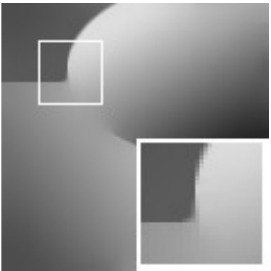}
		\includegraphics[width=1.2in,height=1.2in]{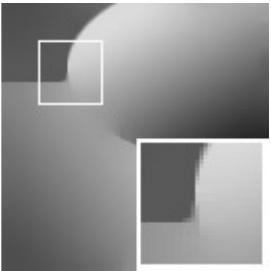}
		\includegraphics[width=1.2in,height=1.2in]{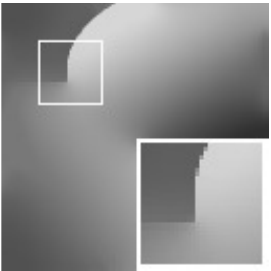}
		\includegraphics[width=1.2in,height=1.2in]{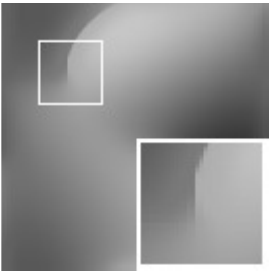}

		$r_3 = 50$ \hspace{0.5in} $r_3 = 100$ \hspace{0.5in} $r_3 = 500$ \hspace{0.5in} $r_3 = 1000$
		
	\end{center}
	\caption{This figure shows the dependence of the THC method computed solution to the augmentation parameters $r_1$, $r_2$ and $r_3$. The default choice being $r_1 = 0.01$, $r_2 = 10$ and $r_3 = 100$, one varied each time only one parameter leaving the other two unchanged (here $a = b = 0.1$ and $tol =1\times 10^{-5}$).   }
	\label{fig:test3}
\end{figure*}

%%%%%%%%%%%%%%%%%%%%%% TEST3 %%%%%%%%%%%%%%%%%%%%%%%%%%%%%%%%%%%%%%%%%%%%%%%%%%%
In Fig. \ref{fig:test3}, we reported the performances of the THC method for different values of $r_1$, $r_2$ and $r_3$. It is clear from this figure that the THC method is quite sensitive to the values of the augmentation parameters, $r_1$ and $r_2$ in particular, implying that augmentation parameter tuning is necessary for the THC method to have good convergence properties.

\subsubsection{Speed of convergence comparisons}\label{sec:4.3.2}
\begin{figure*}[!t]
	\begin{center}
		
		\includegraphics[width=6in,height=4in]{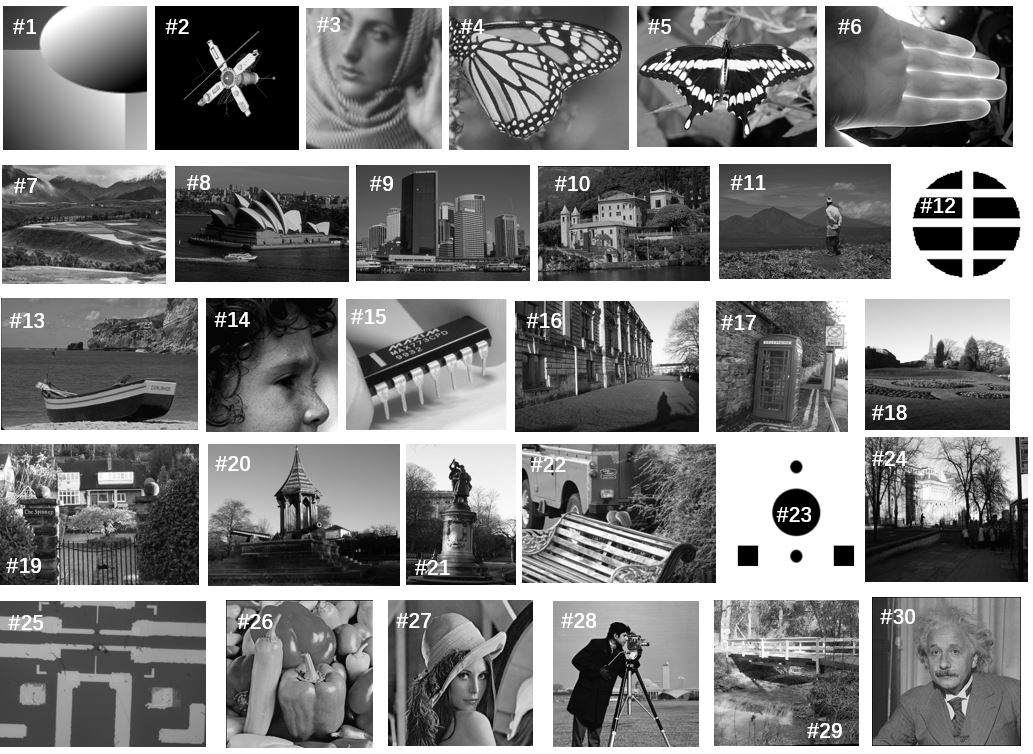}
		
	\end{center}
	\caption{Display of the 30 gray images we used to evaluate and compare the convergence properties of both the method we proposed in this article and the THC method. In order to facilitate the display, we managed to have same size images although they were originally of different sizes.  }
	\label{fig:test5_imgs}
\end{figure*}

We further compare in this subsection the speeds of convergence of the THC and proposed methods. To have fair comparisons, we collected 30 gray images (see Fig. \ref{fig:test5_imgs}), either synthetic or natural, and added Gaussian white noise with zero mean and various standard deviations (\textit{std}) to these images.

\begin{figure*}[!t]
	\begin{center}
		
		\includegraphics[width=2in,height=1.7in]{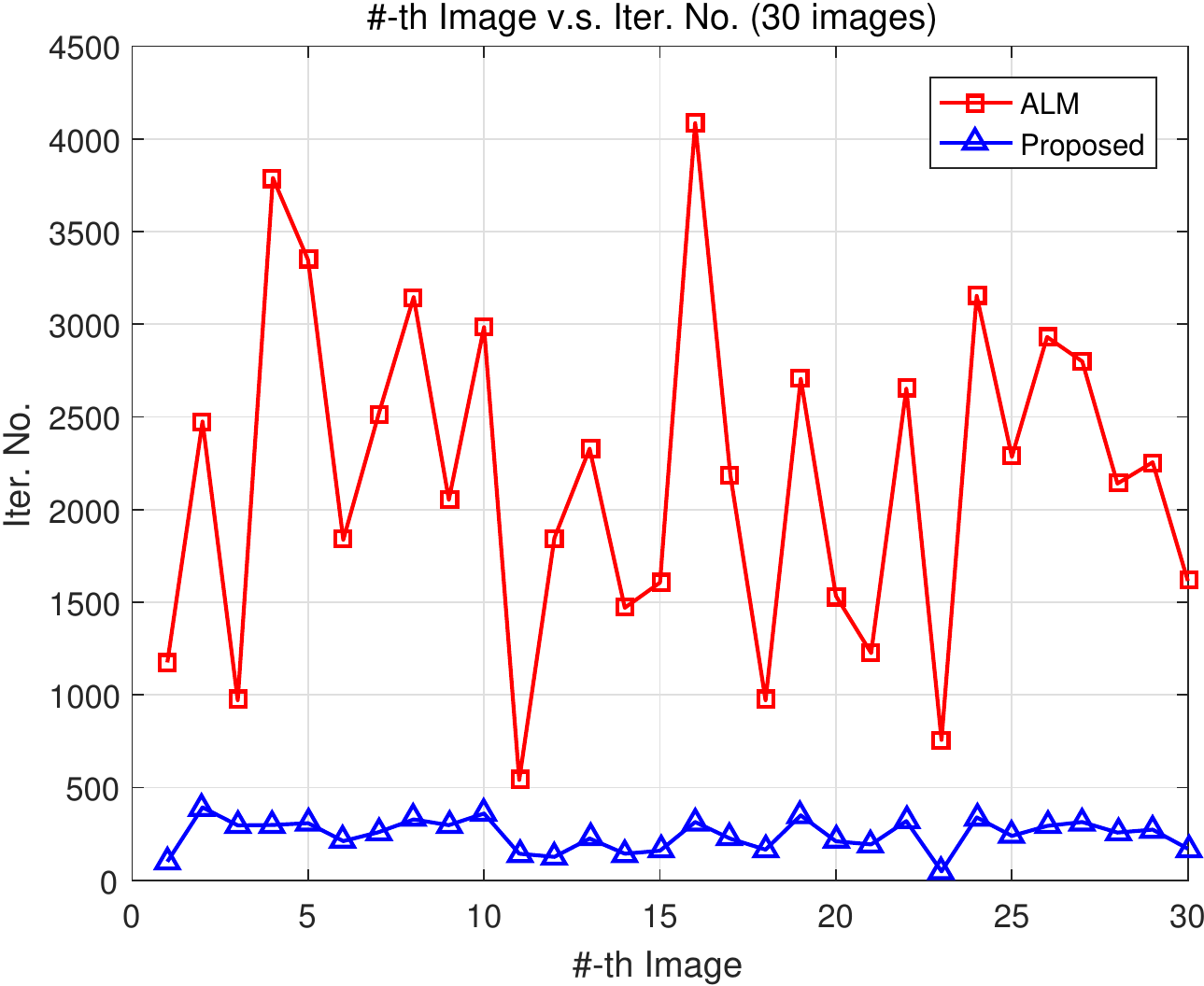}
		\includegraphics[width=2in,height=1.7in]{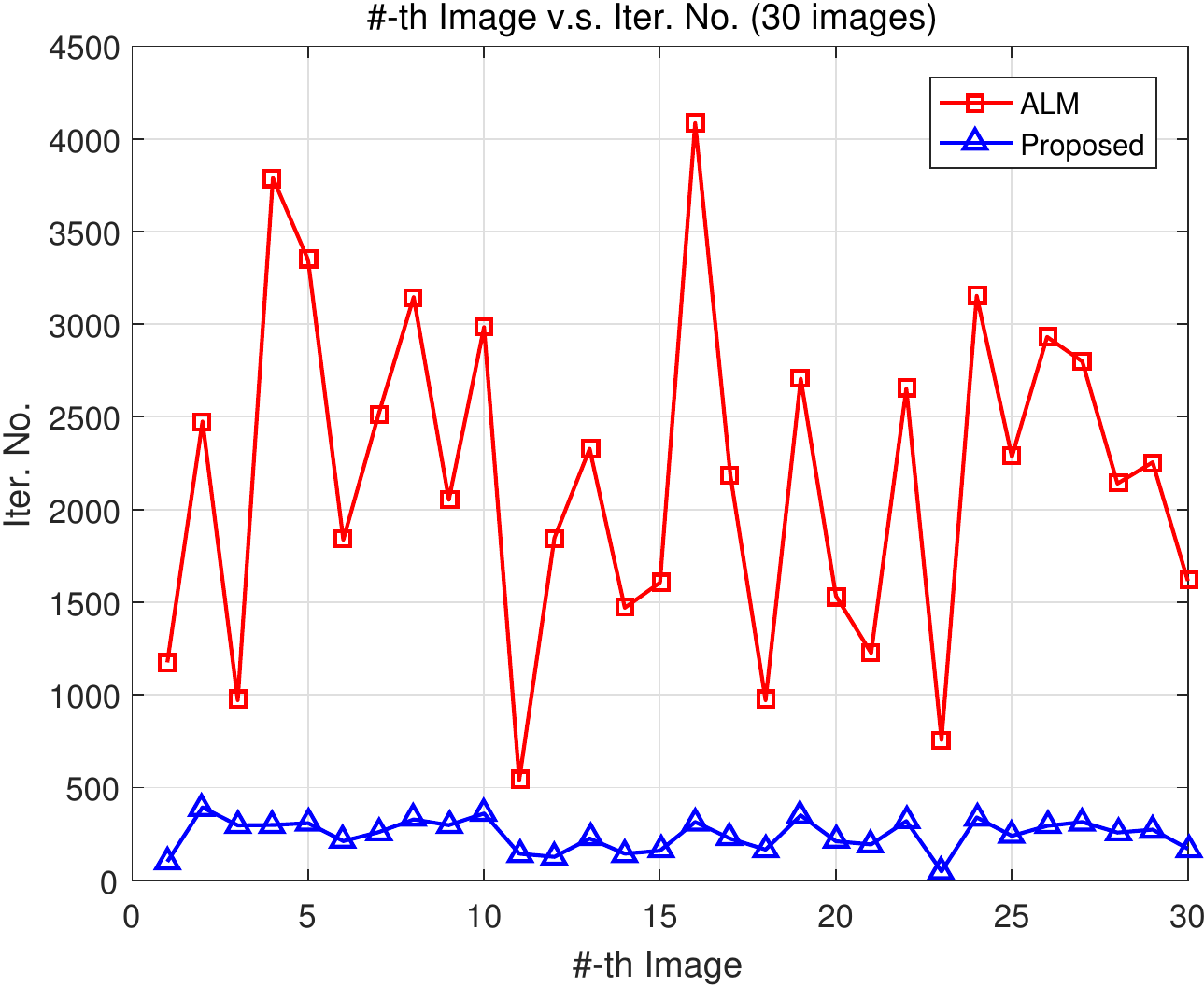}
		\includegraphics[width=2in,height=1.7in]{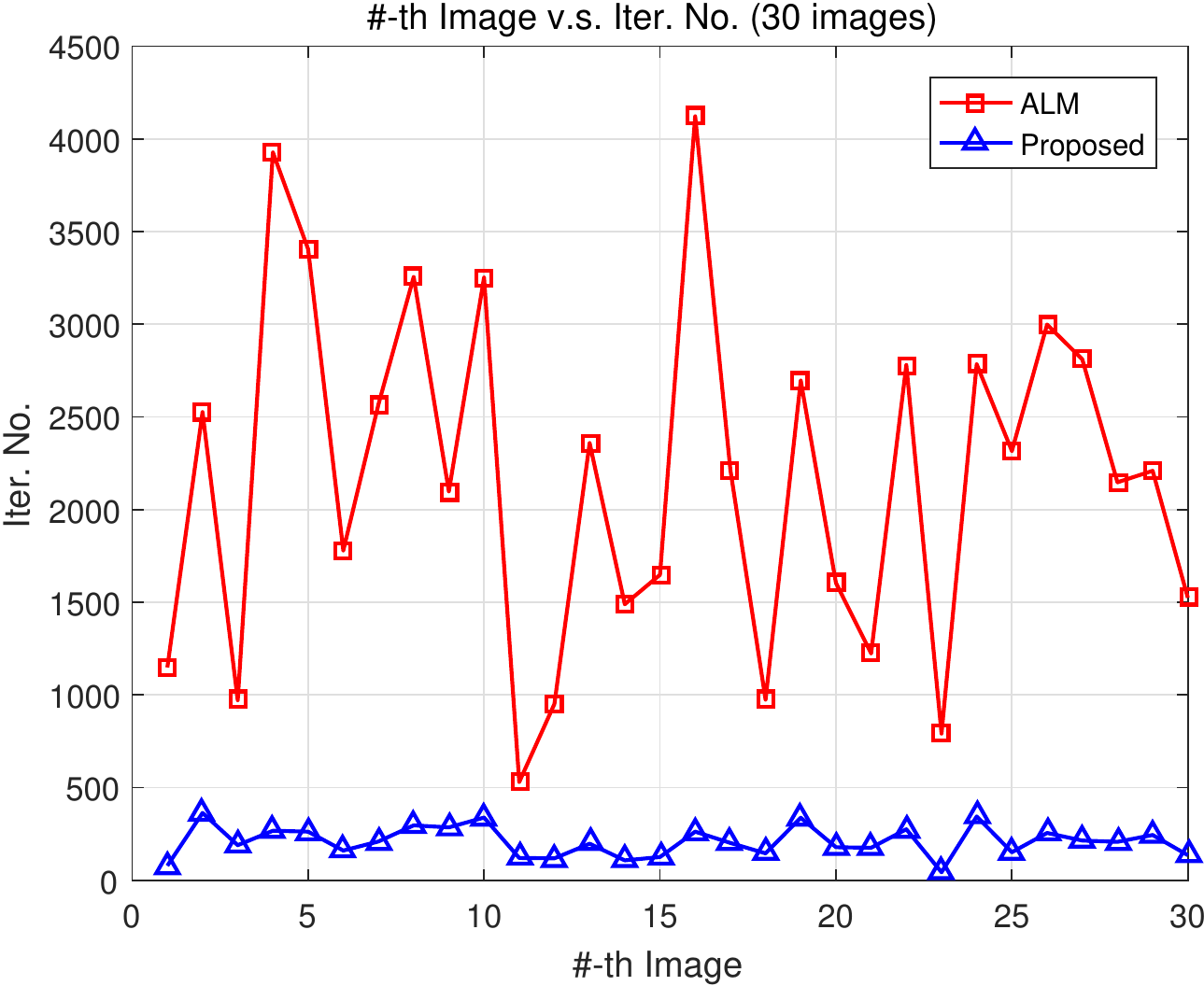}
		
		\includegraphics[width=2in,height=1.7in]{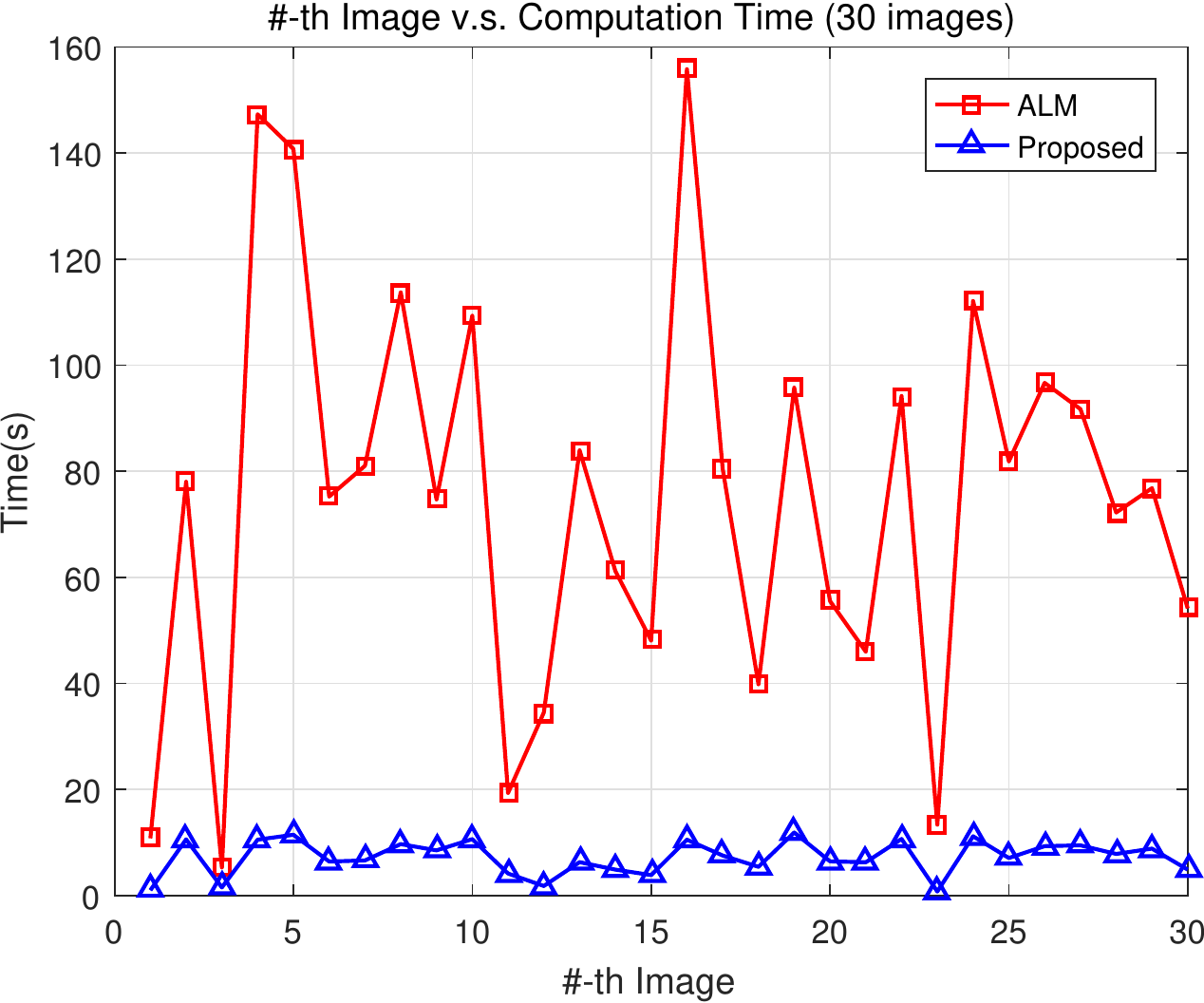}
		\includegraphics[width=2in,height=1.7in]{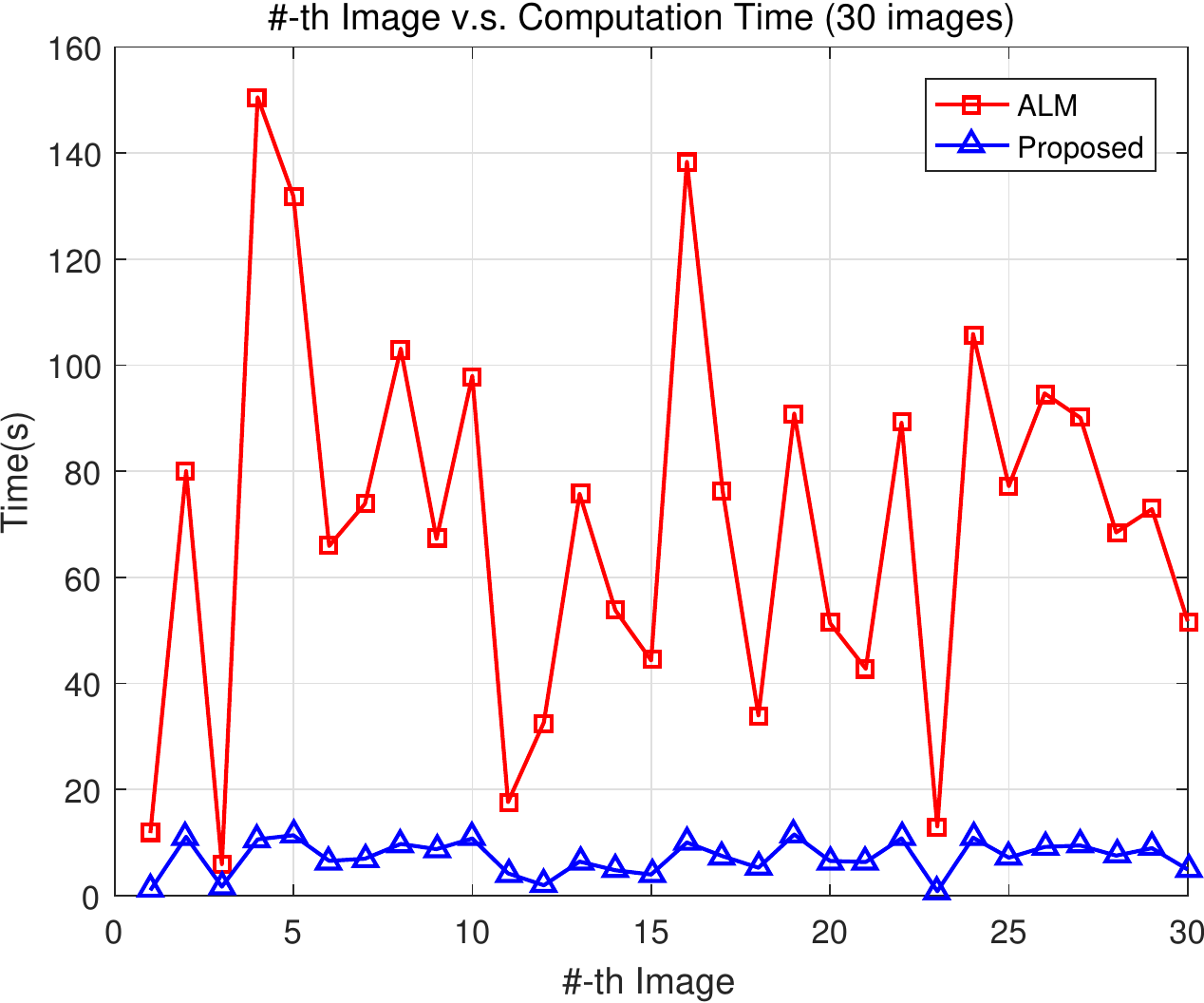}
		\includegraphics[width=2in,height=1.7in]{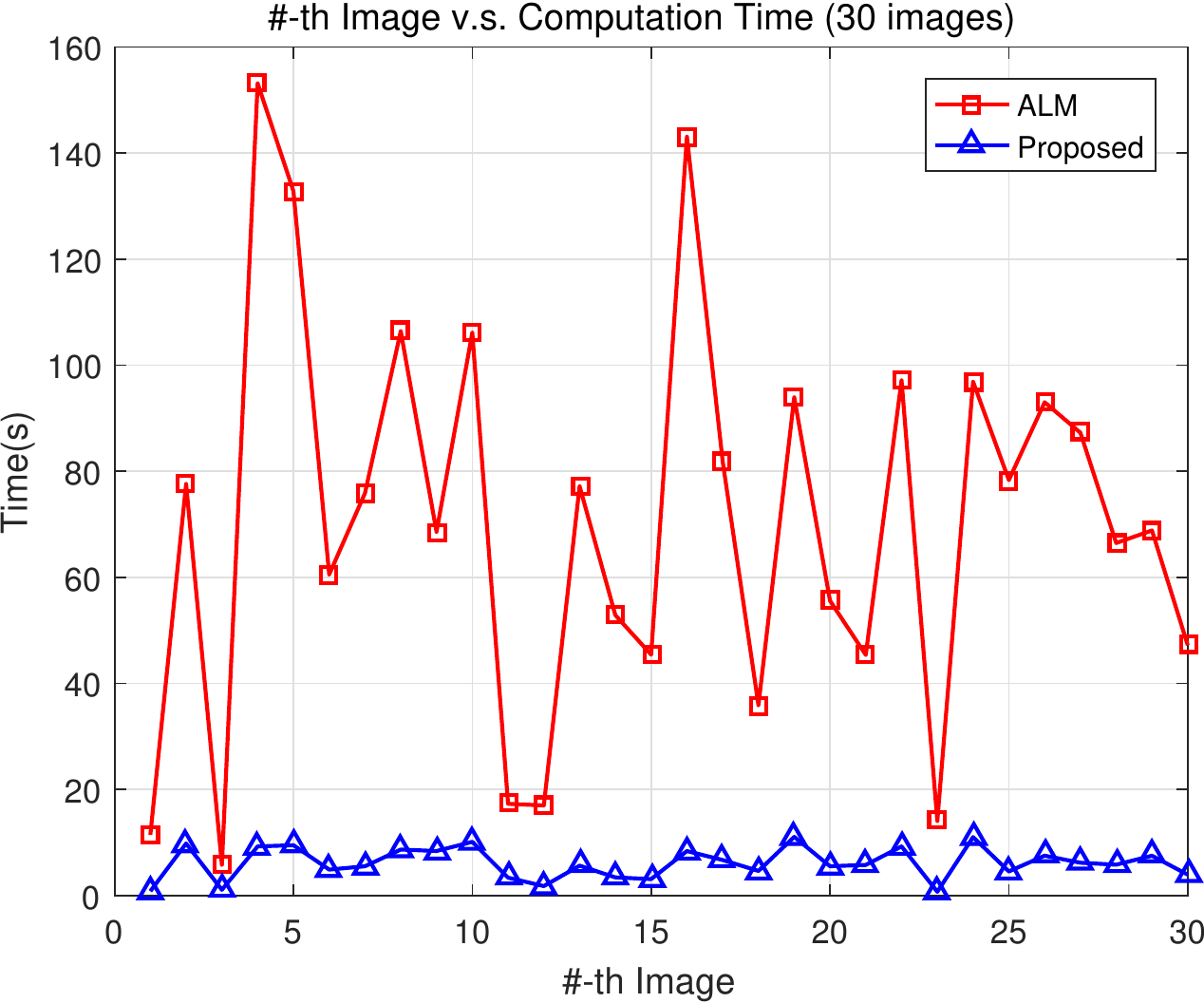}
		
		~~~~			$std = 0.1$ \hspace{1.2in} $std = 0.05$ \hspace{1.2in} $std = 0.02$
		
	\end{center}
	\caption{Number of iterations (first row) and corresponding computational time (second row) for the THC method (red curves) and the method introduced in this article (blue curves), both methods using $tol = 1\times 10^{-5}$ for their respective stopping criterion. The two methods are applied to the 30 images of Figure \ref{fig:test5_imgs}, with added zero mean Gaussian noise, and various standard deviations ($std$) ($std = 0.1$, $0.05$ and $0.02$).  }
	\label{fig:test5_tol-5}
\end{figure*}

\begin{figure*}[!t]
	\begin{center}
		
		\includegraphics[width=2in,height=1.7in]{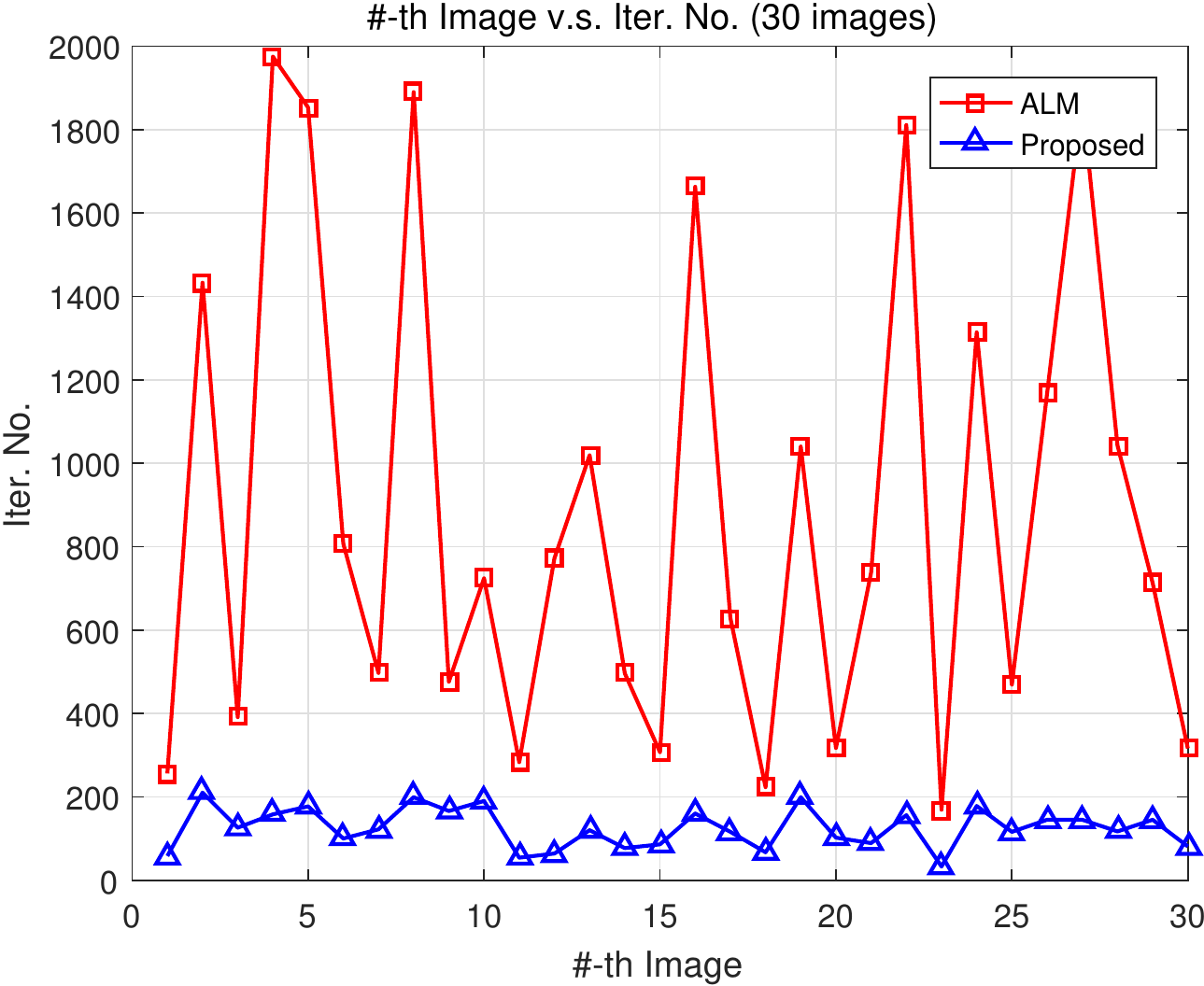}
		\includegraphics[width=2in,height=1.7in]{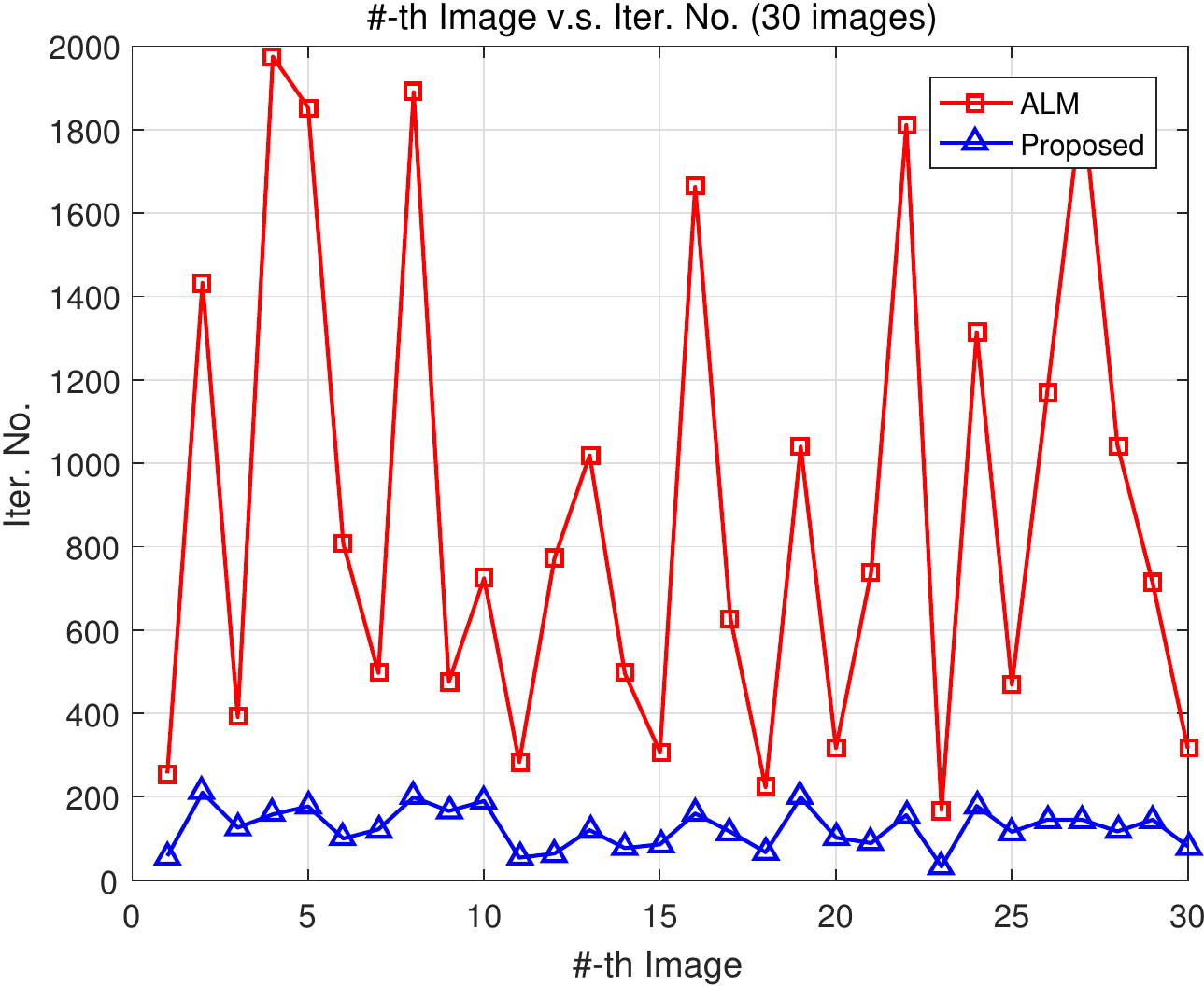}
		\includegraphics[width=2in,height=1.7in]{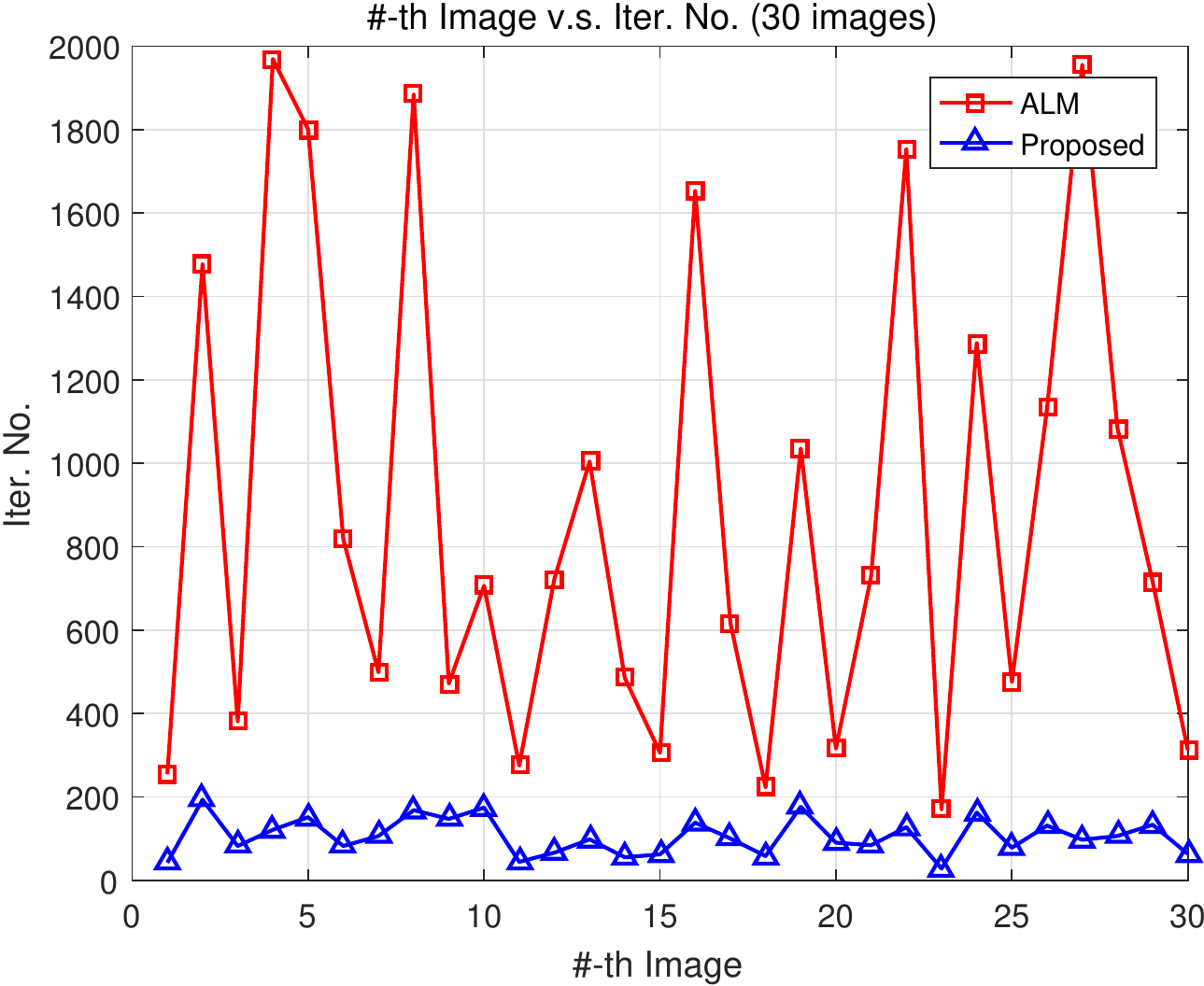}
		
		\includegraphics[width=2in,height=1.7in]{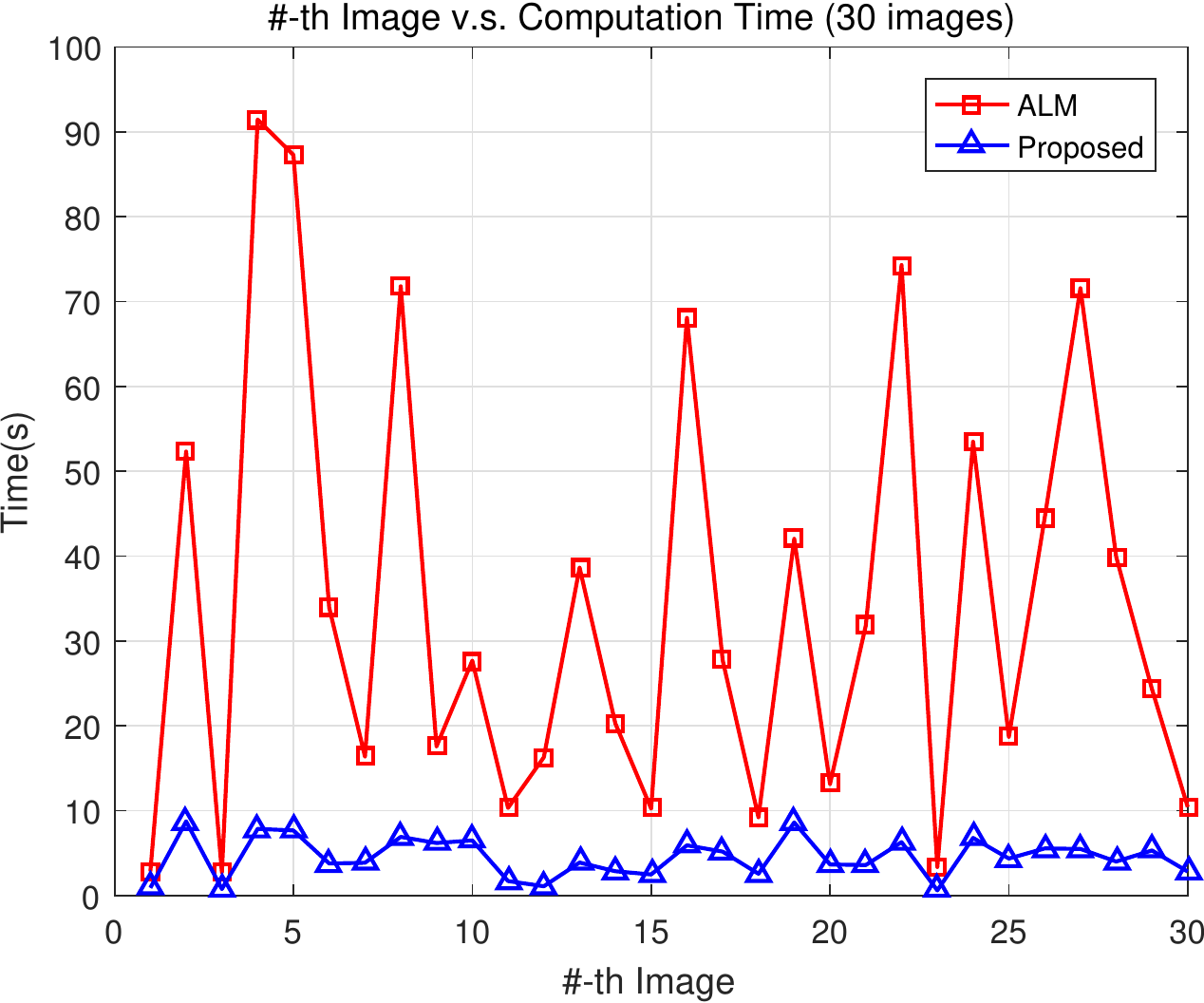}
		\includegraphics[width=2in,height=1.7in]{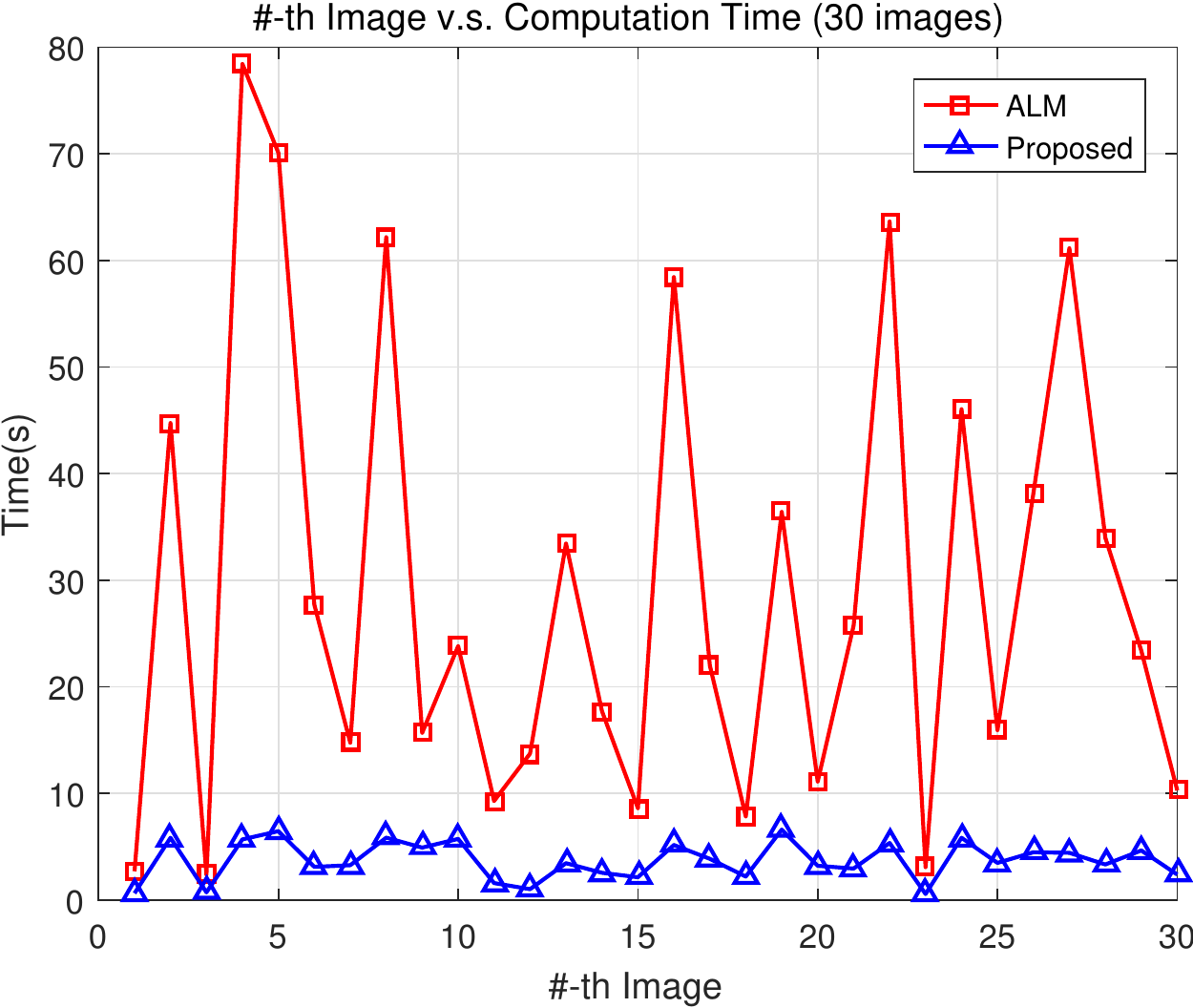}
		\includegraphics[width=2in,height=1.7in]{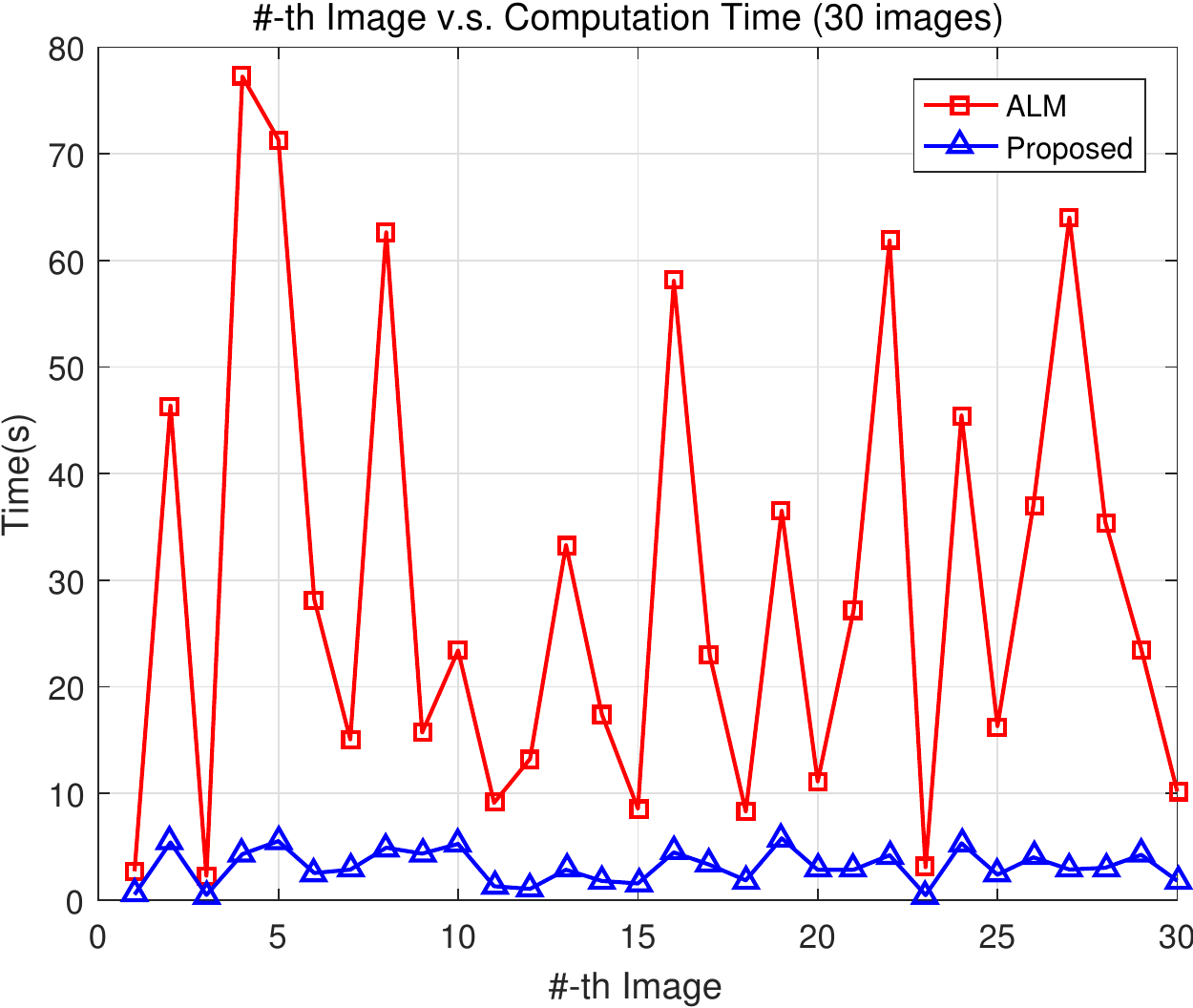}
		
		~~~~			$std = 0.1$ \hspace{1.2in} $std = 0.05$ \hspace{1.2in} $std = 0.02$
		
	\end{center}
	\caption{Number of iterations (first row) and corresponding computational time (second row) for the 30 images of Figure \ref{fig:test5_imgs}, using the THC method (red curves) and the method introduced in this article (blue curves).  Both methods use $tol = 1\times 10^{-5}$ for their respective stopping criterion.   }
	\label{fig:test5_tol5e-5}
\end{figure*}

\begin{figure*}
	\begin{center}
		
		\includegraphics[width=2in,height=1.7in]{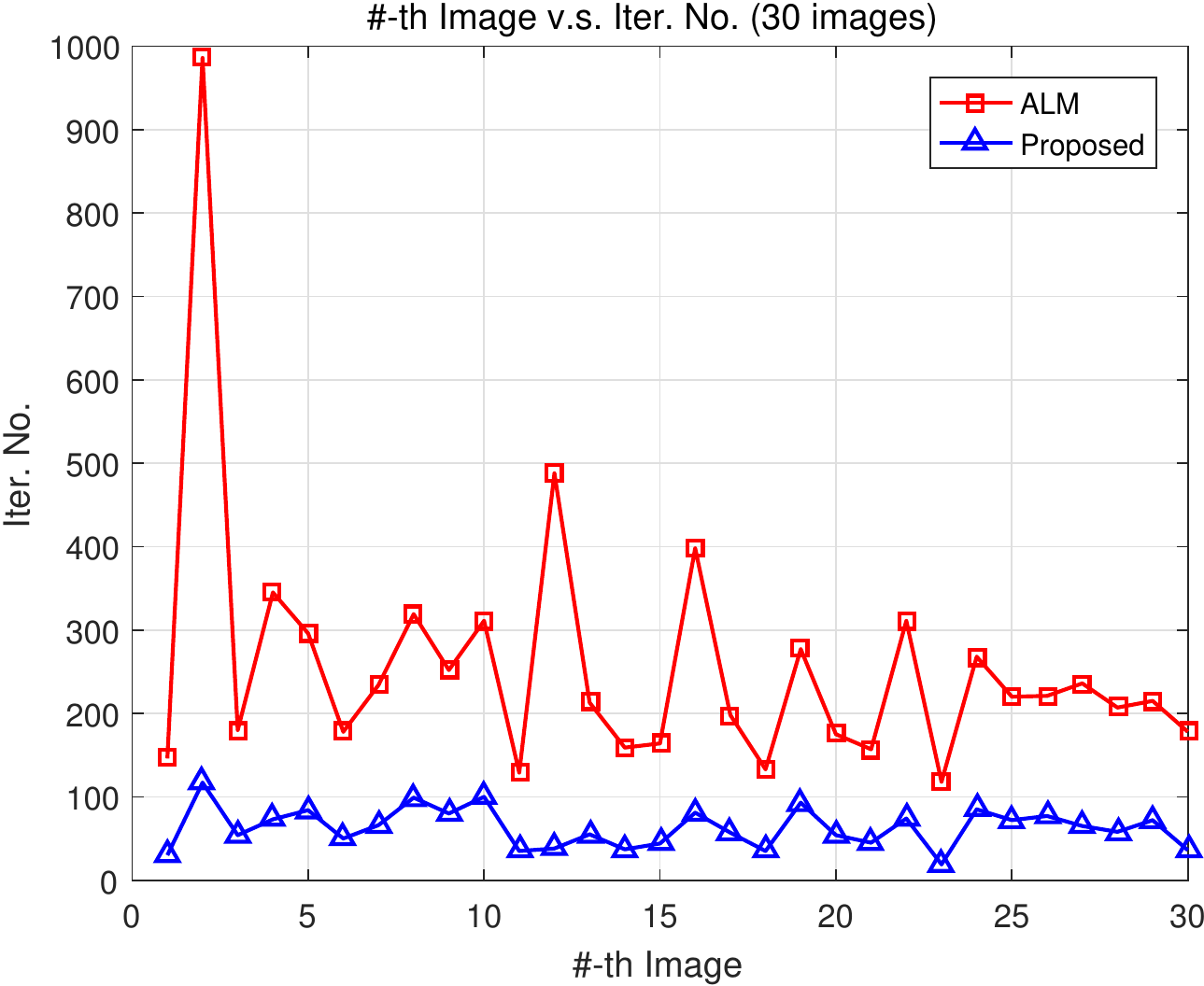}
		\includegraphics[width=2in,height=1.7in]{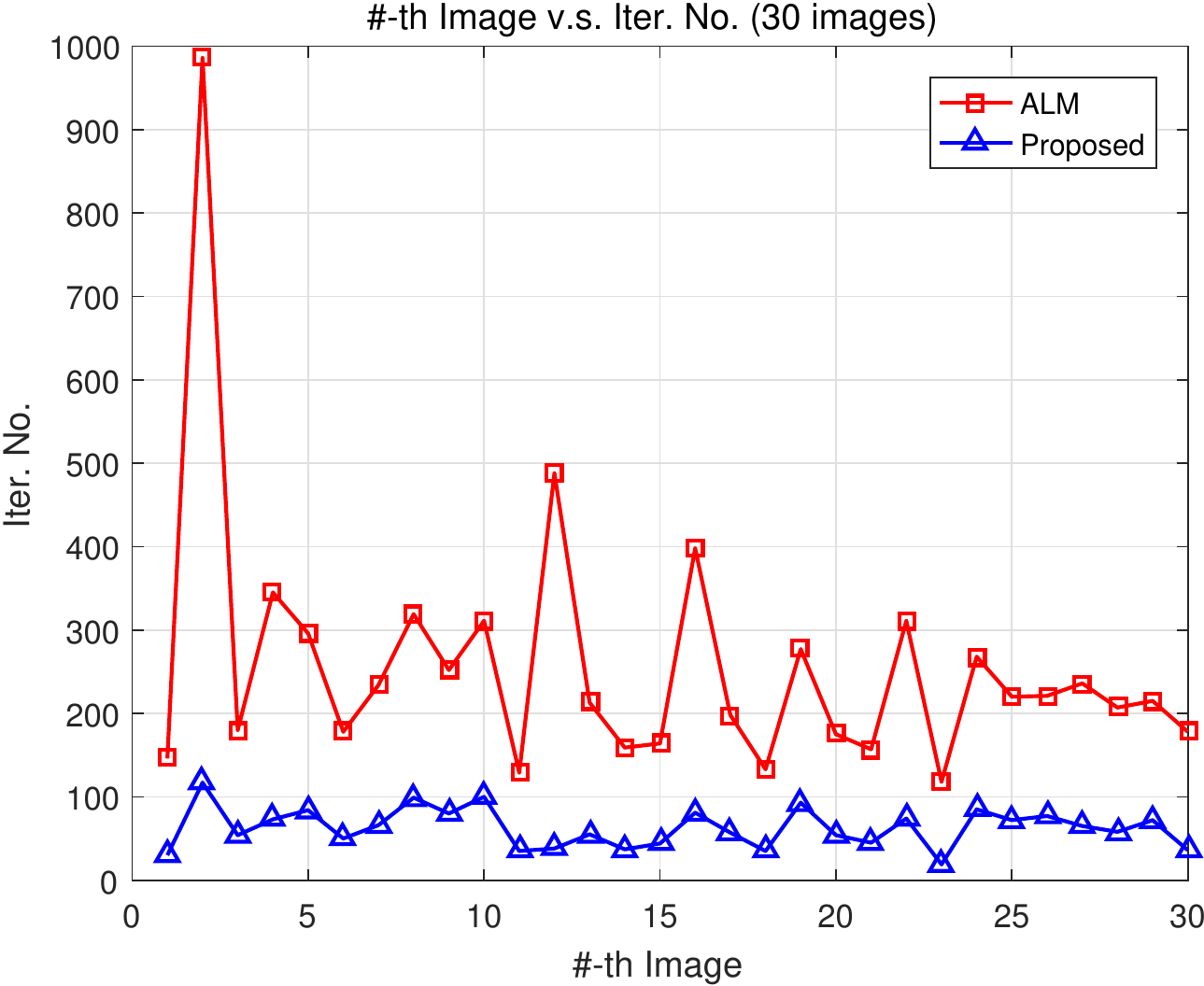}
		\includegraphics[width=2in,height=1.7in]{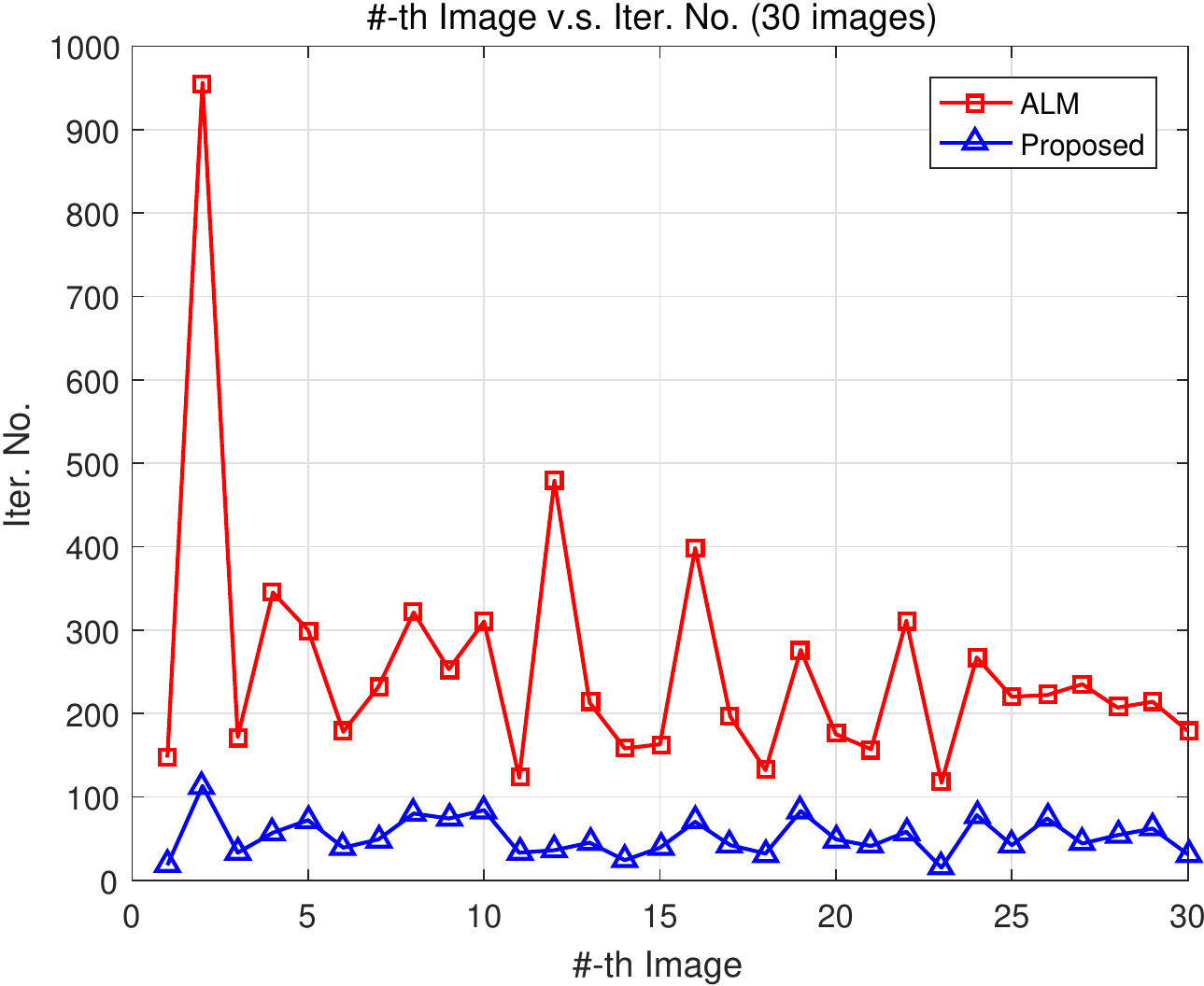}
		
		\includegraphics[width=2in,height=1.7in]{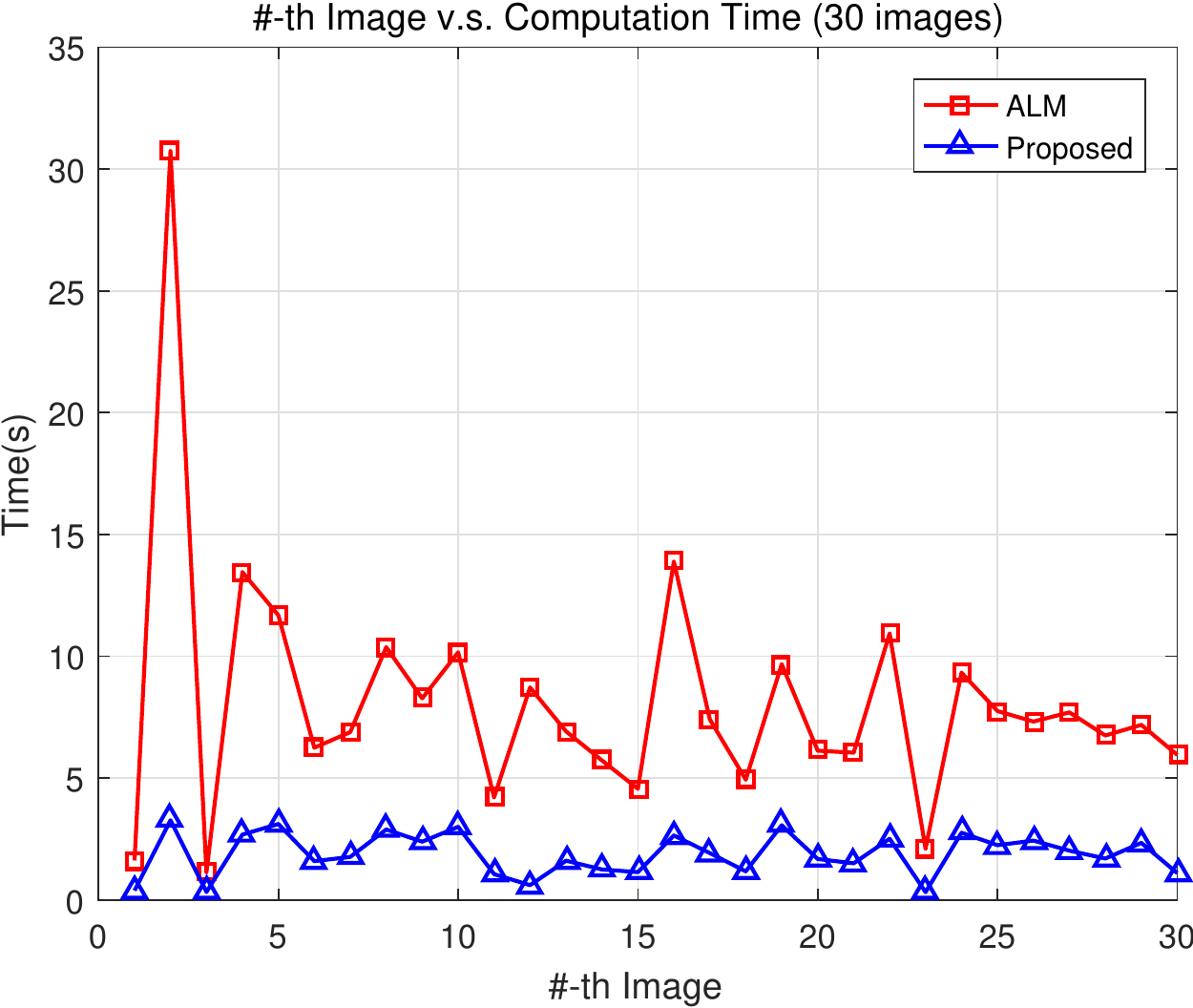}
		\includegraphics[width=2in,height=1.7in]{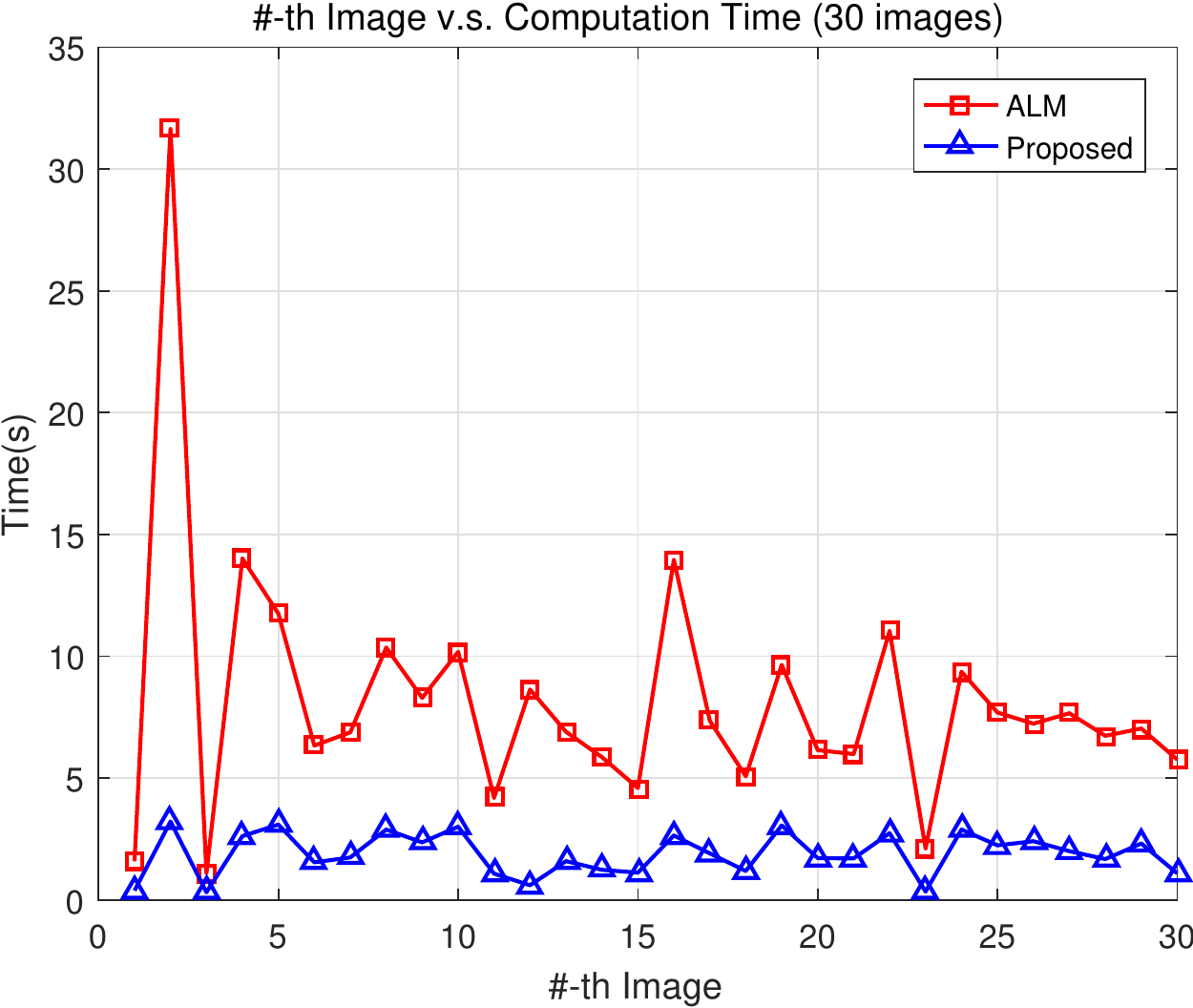}
		\includegraphics[width=2in,height=1.7in]{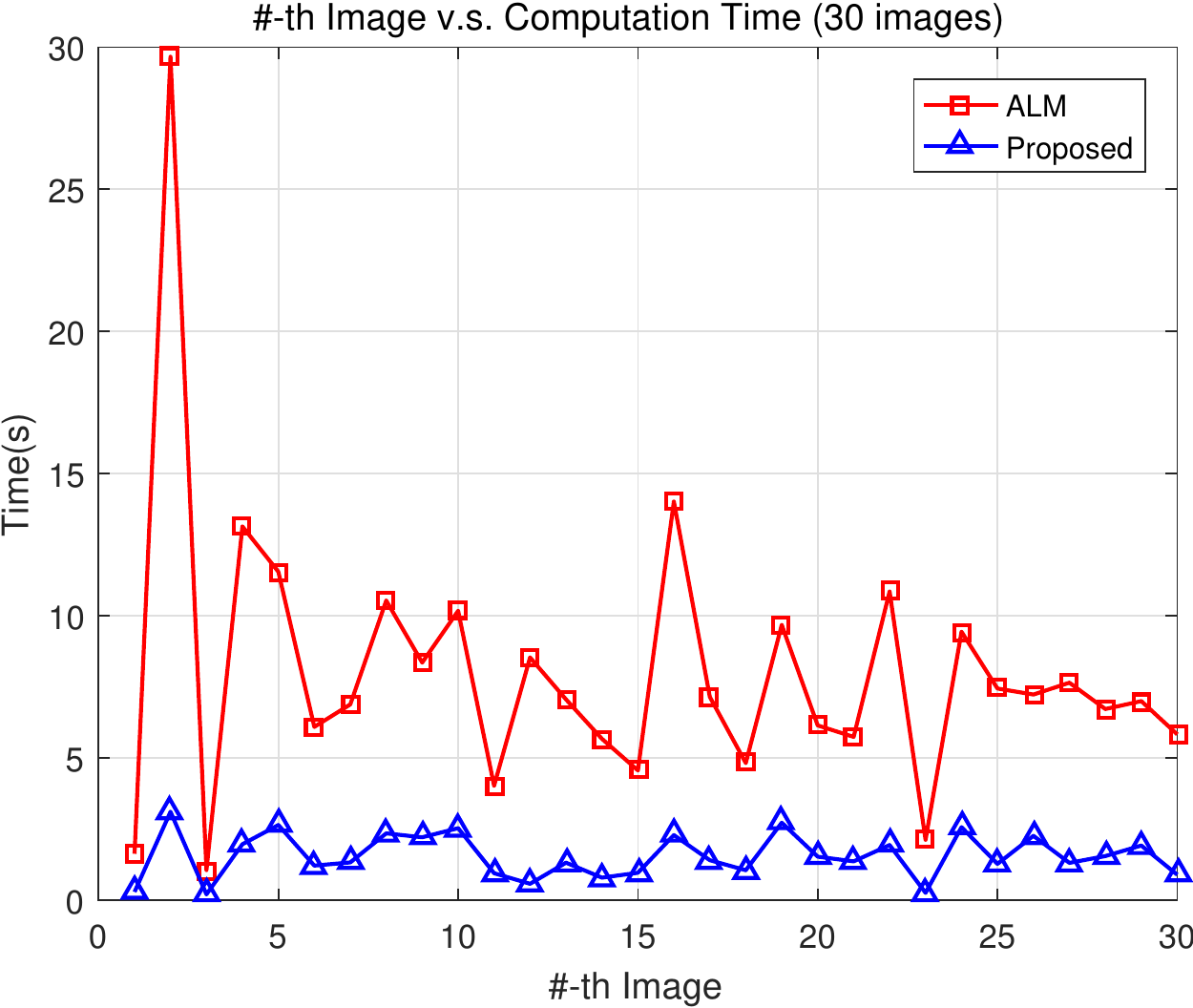}
		
		~~~~			$std = 0.1$ \hspace{1.2in} $std = 0.05$ \hspace{1.2in} $std = 0.02$
		
	\end{center}
	\caption{Number of iterations (first row) and corresponding computational time (second row) for the 30 images of Figure \ref{fig:test5_imgs}, using the THC method (red curves) and the method introduced in this article (blue curves).  Both methods use $tol = 1\times 10^{-4}$ for their respective stopping criterion.   }
	\label{fig:test5_tol-4}
\end{figure*}

On Figures \ref{fig:test5_tol-5}  to \ref{fig:test5_tol-4} we have reported for $tol = 1 \times 10^{-5}$ (Fig. \ref{fig:test5_tol-5}), $5\times 10^{-5}$ (Fig. \ref{fig:test5_tol5e-5}) and $1\times 10^{-4}$ (Fig. \ref{fig:test5_tol-4}), the averaged number of iterations needed to achieve convergence (first row) and the corresponding computational time (second row). These figures leave no doubt about the superiority of the method we introduced in this article over the augmented Lagrangian based THC method. Indeed, the new method outperforms THC's in terms of number of iterations and computational time per iteration (as shown by Tab. \ref{tab:test5}, which displays averaged performances), not to mention its greater simplicity and robustness.

\begin{table}[!t]
	\caption{Averaged computational time (in seconds and per image) of the method we propose in this article and of the THC method, the images being the 30 gray images displayed in Fig. \ref{fig:test5_imgs}.}\label{tab:test5}
	\scriptsize
	\begin{center}
		\begin{tabular}{|c|c|c|c|c|}
			\hline
			$tol$ & Method & $std = 0.1$  & $std = 0.05$  &  $std = 0.02$        \\
			\hline\hline
			\multirow{2}{*}{$tol = 1\times 10^{-5}$}
			&{\textbf{Proposed}}&{7.2}&{7.2}&{6.1}\\
			
			\cline{2-5}
			& \textbf{THC} &  74.9 & 70.3  &  70.4  \\
			
			\hline
			\multirow{2}{*}{$tol = 5\times 10^{-5}$} &{\textbf{Proposed}}&{4.5}&{3.7}&{3.1}\\
			
			\cline{2-5}
			& \textbf{THC} &  34.4 & 29.4  &  29.6   \\
			
			\hline
			\multirow{2}{*}{$tol = 1\times 10^{-4}$}
			&{\textbf{Proposed}}&{1.9}&{1.9}&{1.6}\\
			
			\cline{2-5}
			& \textbf{THC} &  8.1 & 8.2  &  8.0   \\

			\hline
		\end{tabular}
	\end{center}
\end{table}

\section{Conclusions}\label{sec:conclusion}
In the paper, we proposed a simple and efficient operator splitting approach to solve the Euler elastica model, and applied the proposed method to image smoothing. Different from the ALM method, the proposed method only needs to tume one parameter, i.e. the time step.  Numerical experiments demonstrated that the proposed method works well for the Euler elastica energy and produces good results for image smoothing. Moreover, extensive test experiments were also designed and implemented to assess the stability and effectiveness of the proposed method. Furthermore, the comparisons with the THC method demonstrated also that the proposed method is fast, stable and robust.

%%%%%%%%%%%%%%%%%%%%%%%%%%%%%%%%%%%%%%
\section*{Acknowledgments}
The first author would like to thank the supports by NSFC (61702083, 61772003, 61876203). The third author would like to acknowledge the support from HKBU startup grant,  RG(R)- RC/17-18/02-MATH, and FRG2/17-18/033.

%-----------refs----------------%

%%--------------
%\bibliographystyle{IEEEbib}
%%\bibliography{refs}
%\bibliography{references}

\end{document}